# 2004 IMS MEDALLION LECTURE

## LOCAL RADEMACHER COMPLEXITIES AND ORACLE INEQUALITIES IN RISK MINIMIZATION[1,2]

By Vladimir Koltchinskii

*University of New Mexico and Georgia Institute of Technology*

Let $\mathcal{F}$ be a class of measurable functions $f : S \mapsto [0,1]$ defined on a probability space $(S, \mathcal{A}, P)$. Given a sample $(X_1, \ldots, X_n)$ of i.i.d. random variables taking values in $S$ with common distribution $P$, let $P_n$ denote the empirical measure based on $(X_1, \ldots, X_n)$. We study an empirical risk minimization problem $P_n f \to \min$, $f \in \mathcal{F}$. Given a solution $\hat{f}_n$ of this problem, the goal is to obtain very general upper bounds on its excess risk

$$\mathcal{E}_P(\hat{f}_n) := P\hat{f}_n - \inf_{f \in \mathcal{F}} Pf,$$

expressed in terms of relevant geometric parameters of the class $\mathcal{F}$. Using concentration inequalities and other empirical processes tools, we obtain both distribution-dependent and data-dependent upper bounds on the excess risk that are of asymptotically correct order in many examples. The bounds involve localized sup-norms of empirical and Rademacher processes indexed by functions from the class. We use these bounds to develop model selection techniques in abstract risk minimization problems that can be applied to more specialized frameworks of regression and classification.

**1. Introduction.** Let $(S, \mathcal{A}, P)$ be a probability space and let $\mathcal{F}$ be a class of measurable functions $f : S \mapsto [0,1]$. Let $(X_1, \ldots, X_n)$ be a sample of i.i.d. random variables defined on a probability space $(\Omega, \Sigma, \mathbb{P})$ and taking values

Received September 2003; revised July 2005.

[1]Supported in part by NSF Grant DMS-03-04861 and NSA Grant MDA904-02-1-0075.

[2]Discussed in 10.1214/009053606000001028, 10.1214/009053606000001037, 10.1214/009053606000001046, 10.1214/009053606000001055, 10.1214/009053606000001064 and 10.1214/009053606000001073; rejoinder at 10.1214/009053606000001082.

in $S$ with common distribution $P$. Let $P_n$ denote the empirical measure based on the sample $(X_1, \ldots, X_n)$.

We consider the problem of risk minimization

$$Pf \to \min, \qquad f \in \mathcal{F} \tag{1.1}$$

under the assumption that the distribution $P$ is unknown and has to be replaced by its estimate, $P_n$. Thus, the true risk minimization is replaced by the empirical risk minimization:

$$P_n f \to \min, \qquad f \in \mathcal{F}. \tag{1.2}$$

DEFINITION. Let

$$\mathcal{E}(f) := \mathcal{E}_P(f) := \mathcal{E}_P(\mathcal{F}; f) := Pf - \inf_{g \in \mathcal{F}} Pg.$$

This quantity will be called the *excess risk* of $f \in \mathcal{F}$. The set $\mathcal{F}_P(\delta) := \{f \in \mathcal{F} : \mathcal{E}_P(f) \leq \delta\}$ will be called the $\delta$-minimal set of $P$. In particular, $\mathcal{F}_P(0)$ is the minimal set of $P$.

Given a solution (or an approximate solution) $\hat{f} = \hat{f}_n$ of (1.2), *the first problem* of interest is to provide very general upper confidence bounds on the excess risk $\mathcal{E}_P(\hat{f}_n)$ of $\hat{f}_n$ that take into account some relevant geometric parameters of the class $\mathcal{F}$ as well as some measures of accuracy of approximation of $P$ by $P_n$ locally in the class. Namely, based on the $L_2(P)$-diameter $D_P(\mathcal{F}; \delta)$ of the $\delta$-minimal set $\mathcal{F}(\delta)$ and the function

$$\phi_n(\mathcal{F}; \delta) := \mathbb{E} \sup_{f,g \in \mathcal{F}(\delta)} |(P_n - P)(f - g)|,$$

we construct a quantity $\bar{\delta}_n(\mathcal{F}; t)$ such that inequalities of the following type hold:

$$\mathbb{P}\{\mathcal{E}_P(\hat{f}_n) \geq \bar{\delta}_n(\mathcal{F}; t)\} \leq \log \frac{n}{t} e^{-t}, \qquad t > 0$$

(see Section 3). The bound $\bar{\delta}_n(\mathcal{F}; t)$ has an asymptotically correct order (with respect to $n$) in many particular examples of risk minimization problems occurring in regression, classification and machine learning. However, if the diameter $D_P(\mathcal{F}; \delta)$ does not tend to $0$ as $\delta \to 0$ (which is the case when the risk minimization problem has multiple solutions), it happens that the bound $\bar{\delta}_n(\mathcal{F}; t)$ is no longer tight, and one has to redefine it using more subtle characteristics of geometry of the class than $D_P(\mathcal{F}; \delta)$ (see Section 4).

We will now describe a heuristic way to derive such bounds. It is based on *iterative localization* of the bound and it can be made precise (see the remark after the proof of Theorem 2 in Section 9 and also [27] where this type of argument was introduced in a more specialized setting). Define

$$\bar{U}_n(\delta;t) := K\left(\phi_n(\mathcal{F};\delta) + D(\mathcal{F};\delta)\sqrt{\frac{t}{n}} + \frac{t}{n}\right).$$

It follows from Talagrand's concentration inequality (see Section 2.1) that with some constant $K > 0$ for all $t > 0$

$$\mathbb{P}\left\{\sup_{f,g \in \mathcal{F}(\delta)} |(P_n - P)(f - g)| \geq \bar{U}_n(\delta;t)\right\} \leq e^{-t}.$$

Take $\delta^{(0)} = 1$, so that $\mathcal{F}(\delta^{(0)}) = \mathcal{F}$ (recall that functions in $\mathcal{F}$ take values in $[0,1]$). Assume, for simplicity, that the minimum of $Pf$ is attained at $\bar{f} \in \mathcal{F}$. Since $\hat{f}, \bar{f} \in \mathcal{F}(\delta^{(0)})$ and $P_n\hat{f} \leq P_n\bar{f}$, we have with probability at least $1 - e^{-t}$

$$\mathcal{E}_P(\hat{f}) = P\hat{f} - P\bar{f} = P_n\hat{f} - P_n\bar{f} + (P - P_n)(\hat{f} - \bar{f})$$
$$\leq \sup_{f,g \in \mathcal{F}(\delta)} |(P_n - P)(f - g)| \leq \bar{U}_n(\delta^{(0)};t) \wedge 1 =: \delta^{(1)}.$$

This implies that $\hat{f}, \bar{f} \in \mathcal{F}(\delta^{(1)})$ and we can repeat the above argument to show that with probability at least $1 - 2e^{-t}$, $\mathcal{E}_P(\hat{f}) \leq \bar{U}_n(\delta^{(1)};t) \wedge 1 =: \delta^{(2)}$. Iterating the argument $N$ times shows that with probability at least $1 - Ne^{-t}$ we have $\mathcal{E}_P(\hat{f}) \leq \delta^{(N)}$, where $\delta^{(N)} := \bar{U}_n(\delta^{(N-1)};t) \wedge 1$. If the sequence $\delta^{(N)}$ converges to the solution $\bar{\delta}$ of the fixed point equation $\delta = \bar{U}_n(\delta;t) \wedge 1$ and if the convergence is fast enough so that with some $C > 1$ for relatively small $N$ we have $\delta^{(N)} \leq C\bar{\delta}$, the above argument shows that $\mathcal{E}_P(\hat{f}) \leq C\bar{\delta}$ with probability at least $1 - Ne^{-t}$. Both with and without this iterative argument, we show in Section 3 (and prove in Section 9) that the construction of good upper bounds on the excess risk of $\hat{f}$ is related to *fixed point-type equations* for function $\bar{U}_n(\delta;t)$. The fixed point method has been developed in recent years in Massart [36], Koltchinskii and Panchenko [27] and Bartlett, Bousquet and Mendelson [5] (and in several other papers of these authors).

The second problem is to develop *ratio-type inequalities* for the excess risk, namely, to bound the following probabilities:

$$\mathbb{P}\left\{\sup_{f \in \mathcal{F}, \mathcal{E}_P(f) \geq \delta} \left|\frac{\mathcal{E}_{P_n}(f)}{\mathcal{E}_P(f)} - 1\right| \geq \varepsilon\right\}$$

(see Section 3). This problem is an important ingredient of the analysis of empirical risk minimization [in particular, we will use inequalities for such probabilities in our construction of data-dependent bounds on the excess risk $\mathcal{E}_P(\hat{f})$] and it is related to the study of ratio-type empirical processes (see [19, 20] for recent results on this subject).

The third problem is to construct data-dependent upper confidence bounds on $\mathcal{E}_P(\hat{f}_n)$. To this end, we replace the geometric parameters of the class [such as $D_P(\mathcal{F};\delta)$] by their empirical versions and the empirical process involved in the definition of data-dependent bounds by the Rademacher process (Section 3). The idea to use sup-norms or localized sup-norms of the Rademacher process as bootstrap-type estimates of the size of corresponding suprema of the empirical process has originated in machine learning literature (see [4, 5, 14, 26, 27, 34]). The current paper continues this line of research. Very recently, Bartlett and Mendelson [7] developed an interesting new definition of localized Rademacher complexities and gave a curious example in which this complexity provides a sharper bound on the risk of empirical risk minimizers than the complexities studied so far. It is not clear yet whether the phenomenon they studied occurs in actual machine learning or statistical problems. Because of this, we do not pursue this approach in the current paper.

The fourth problem is to develop rather general model selection techniques in risk minimization that utilize our data-dependent bounds on the excess risk (Sections 5, 6). More precisely, we study a version of structural risk minimization in which the class $\mathcal{F}$ is approximated by a family of classes $\mathcal{F}_k$, $k \geq 1$ (they are often associated with certain models, e.g., in regression or classification) and the empirical risk minimization problem (1.2) is replaced by a family of problems

$$(1.3) \qquad P_n f \to \min, \qquad f \in \mathcal{F}_k, \ k \geq 1.$$

The goal now is, based on solutions $\hat{f}_{n,k}$ of problems (1.3) and on the data, to construct an estimate $\hat{k}$ of index $k(P)$ of the "correct" model (i.e., a value of $k$ such that the solution of risk minimization problem (1.1) belongs to $\mathcal{F}_k$, or at least is well approximated by this class) and an "adaptive" solution $\hat{f} = \hat{f}_{n,\hat{k}}$ whose excess risk is close to being "optimal." The optimality of the solution is typically expressed by so-called *oracle inequalities* which, very roughly, show that the excess risk of $\hat{f}$ is within a constant from the excess risk of the solution one would have obtained with the help of an "oracle" who knows precisely to which of the classes $\mathcal{F}_k$ the true risk minimizer belongs [knows $k(P)$]. This way of thinking has become rather common in nonparametric statistics literature where various types of oracle inequalities have been proved, most often, in specialized settings (see [23] for a discussion on the subject).

The first general theory of empirical risk minimization was systematically developed by Vapnik and Chervonenkis [49] (see also [48] and references therein) in the late 1970s and early 1980s (although a number of more special results had been obtained much earlier, in particular, in connection with the development of the theory of maximum likelihood and $M$-estimation).

They obtained a number of bounds on $\mathcal{E}_P(\hat{f}_n)$ based on the inequality $\mathcal{E}_P(\hat{f}_n) \leq 2\|P_n - P\|_{\mathcal{F}}$ and on further bounding the sup-norm $\|P_n - P\|_{\mathcal{F}}$ in terms of random entropies or, now famous, VC-dimensions of the class $\mathcal{F}$ [here and in what follows $\|Y\|_{\mathcal{F}} := \sup_{f \in \mathcal{F}} |Y(f)|$ for $Y : \mathcal{F} \mapsto \mathbb{R}$]. They also developed more subtle bounds that provide an improvement in the case of small (in particular, zero) risk. These results played a significant role in the development of the general theory of empirical processes (see [16, 47]).

New developments in nonparametric statistics and, especially, in machine learning have motivated a number of improvements in the Vapnik–Chervonenkis theory of empirical risk minimization. Our approach largely relies on well-known papers of Birgé and Massart [8], Barron, Birgé and Massart [3], and on the more recent paper of Massart [36]. These authors proved a number of oracle inequalities for regression, density estimation and other nonparametric problems. More importantly, they suggested a rather general methodology of dealing with model selection for minimum contrast estimators that is based on Talagrand's concentration and deviation inequalities for empirical processes [42, 43], a new probabilistic tool at the time when these papers were written. Despite the fact that in many special statistical problems the use of Talagrand's inequalities can be avoided and oracle inequalities can be proved relying on more elementary probabilistic methods, one could hardly deny that concentration inequalities are the only universal tool in probability that suits the needs of model selection and oracle inequalities problems extremely well and are, probably, unavoidable when these problems are being dealt with in their full generality (e.g., in a machine learning setting). Talagrand's inequalities will be the main tool in this paper. Another important piece of work is the paper by Shen and Wong [39] where empirical processes methods were used to analyze empirical risk minimization on sieves (and, in particular, a version of iterative localization of excess risk bounds close to the approach discussed above was developed in a more specialized framework).

One of our main motivations was to understand better the results of Mammen and Tsybakov [35] on fast convergence rates in classification as well as more recent results of Tsybakov [44] and Tsybakov and van de Geer [45] on adaptation strategies for which these rates are attained. Our goal is to include these types of results in a more general framework of abstract empirical risk minimization (see Section 6). Another goal is to include into the same framework some other recent model selection results, especially in learning theory, where there is a definite need to develop general data-driven complexity penalization techniques suitable for neural networks, kernel machines and ensemble methods (see [28, 29, 30]). The analysis of convergence rates and the development of adaptive strategies for classification are currently at early stages (even consistency of boosting and kernel machines classification algorithms was established only recently; see [33, 40, 50]). Very recently, Bartlett, Jordan and McAuliffe [6] and Blanchard, Lugosi and Vayatis [10] obtained convergence rates of boosting-type classification methods based on convex risk minimization. Blanchard, Bousquet and Massart [9] obtained interesting oracle inequalities for penalized empirical risk minimization in kernel machines. It is of importance to develop better general ingredients of the proofs of such results so that it would be possible to concentrate on more specific difficulties related to the nature of the classification problem. These types of problems as well as a somewhat more general framework of convex risk minimization, including regression problems, are also within the scope of the methods of this paper (Sections 7, 8).

The proofs of all main results in the paper are given in Section 9.

## 2. Preliminaries.

2.1. *Talagrand's concentration inequalities.* Most of the results of the paper are based on famous concentration inequalities for empirical processes due to Talagrand [42, 43] (that provide uniform versions of classical Bernstein's-type inequalities for sums of i.i.d. random variables). We use the versions of these inequalities proved by Bousquet [13] and Klein [24] (see [11] for some other relevant inequalities). Namely, for a class $\mathcal{F}$ of measurable functions from $S$ into $[0, 1]$ (by a simple rescaling $[0, 1]$ can be replaced by any bounded interval) the following bounds hold for all $t > 0$:

- Bousquet's bound:

$$\mathbb{P}\left\{\|P_n - P\|_{\mathcal{F}} \geq \mathbb{E}\|P_n - P\|_{\mathcal{F}} + \sqrt{2\frac{t}{n}(\sigma_P^2(\mathcal{F}) + 2\mathbb{E}\|P_n - P\|_{\mathcal{F}})} + \frac{t}{3n}\right\} \leq e^{-t}.$$

- Klein's bound:

$$\mathbb{P}\left\{\|P_n - P\|_{\mathcal{F}} \leq \mathbb{E}\|P_n - P\|_{\mathcal{F}} - \sqrt{2\frac{t}{n}(\sigma_P^2(\mathcal{F}) + 2\mathbb{E}\|P_n - P\|_{\mathcal{F}})} - \frac{8t}{3n}\right\} \leq e^{-t}$$

(we modified Klein's bound slightly). Here $\sigma_P^2(\mathcal{F}) := \sup_{f \in \mathcal{F}}(Pf^2 - (Pf)^2)$.

2.2. *Empirical and Rademacher processes.* The empirical process is commonly defined as $n^{1/2}(P_n - P)$ and it is most often viewed as a stochastic process indexed by a function class $\mathcal{F}: n^{1/2}(P_n - P)(f)$, $f \in \mathcal{F}$ (see [16] or [47]). The Rademacher process indexed by a class $\mathcal{F}$ is defined as

$$R_n(f) := n^{-1} \sum_{i=1}^{n} \varepsilon_i f(X_i), \qquad f \in \mathcal{F},$$

$\{\varepsilon_i\}$ being i.i.d. Rademacher random variables (i.e., $\varepsilon_i$ takes the values $+1$ and $-1$ with probability $1/2$ each) independent of $\{X_i\}$. Roughly, $R_n(f)$ is the value of empirical correlation coefficient between $f(X_i)$, $i = 1, \ldots, n$ and Rademacher random noise. If $\|R_n\|_{\mathcal{F}}$ is large, it means that there exists $f \in \mathcal{F}$ for which $f(X_i)$ fits the noise well. Using such a class $\mathcal{F}$ in empirical risk minimization is likely to result in overfitting, which provides an intuitive explanation of the role of $\|R_n\|_{\mathcal{F}}$ as a complexity penalty in empirical risk minimization problems.

Rademacher processes have been widely used in the theory of empirical processes because of the following important inequality:

$$\frac{1}{2}\mathbb{E}\|R_n\|_{\mathcal{F}_c} \leq \mathbb{E}\|P_n - P\|_{\mathcal{F}} \leq 2\mathbb{E}\|R_n\|_{\mathcal{F}},$$

where $\mathcal{F}_c := \{f - Pf : f \in \mathcal{F}\}$. The upper bound is often referred to as a *symmetrization inequality* and the lower bound as a *desymmetrization inequality*. We will use this terminology in the future. These inequalities were brought into the theory of empirical processes by Giné and Zinn [21]. It is often convenient to use the desymmetrization inequality in combination with the following elementary lower bound:

$$\begin{aligned}
\mathbb{E}\|R_n\|_{\mathcal{F}_c} &\geq \mathbb{E}\|R_n\|_{\mathcal{F}} - \sup_{f \in \mathcal{F}} |Pf|\, \mathbb{E}|R_n(1)| \\
&\geq \mathbb{E}\|R_n\|_{\mathcal{F}} - \sup_{f \in \mathcal{F}} |Pf| \mathbb{E}^{1/2} \left| n^{-1} \sum_{j=1}^{n} \varepsilon_j \right|^2 \\
&\geq \mathbb{E}\|R_n\|_{\mathcal{F}} - \frac{\sup_{f \in \mathcal{F}} |Pf|}{\sqrt{n}}.
\end{aligned}$$

Rademacher processes possess many remarkable properties. In particular, they satisfy the following beautiful *contraction inequality*: if $\mathcal{F}$ is a class of functions with values in $[-1, 1]$, $\varphi$ is a function on $[-1, 1]$ with $\varphi(0) = 0$ and of Lipschitz norm bounded by 1, and $\varphi \circ \mathcal{F} := \{\varphi \circ f : f \in \mathcal{F}\}$, then $\mathbb{E}\|R_n\|_{\varphi \circ \mathcal{F}} \leq 2\mathbb{E}\|R_n\|_{\mathcal{F}}$ (follows from [31], Theorem 4.12). This implies, for instance, that

$$\mathbb{E} \sup_{f \in \mathcal{F}} \left| n^{-1} \sum_{i=1}^{n} \varepsilon_i f^2(X_i) \right| \leq 4 \mathbb{E} \sup_{f \in \mathcal{F}} \left| n^{-1} \sum_{i=1}^{n} \varepsilon_i f(X_i) \right|.$$

Concentration inequalities also apply to the Rademacher process since it can be viewed as an empirical process based on the sample $(X_1, \varepsilon_1), \ldots, (X_n, \varepsilon_n)$.

Often one needs to bound expected suprema of empirical and Rademacher processes. This can be done using various types of covering numbers (such as uniform covering numbers, random covering numbers, bracketing numbers, etc.) and the corresponding Dudley's entropy integrals. For instance, let $N(\mathcal{F}; L_2(P_n); \varepsilon)$ denote the minimal number of $L_2(P_n)$-balls of radius $\varepsilon$ covering $\mathcal{F}$. Suppose that $\forall f \in \mathcal{F}, \forall x \in S : |f(x)| \leq F(x) \leq U$, where $U > 0$ and $F$ is a measurable function (called an envelope of $\mathcal{F}$). Let $\sigma^2 := \sup_{f \in \mathcal{F}} Pf^2$. If for some $A > 0$, $V > 0$

$$(2.1) \qquad \forall \varepsilon > 0 \qquad N(\mathcal{F}; L_2(P_n); \varepsilon) \leq \left(\frac{A\|F\|_{L_2(P_n)}}{\varepsilon}\right)^V,$$

then with some universal constant $C > 0$ (for $\sigma^2 \geq \text{const } n^{-1}$)

$$(2.2) \quad \mathbb{E}\|R_n\|_{\mathcal{F}} \leq C\left[\sqrt{\frac{V}{n}}\sigma\sqrt{\log\frac{A\|F\|_{L_2(P)}}{\sigma}} \vee \frac{VU}{n}\log\frac{A\|F\|_{L_2(P)}}{\sigma}\right].$$

If for some $A > 0, \rho \in (0, 1)$

$$(2.3) \qquad \forall \varepsilon > 0 \qquad \log N(\mathcal{F}; L_2(P_n); \varepsilon) \leq \left(\frac{A\|F\|_{L_2(P_n)}}{\varepsilon}\right)^{2\rho},$$

then

$$(2.4) \quad \mathbb{E}\|R_n\|_{\mathcal{F}} \leq C\left[\frac{A^\rho\|F\|^\rho_{L_2(P)}}{\sqrt{n}}\sigma^{1-\rho} \vee \frac{A^{2\rho/(\rho+1)}\|F\|^{2\rho/(\rho+1)}_{L_2(P)}U^{(1-\rho)/(1+\rho)}}{n^{1/(1+\rho)}}\right].$$

The proofs of these types of bounds can be found in [17, 18, 20, 37, 41]; the current version of (2.4) is due to Giné and Koltchinskii [19]).

In particular, if $\mathcal{F}$ is a VC-subgraph class, then the condition (2.1) holds (in fact, the condition holds even for the uniform covering numbers) and one can use the bound (2.2). We will call the function classes satisfying (2.1) *VC-type classes*. If $\mathcal{H}$ is VC-type, then its convex hull $\text{conv}(\mathcal{H})$ satisfies (2.3) with $\rho := \frac{V}{V+2}$ (see [47]), so one can use the bound (2.4) for $\mathcal{F} \subset \text{conv}(\mathcal{H})$ (note that one should use the envelope $F$ of the class $\mathcal{H}$ itself for its convex hull as well). Many other useful bounds on expected suprema of empirical and Rademacher processes (in particular, in terms of bracketing numbers) can be found in [47] and [16].

2.3. *The $\sharp$-transform and related questions.* In this section, we introduce and discuss some useful transformations, involved in the definitions of various complexity measures of function classes in empirical risk minimization. As it has been already pointed out in the Introduction, the excess risk bounds are often based on solving the fixed point equation, or, more generally, equations of the type $\psi(\delta) = \varepsilon\delta$, for $\psi(\cdot) = U_n(\cdot; t)$. This naturally leads to the following definitions.

For a function $\psi: \mathbb{R}_+ \mapsto \mathbb{R}_+$, define

$$\psi^\flat(\delta) := \sup_{\sigma \geq \delta} \frac{\psi(\sigma)}{\sigma} \quad \text{and} \quad \psi^\sharp(\varepsilon) := \inf\{\delta > 0 : \psi^\flat(\delta) \leq \varepsilon\}.$$

We will call these transformations, respectively, the $\flat$-transform and the $\sharp$-transform of $\psi$. We are mainly interested in the $\sharp$-transform. It has the following properties whose proofs are elementary and straightforward:

1. Suppose that $\psi(u) = o(u)$ as $u \to \infty$. Then the function $\psi^\sharp$ is defined on $(0, +\infty)$ and is a nonincreasing function on this interval.

2. If $\psi_1 \leq \psi_2$, then $\psi_1^\sharp \leq \psi_2^\sharp$. Moreover, it is enough to assume that $\psi_1(\delta) \leq \psi_2(\delta)$ either for all $\delta \geq \psi_2^\sharp(\varepsilon)$, or for all $\delta \geq \psi_1^\sharp(\varepsilon) - \tau$ with an arbitrary $\tau > 0$, to conclude that $\psi_1^\sharp(\varepsilon) \leq \psi_2^\sharp(\varepsilon)$.

3. For $a > 0$, $(a\psi)^\sharp(\varepsilon) = \psi^\sharp(\varepsilon/a)$.

4. If $\varepsilon = \varepsilon_1 + \cdots + \varepsilon_m$, then

$$\psi_1^\sharp(\varepsilon) \vee \cdots \vee \psi_m^\sharp(\varepsilon) \leq (\psi_1 + \cdots + \psi_m)^\sharp(\varepsilon) \leq \psi_1^\sharp(\varepsilon_1) \vee \cdots \vee \psi_m^\sharp(\varepsilon_m).$$

5. If $\psi(u) \equiv c$, then $\psi^\sharp(\varepsilon) = c/\varepsilon$.

6. If $\psi(u) := u^\alpha$ with $\alpha \leq 1$, then $\psi^\sharp(\varepsilon) := \varepsilon^{-1/(1-\alpha)}$.

7. For $c > 0$, let $\psi_c(\delta) := \psi(c\delta)$. Then $\psi_c^\sharp(\varepsilon) = \frac{1}{c}\psi^\sharp(\varepsilon/c)$. If $\psi$ is nondecreasing and $c \geq 1$, then this easily implies that $c\psi^\sharp(u) \leq \psi^\sharp(u/c)$.

8. For $c > 0$, let now $\psi_c(\delta) := \psi(\delta + c)$. Then for all $u > 0$, $\varepsilon \in (0, 1]$, $\psi_c^\sharp(u) \leq \psi^\sharp(\varepsilon u/2) - c \vee c\varepsilon$.

Let us call $\psi : \mathbb{R}_+ \mapsto \mathbb{R}_+$ a function of concave type if it is nondecreasing and $u \mapsto \frac{\psi(u)}{u}$ is decreasing. If, in addition, for some $\gamma \in (0, 1)$, $u \mapsto \frac{\psi(u)}{u^\gamma}$ is decreasing, $\psi$ will be called a function of strictly concave type (with exponent $\gamma$). In particular, if $\psi(u) := \varphi(u^\gamma)$, or $\psi(u) := \varphi^\gamma(u)$, where $\varphi$ is a nondecreasing strictly concave function with $\varphi(0) = 0$, then $\psi$ is of concave type for $\gamma = 1$ and of strictly concave type for $\gamma < 1$.

9. If $\psi$ is of concave type, then $\psi^\sharp$ is the inverse of the function $\delta \mapsto \frac{\psi(\delta)}{\delta}$. In this case, $\psi^\sharp(cu) \geq \psi^\sharp(u)/c$ for $c \leq 1$ and $\psi^\sharp(cu) \leq \psi^\sharp(u)/c$ for $c \geq 1$.

10. If $\psi$ is of strictly concave type with exponent $\gamma$, then for $c \leq 1$, $\psi^\sharp(cu) \leq \psi^\sharp(u) c^{-\frac{1}{1-\gamma}}$.

It will be convenient sometimes to discretize the supremum in the definition of $\psi^\flat$. Namely, let $q > 1$ and $\delta_j := q^{-j}$, $j \in \mathbb{Z}$. Define

$$\psi^{\flat,q}(\delta) := \sup_{\delta_j \geq \delta} \frac{\psi(\delta_j)}{\delta_j}, \qquad \psi^{\sharp,q}(\varepsilon) := \inf\{\delta > 0 : \psi^{\flat,q}(\delta) \leq \varepsilon\}$$

and

$$\psi_{[0,1]}^{\flat,q}(\delta) := \sup_{1 \geq \delta_j \geq \delta} \frac{\psi(\delta_j)}{\delta_j}, \qquad \psi_{[0,1]}^{\sharp,q}(\varepsilon) := \inf\{\delta \in (0,1] : \psi_{[0,1]}^{\flat,q}(\delta) \leq \varepsilon\}$$

(if in the last definition $\psi_{[0,1]}^{\flat,q}(\delta)$ is larger than $\varepsilon$ for all $\delta \leq 1$, then we set $\psi_{[0,1]}^{\sharp,q}(\varepsilon) := 1$).

Properties 1–4 and 7 hold for $\psi^{\sharp,q}$ with the following obvious changes. In property 2, it is enough to assume that $\psi_1(\delta) \leq \psi_2(\delta)$ only for $\delta = \delta_j$ and the second part of this property should be formulated as follows: if $\psi_1(\delta) \leq \psi_2(\delta)$ either for all $\delta \geq \psi_2^{\sharp,q}(\varepsilon)$, or for all $\delta \geq q^{-1}\psi_1^{\sharp,q}(\varepsilon)$, then $\psi_1^{\sharp,q}(\varepsilon) \leq \psi_2^{\sharp,q}(\varepsilon)$. Property 7 holds with $c = q^j$ for any $j$. We will refer to these properties as $1'$–$4'$ and $7'$ in what follows.

Also, the following simple fact is true:

11. If $\psi$ is nondecreasing, then $\psi^{\flat,q}(\varepsilon) \leq \psi^{\sharp,q}(\varepsilon) \leq \psi^{\sharp}(\varepsilon) \leq \psi^{\sharp,q}(\varepsilon/q)$. In addition, if $\psi(\delta) = \text{const}$ for $\delta \geq 1$ (which will be the case in many situations), then $\psi_{[0,1]}^{\sharp,q}(\varepsilon) = \psi^{\sharp,q}(\varepsilon)$.

We conclude this section with a simple proposition, describing useful properties of functions of strictly concave type.

PROPOSITION 1. (i) *If $\psi$ is a function of strictly concave type with some exponent $\gamma \in (0,1)$, then*

$$\sum_{j:\delta_j \geq \delta} \frac{\psi(\delta_j)}{\delta_j} \leq c_{\gamma,q} \frac{\psi(\delta)}{\delta},$$

*where $c_{\gamma,q}$ is a constant depending only on $q, \gamma$.*

(ii) *Under the same assumptions, the equation $\psi(\delta) = \delta$ has unique solution $\bar{\delta}$. Suppose $\bar{\delta} \leq 1$ and define $\bar{\delta}_0 := 1$, $\bar{\delta}_{k+1} := \psi(\bar{\delta}_k) \wedge 1$. Then $\{\bar{\delta}_k\}$ is a nonincreasing sequence converging to $\bar{\delta}$ and, for all $k$, $\bar{\delta}_k - \bar{\delta} \leq \bar{\delta}^{1-\gamma^k}(1-\bar{\delta})^{\gamma^k}$.*

2.4. *Empirical and Rademacher complexities.* The most natural complexity penalties in risk minimization problems are based on expected supnorms of the empirical process over the whole class $\mathcal{F}$ or its subsets. However, such complexities are *distribution dependent*, so it is hard to use them in model selection. The idea to use Rademacher processes to construct *data-dependent* complexity penalties in model selection problems of learning theory was suggested independently by Koltchinskii [26] and Bartlett, Boucheron and Lugosi [4]. It is based on the following simple observation: if one combines the symmetrization inequality with concentration inequalities for empirical and Rademacher processes (in fact, with simpler Hoeffding-type concentration inequalities based on the martingale difference approach), one can get the following bound:

$$\mathbb{P}\bigg\{\|P_n - P\|_{\mathcal{F}} \geq 2\|R_n\|_{\mathcal{F}} + \frac{3t}{\sqrt{n}}\bigg\} \leq \exp\bigg\{-\frac{2t^2}{3}\bigg\}, \qquad t > 0.$$

Quite similarly, using instead the desymmetrization inequality one can get a simple lower confidence bound on $\|P_n - P\|_{\mathcal{F}}$ in terms of $\|R_n\|_{\mathcal{F}}$. Since the Rademacher process does not involve the unknown distribution directly and can be computed based only on the data, one can use $\|R_n\|_{\mathcal{F}}$ as a datadependent measure of the accuracy of approximation of the true distribution $P$ by the empirical distribution $P_n$ uniformly over the class. Essentially, this

justifies using $\|R_n\|_{\mathcal{F}}$ as a bootstrap-type complexity penalty associated with the class $\mathcal{F}$ (although Rademacher bootstrap is not asymptotically correct). The main problem, however, is that such *global* complexities as $\|R_n\|_{\mathcal{F}}$ do not allow one to recover the convergence rates in risk minimization problems. Typically, $\|R_n\|_{\mathcal{F}}$ would be of the order $O(n^{-1/2})$ (this is the case, e.g., for VC-classes and, more generally, for Donsker classes of functions). The convergence rates in many risk minimization problems are often faster than this and they are related to the behavior of the *continuity modulus* of the empirical process $n^{1/2}(P_n - P)$ rather than to the behavior of its sup-norm (see [36]). Thus, relevant data-dependent complexities could be based on the continuity modulus of the Rademacher process that mimics the properties of the empirical process. As we will see later, the complexities of this type are defined as the $\sharp$-transform of the corresponding (expected) continuity modulus.

Let $\rho_P : L_2(P) \times L_2(P) \mapsto [0, +\infty)$ be a function such that

$$\rho_P^2(f, g) \geq P(f-g)^2 - (P(f-g))^2, \qquad f, g \in L_2(P).$$

Typically $\rho_P$ will be also a (pseudo)metric, for instance, $\rho_P^2(f,g) = P(f-g)^2$ or $\rho_P^2(f,g) = P(f-g)^2 - (P(f-g))^2$.

Given a function $Y : \mathcal{F} \mapsto \mathbb{R}$, define its continuity moduli (local and global) as follows:

$$\omega_{\rho_P}(Y; f; \delta) := \sup_{g \in \mathcal{F}, \rho_P(g,f) \leq \delta} |Y(g) - Y(f)| \quad \text{and}$$

$$\omega_{\rho_P}(Y; \delta) := \sup_{f,g \in \mathcal{F}, \rho_P(f,g) \leq \delta} |Y(f) - Y(g)|.$$

Assume, for simplicity, that the infimum of $Pf$ over $\mathcal{F}$ is attained at a function $\bar{f} \in \mathcal{F}$ (we are assuming this in what follows whenever it is needed; otherwise, the definitions can be easily modified). Let

$$\theta_n(\delta) := \theta_n(\mathcal{F}; \bar{f}; \delta) := \mathbb{E}\omega_{\rho_P}(P_n - P; \bar{f}; \sqrt{\delta}).$$

The empirical complexity, such as the ones previously used in [5, 14, 27, 36], can be now defined as $\theta_n^{\sharp}(\varepsilon)$ where $\varepsilon$ is a numerical constant (often, $\varepsilon = 1$, which corresponds to the fixed point equation, but sometimes the dependence on $\varepsilon$ is of importance). The function $\theta_n(\delta)$ in this definition can be replaced by $\sup_{f \in \mathcal{F}} \mathbb{E}\omega_{\rho_P}(P_n - P; f; \sqrt{\delta})$, or even by $\mathbb{E}\omega_{\rho_P}(P_n - P; \sqrt{\delta})$, without increasing the complexity significantly (at least, in most of the relevant examples).

It will be shown in the next sections how to use these types of quantities to provide upper bounds on the excess risk. Now, we utilize the Rademacher process to construct data-dependent bounds on $\theta_n^{\sharp}(\varepsilon)$. Suppose that $\rho_P^2(f, g) := P(f-g)^2$. Define

$$\bar{\omega}_n(\delta) := \mathbb{E}\omega_{\rho_P}(R_n; \sqrt{\delta}), \qquad \hat{\omega}_n(\delta) := \omega_{\rho_{P_n}}(R_n; \sqrt{\delta}),$$

$$\hat{\omega}_{n,r}(\delta) := \mathbb{E}_\varepsilon \omega_{\rho_{P_n}}(R_n; \sqrt{\delta}),$$

where $\mathbb{E}_\varepsilon$ denotes the expectation only with respect to the Rademacher sequence $\{\varepsilon_i\}$.

The next lemma is pretty much akin to some statements in [5]. Koltchinskii and Panchenko [27] proved some results in this direction in a more specialized setting of function learning (in zero error case). We give its proof in Section 9 for completeness and also because a similar approach is used in the proofs of several other results given below.

LEMMA 1. *For $q > 1$, there exist constants $C, c > 0$ (depending only on $q$) such that*

$$\forall \varepsilon > 0 \qquad \theta_n^\sharp(\varepsilon) \leq \bar{\omega}_n^\sharp(\varepsilon/2)$$

*and for all $\varepsilon \in (0, 1]$*

$$\mathbb{P}\left\{\bar{\omega}_n^\sharp(\varepsilon) \geq C\left(\hat{\omega}_n^\sharp(c\varepsilon) + \frac{t}{n\varepsilon^2}\right)\right\} \leq 2\log_q \frac{qn}{t} e^{-t},$$

$$\mathbb{P}\left\{\hat{\omega}_n^\sharp(\varepsilon) \geq C\left(\bar{\omega}_n^\sharp(c\varepsilon) + \frac{t}{n\varepsilon^2}\right)\right\} \leq 2\log_q \frac{qn}{t} e^{-t}.$$

*The same is true with $\hat{\omega}_n^\sharp$ replaced by $\hat{\omega}_{n,r}^\sharp$.*

2.5. *Examples.* We give below several simple bounds on local Rademacher complexities $\theta_n^\sharp(\varepsilon)$, $\varepsilon \in (0, 1]$ that are of interest in applications and have been discussed, for example, in [5, 6, 10, 36].

EXAMPLE 1 (Finite-dimensional classes). Suppose that $\mathcal{F}$ is a subset of a *finite-dimensional* subspace $L$ of $L_2(P)$ with $\dim(L) = d$. Then $\theta_n(\delta) \leq (\delta d/n)^{1/2}$ and $\theta_n^\sharp(\varepsilon) \leq d/(n\varepsilon^2)$. Indeed, if $e_1, \ldots, e_d$ is an orthonormal basis of $L$, and $g, \bar{g} \in L$, $g = \sum_{i=1}^d \alpha_i e_i$, $\bar{g} = \sum_{i=1}^d \bar{\alpha}_i e_i$, then $\|g - \bar{g}\|_{L_2(\Pi)}^2 = \sum_{i=1}^d (\alpha_i - \bar{\alpha}_i)^2$. Therefore, using the Cauchy–Schwarz inequality,

$$\theta_n(\delta) = \mathbb{E} \sup_{g \in \mathcal{F}, \|g - \bar{g}\|_{L_2(P)} \leq \sqrt{\delta}} |(P_n - P)(g - \bar{g})|$$

$$\leq \mathbb{E} \sup_{\sum_{i=1}^d (\alpha_i - \bar{\alpha}_i)^2 \leq \delta} \left|\sum_{i=1}^d (\alpha_i - \bar{\alpha}_i)(P_n - P)(e_i)\right|$$

$$\leq \sqrt{\delta}\left(\sum_{i=1}^d \mathbb{E}(P_n - P)^2(e_i)\right)^{1/2} \leq \sqrt{\frac{\delta d}{n}},$$

and the second bound on $\theta_n^\sharp(\varepsilon)$ is now immediate due to the properties of $\sharp$-transform.

EXAMPLE 2 (Ellipsoids in $L_2$). This is a simple generalization of the previous example. Suppose that $\mathcal{F} := \{Tg : \|g\|_{L_2(P)} \leq 1\}$, where $T : L_2(P) \mapsto L_2(P)$ is a Hilbert–Schmidt operator with Hilbert–Schmidt norm $\|T\|_{HS}$ and such that its operator norm $\|T\| \leq 1$. Thus, $\mathcal{F}$ is an ellipsoid in Hilbert space $L_2(P)$. Suppose also that $\text{Ker}(T) = \{0\}$, and, for $f_1 = Tg_1$, $f_2 = Tg_2$, we define $\rho_P(f_1, f_2) = \|g_1 - g_2\|_{L_2(P)}$. Then, the same argument as in the previous example yields $\theta_n(\delta) \leq (\delta \|T\|_{HS}^2/n)^{1/2}$ and $\theta_n^\sharp(\varepsilon) \leq \|T\|_{HS}^2/(n\varepsilon^2)$.

Often, it is natural to use Dudley's entropy integral to bound the function $\theta_n(\delta)$ and then to derive a bound on $\theta_n^\sharp(\varepsilon)$. Various notions of the entropy of function class $\mathcal{F}$ can be used for this purpose (entropy with bracketing, random entropy, uniform entropy, etc.). This technique is standard in the theory of empirical processes and can be found, for example, in the book of Van der Vaart and Wellner [47]. Here are some examples of the bounds based on this approach.

EXAMPLE 3 (VC-type classes). Suppose that $\mathcal{F}$ is a VC-type class, that is, the condition (2.1) is satisfied (in particular, $\mathcal{F}$ might be a VC-subgraph class). Assume for simplicity that $F \equiv U = 1$. Then it follows from (2.2) that

$$\theta_n(\delta) \leq K\left(\sqrt{\frac{V\delta}{n}}\sqrt{\log \frac{1}{\delta}} \vee \frac{V}{n}\log \frac{1}{\delta}\right),$$

which leads to the following bound: $\theta_n^\sharp(\varepsilon) \leq CV/(n\varepsilon^2) \log(n\varepsilon^2/V)$.

EXAMPLE 4 (Entropy conditions). In the case when the entropy of the class (uniform, bracketing, etc.) is bounded by $O(\varepsilon^{-2\rho})$ for some $\rho \in (0,1)$, we typically have $\theta_n^\sharp(\varepsilon) = O(n^{-1/(1+\rho)})$. For instance, if (2.3) holds, then it follows from (2.4) (with $F \equiv U = 1$ for simplicity) that

$$\theta_n(\delta) \leq K\left(\frac{A^\rho}{\sqrt{n}}\delta^{(1-\rho)/2} \vee \frac{A^{2\rho/(\rho+1)}}{n^{1/(1+\rho)}}\right).$$

Therefore, $\theta_n^\sharp(\varepsilon) \leq CA^{2\rho/(1+\rho)}/(n\varepsilon^2)^{1/(1+\rho)}$.

EXAMPLE 5 (Convex hulls). If $\mathcal{F} := \text{conv}(\mathcal{H}) := \{\sum_j \lambda_j h_j : \sum_j |\lambda_j| \leq 1, h_j \in \mathcal{H}\}$ is the symmetric convex hull of a given VC-type class $\mathcal{H}$ of measurable functions from $S$ into $[0,1]$, then the condition of the previous example is satisfied with $\rho := \frac{V}{V+2}$. This yields $\theta_n^\sharp(\varepsilon) \leq (K(V)/(n\varepsilon^2))^{\frac{1}{2}\frac{2+V}{1+V}}$.

EXAMPLE 6 (Shattering numbers for classes of binary functions). Let $\mathcal{F}$ be a class of binary functions, that is, functions $f : S \mapsto \{0,1\}$. Let

$$\Delta^\mathcal{F}(X_1, \ldots, X_n) := \text{card}(\{(f(X_1), \ldots, f(X_n)) : f \in \mathcal{F}\})$$

be the shattering number of the class $\mathcal{F}$ on the sample $(X_1,\ldots,X_n)$. Using a bound that can be found in [36], we get

$$\theta_n(\delta) \leq K\left[\sqrt{\delta \frac{\mathbb{E}\log\Delta^{\mathcal{F}}(X_1,\ldots,X_n)}{n}} + \frac{\mathbb{E}\log\Delta^{\mathcal{F}}(X_1,\ldots,X_n)}{n}\right],$$

which easily yields

$$\theta_n^{\sharp}(\varepsilon) \leq C\frac{\mathbb{E}\log\Delta^{\mathcal{F}}(X_1,\ldots,X_n)}{n\varepsilon^2}.$$

EXAMPLE 7 (Mendelson's complexities for kernel machines). Let $K$ be a symmetric nonnegatively definite kernel on $S \times S$ and let $H_K$ be the corresponding reproducing kernel Hilbert space, that is, $H_K$ is the closure of the set of linear combinations $\sum_i \alpha_i K(x_i,\cdot)$, $x_i \in S$, $\alpha_i \in \mathbb{R}$ with respect to the norm $\|\cdot\|_K$ defined as

$$\left\|\sum_i \alpha_i K(x_i,\cdot)\right\|_K^2 = \sum_{i,j} \alpha_i \alpha_j K(x_i, x_j).$$

Suppose that $\mathcal{F} := B_K$ is the unit ball in $H_K$. Such classes are frequently used in learning theory for kernel machines. Let $\lambda_i$ be the eigenvalues of the integral operator generated by $K$ in space $L_2(P)$. The following is a version of bounds of Mendelson [37]:

$$C_1\left(n^{-1}\sum_{j=1}^{\infty}\lambda_j \wedge \delta\right)^{1/2} \leq \bar{\omega}_n(\delta) = \mathbb{E}\sup_{P(f-g)^2\leq\delta, f,g\in\mathcal{F}}|R_n(f-g)|$$

$$\leq C_2\left(n^{-1}\sum_{j=1}^{\infty}\lambda_j \wedge \delta\right)^{1/2}$$

with some numerical constants $C_1, C_2 > 0$. Similarly, if $\lambda_i^{(n)}$, $i=1,\ldots,n$ are the eigenvalues of the matrix $(n^{-1}K(X_i,X_j): 1\leq i,j\leq n)$, then Mendelson's argument also gives

$$C_1\left(n^{-1}\sum_{j=1}^{n}\lambda_j^{(n)} \wedge \delta\right)^{1/2} \leq \hat{\omega}_{n,r}(\delta) = \mathbb{E}_\varepsilon \sup_{P_n(f-g)^2\leq\delta, f,g\in\mathcal{F}}|R_n(f-g)|$$

$$\leq C_2\left(n^{-1}\sum_{j=1}^{n}\lambda_j^{(n)} \wedge \delta\right)^{1/2}.$$

Denote the true and empirical Mendelson's complexities by

$$\bar{\gamma}_n(\delta) = \gamma_n(\mathcal{F};\delta) = \left(n^{-1}\sum_{j=1}^{\infty}\lambda_j \wedge \delta\right)^{1/2} \quad \text{and}$$

$$\hat{\gamma}_n(\delta) = \hat{\gamma}_n(\mathcal{F};\delta) = \left(n^{-1}\sum_{j=1}^n \lambda_j^{(n)} \wedge \delta\right)^{1/2}.$$

Note that these functions are strictly concave, nondecreasing and are equal to 0 for $\delta = 0$. Moreover, they are both square roots of concave functions and, hence, they are of strictly concave type. The properties of $\sharp$-transform imply that with some constants $c_1, c_2$

$$\bar{\gamma}_n^\sharp(c_1\varepsilon) \leq \bar{\omega}_n^\sharp(\varepsilon) \leq \bar{\gamma}_n^\sharp(c_2\varepsilon) \quad \text{and} \quad \hat{\gamma}_n^\sharp(c_1\varepsilon) \leq \hat{\omega}_{n,r}^\sharp(\varepsilon) \leq \hat{\gamma}_n^\sharp(c_2\varepsilon).$$

Together with Lemma 1, this allows one to use empirical Mendelson's complexity as an estimate of true Mendelson's complexity.

**3. First excess risk bounds.** The idea to express excess risk bounds in terms of solutions of fixed point equations for continuity modulus of empirical or Rademacher processes and also to relate them to ratio-type inequalities has been around for a while (see [5, 27, 36]). Comparing with the recent work of Bartlett, Bousquet and Mendelson [5], our approach in this section relates the excess risk bounds more directly to the diameter of the $\delta$-minimal set of $P$ (recall the definitions in Section 1) and also provides ratio-type inequalities for the empirical excess risk expressed in terms of $\sharp$-transform of the function $\bar{U}_n(\delta;t)$ involved in Talagrand's inequality. The excess empirical risk is defined as $\hat{\mathcal{E}}_n(f) := \mathcal{E}_{P_n}(f)$ and the $\delta$-minimal set of $P_n$ as $\hat{\mathcal{F}}_n(\delta) := \mathcal{F}_{P_n}(\delta)$. Also, denote $\mathcal{F}(s,r] := \mathcal{F}_P(s,r] := \mathcal{F}(r) \setminus \mathcal{F}(s)$.

Let $\hat{f}_n := \operatorname{argmin}_{f \in \mathcal{F}} P_n f$ be an empirical risk minimizer [i.e., a solution of (1.2)]. For simplicity, we assume that it exists, although the results can be easily modified for approximate solutions of (1.2). Recall that $D(\delta) := D_P(\mathcal{F};\delta) := \sup_{f,g \in \mathcal{F}(\delta)} \rho_P(f,g)$ denotes the $\rho_P$-diameter of the $\delta$-minimal set and also that

$$\phi_n(\delta) := \phi_n(\mathcal{F};P;\delta) := \mathbb{E} \sup_{f,g \in \mathcal{F}(\delta)} |(P_n - P)(f-g)|.$$

Let

$$U_n(\delta;t) := U_{n,t}(\delta) := \phi_n(\delta) + \sqrt{2\frac{t}{n}(D^2(\delta) + 2\phi_n(\delta))} + \frac{t}{2n}.$$

Finally, let us fix $q > 1$ and define $V_n$ and $\delta_n(t)$ as follows:

$$V_n(\delta;t) := V_{n,t}(\delta) := U_{n,t}^{\flat,q}(\delta) \quad \text{and} \quad \delta_n(t) := U_{n,t}^{\sharp,q}\left(\frac{1}{2q}\right).$$

Whenever it is needed, we will write $\delta_n(\mathcal{F};t)$ or $\delta_n(\mathcal{F};P;t)$ to emphasize the dependence of these types of quantities on function class and on distribution. The following result gives an upper bound on the excess risk of $\hat{f}_n$ and also provides uniform bounds on the ratios of the empirical excess risk of a function $f \in \mathcal{F}$ to its true excess risk.

THEOREM 1. *For all $t > 0$ and all $\delta \geq \delta_n(t)$*

$$\mathbb{P}\{\mathcal{E}(\hat{f}_n) \geq \delta\} \leq \log_q \frac{q}{\delta} e^{-t} \quad and$$

$$\mathbb{P}\left\{\sup_{f \in \mathcal{F}, \mathcal{E}(f) \geq \delta} \left| \frac{\hat{\mathcal{E}}_n(f)}{\mathcal{E}(f)} - 1 \right| \geq q V_n(\delta; t) \right\} \leq \log_q \frac{q}{\delta} e^{-t}.$$

Almost as in Section 2, define the expected continuity modulus

$$\omega_n(\mathcal{F}; \delta) := \mathbb{E} \sup_{\rho_P(f,g) \leq \delta, f, g \in \mathcal{F}} |(P_n - P)(f - g)|.$$

Since $\phi_n(\delta) \leq \omega_n(\mathcal{F}; D(\delta))$, the behavior of $\phi_n$ can be determined by $\omega_n$ and $D$. If $\mathcal{F}$ is a $P$-Donsker class, then, by asymptotic equicontinuity of empirical processes,

$$\lim_{\delta \to 0} \limsup_n n^{1/2} \omega_n(\mathcal{F}; \delta) = 0.$$

This fact and the definition of $\delta_n(t)$ immediately imply that $\delta_n(t) = o(n^{-1/2})$ as soon as $\mathcal{F}$ is $P$-Donsker and $D(\delta) \to 0$. The last condition is natural if the risk minimization problem (1.1) has unique solution. Moreover, there exists a sequence $t_n \to \infty$ such that $\delta_n(t_n) = o(n^{-1/2})$. Thus, by Theorem 1, we can conclude that $\mathcal{E}_P(\hat{f}_n) = o_P(n^{-1/2})$ whenever the empirical risk minimization occurs over a $P$-Donsker class and $D(\delta) \to 0$. This observation shows that convergence rates of the excess risk faster than $n^{-1/2}$ (that came as a surprise in classification problems in nonzero error case several years ago) are, in fact, typical in general empirical risk minimization over Donsker classes.

In the case when the function $f \mapsto Pf$ has the unique minimum in $\mathcal{F}$ (i.e., the minimal set $\mathcal{F}(0)$ consists precisely of one element), the quantity $\delta_n(t)$ often gives correct (in a minimax sense) convergence rate in risk minimization problems (see Section 6.1). However, if $\mathcal{F}(0)$ consists of more than one function, then the diameter $D(\delta)$ of the $\delta$-minimal set becomes bounded away from 0 and as a result $\delta_n(t)$ cannot be smaller than $c\sqrt{\frac{t}{n}}$ (and the optimal convergence rate is often better than this, e.g., in classification problems). In the next section, we study more subtle geometric characteristics of the class $\mathcal{F}$ that might be used in such cases to recover the correct convergence rates.

An important consequence of Theorem 1 is the following lemma that shows that $\delta$-minimal sets can be estimated by empirical $\delta$-minimal sets provided that $\delta$ is not too small.

LEMMA 2. *For all $t > 0$, there exists an event of probability at least $1 - \log_q \frac{q^2}{\delta_n(t)} e^{-t}$ such that on this event $\forall \delta \geq \delta_n(t) : \mathcal{F}(\delta) \subset \hat{\mathcal{F}}_n(3\delta/2)$ and $\hat{\mathcal{F}}_n(\delta) \subset \mathcal{F}(2\delta)$.*

Note that, as follows from the definition, $\delta_n(t) \geq \frac{t}{n}$, so the probabilities in Theorem 1 are, in fact, upper bounded by $\log_q \frac{n}{t} \exp\{-t\}$ (which depends neither on the class $\mathcal{F}$, nor on $P$). The logarithmic factor in front of the exponent, most often, does not spoil the bound since in typical applications $\delta_n(t)$ is upper bounded by $\delta_n + \frac{t}{n}$, where $\delta_n$ is larger than $\frac{\log \log n}{n}$. Adding $\log \log n$ to $t$ is enough to eliminate the influence of the logarithm. However, if $\delta_n = O(n^{-1})$, the logarithmic factor would create a problem. It is good to know that it can be eliminated under extra conditions on $\phi_n(\delta)$ and $D(\delta)$. More precisely, assume that $\phi_n(\delta) \leq \check{\phi}_n(\delta)$ and $D(\delta) \leq \check{D}(\delta)$, $\delta > 0$, where $\check{\phi}_n$ is a function of strictly concave type with some exponent $\gamma \in (0, 1)$ and $\check{D}$ is a concave-type function (see the definitions in Section 2.3). Define

$$\check{U}_n(\delta; t) := \check{U}_{n,t}(\delta) := \check{K}\left(\check{\phi}_n(\delta) + \check{D}(\delta)\sqrt{\frac{t}{n}} + \frac{t}{n}\right)$$

with some numerical constant $\check{K}$. Then $\check{U}_n(\cdot; t)$ is a concave-type function. In this case, it is natural to define

$$\check{V}_n(\delta; t) := \check{U}_{n,t}^{\flat}(\delta) = \frac{\check{U}_n(\delta; t)}{\delta} \quad \text{and} \quad \check{\delta}_n(t) := \check{U}_{n,t}^{\sharp}\left(\frac{1}{q}\right).$$

THEOREM 2. *There exists a constant $\check{K}$ such that for all $t > 0$ and for all $\delta \geq \check{\delta}_n(t)$,*

$$\mathbb{P}\{\mathcal{E}(\hat{f}_n) \geq \delta\} \leq e^{-t} \quad \text{and} \quad \mathbb{P}\left\{\sup_{f \in \mathcal{F}, \mathcal{E}(f) \geq \delta} \left|\frac{\hat{\mathcal{E}}_n(f)}{\mathcal{E}(f)} - 1\right| \geq q\check{V}_n(\delta; t)\right\} \leq e^{-t}.$$

In what follows we do not use this refinement except in several cases when it is really needed.

Now we outline a way to define the empirical version of $\delta_n(t)$. To this end, it will be convenient to choose $\rho_P^2(f, g) := P(f - g)^2$. Note that

$$U_n(\delta; t) \leq \bar{U}_n(\delta; t) := \bar{U}_{n,t}(\delta) := \bar{K}\left(\phi_n(\delta) + D(\delta)\sqrt{\frac{t}{n}} + \frac{t}{n}\right),$$

where $\bar{K} = 2$. Hence, if we define $\bar{\delta}_n(t) := \bar{U}_{n,t}^{\sharp;q}(1/2q^3)$, then it follows from the definitions that $\delta_n(t) \leq \bar{\delta}_n(t)$.

Define the empirical versions of the functions $D$ and $\phi_n$ as follows:

$$\hat{D}_n(\delta) := \sup_{f,g \in \hat{\mathcal{F}}_n(\delta)} \rho_{P_n}(f, g) \quad \text{and} \quad \hat{\phi}_n(\delta) := \sup_{f,g \in \hat{\mathcal{F}}_n(\delta)} |R_n(f - g)|.$$

Let

$$\hat{U}_n(\delta;t) := \hat{U}_{n,t}(\delta) := \hat{K}\left(\hat{\phi}_n(\hat{c}\delta) + \hat{D}_n(\hat{c}\delta)\sqrt{\frac{t}{n}} + \frac{t}{n}\right),$$

$$\tilde{U}_n(\delta;t) := \tilde{U}_{n,t}(\delta) := \tilde{K}\left(\phi_n(\tilde{c}\delta) + D(\tilde{c}\delta)\sqrt{\frac{t}{n}} + \frac{t}{n}\right),$$

where $2 \leq \hat{K} \leq \tilde{K}$, $\hat{c},\tilde{c} \geq 1$ are numerical constants. It happens that $\hat{U}_n$ is a data-dependent function that upper bounds $\bar{U}_n$ with a high probability. $\tilde{U}_n$ is a distribution-dependent function that provides an upper bound on $\hat{U}_n$ (again, with a high probability). We now construct $\bar{V}_n, \hat{V}_n, \tilde{V}_n$ from $\bar{U}_n, \hat{U}_n, \tilde{U}_n$ the same way as we have constructed $V_n$ from $U_n$ and set $\hat{\delta}_n(t) := \hat{U}_{n,t}^{\sharp,q}(1/2q^3)$, $\tilde{\delta}_n(t) := \tilde{U}_{n,t}^{\sharp,q}(1/2q^3)$.

We will prove the following theorem.

THEOREM 3. *For all $t > 0$*

$$\mathbb{P}\{\bar{\delta}_n(t) \leq \hat{\delta}_n(t) \leq \tilde{\delta}_n(t)\} \geq 1 - \left(\log_q \frac{q^2}{\delta_n(t)} + 4\log_q \frac{q}{\delta_n(t)}\right)\exp\{-t\}.$$

In many situations, $\delta_n(t)$ and $\tilde{\delta}_n(t)$ are asymptotically within a constant one from another as $n \to \infty$. The above theorem suggests that $\hat{\delta}_n(t)$ can be used as an estimate (up to a constant) of $\delta_n(t)$ and this allows one to use this quantity as a data-dependent penalty in a model selection setting.

**4. Toward sharper inequalities for excess risk.** Suppose that risk minimization problem (1.1) has multiple solutions. This is a possibility, for instance, in risk minimization with nonconvex loss functions. Also, in a model selection framework (see Section 5) one deals with a family of risk minimization problems over classes $\mathcal{F}_k \subset \mathcal{F}$ that approximate problem (1.1). It is possible then that the global minimum of risk over the class $\mathcal{F}$ is attained at a number of different competing classes (models) $\mathcal{F}_k$. Anyway, the multiple minima case has to be understood as a part of comprehensive theory of empirical risk minimization. In such cases, the diameter $D(\delta) = D_P(\mathcal{F};\delta)$ of the $\delta$-minimal set does not tend to 0 as $\delta \to 0$, and it is easy to see that the quantity $\delta_n(t)$ defined in the previous section is going to be at least as large as $O(n^{-1/2})$. As a result, the bounds we have proved so far are not necessarily optimal. The question is whether it is possible to replace the diameter $D(\delta)$ by a more sophisticated geometric characteristic that would allow us to construct tighter bounds on the excess risk. We explore in this section one possible approach to this problem. Namely, we define the following quantity:

$$\check{r}(\sigma;\delta) := \sup_{f\in\mathcal{F}(\delta)} \inf_{g\in\mathcal{F}(\sigma)} \rho_P(f,g), \qquad 0 < \sigma \leq \delta,$$

that characterizes the accuracy of approximation of the functions from the $\delta$-minimal sets by the functions from the $\sigma$-minimal set for two different levels $\delta$ and $\sigma$. If $\mathcal{F}(0) \neq \emptyset$ (i.e., the minimum of $Pf$ is attained on $\mathcal{F}$), $\check{r}$ is also well defined for $\sigma = 0$, $\delta \geq \sigma$.

The function $\check{r}(\sigma, \delta)$ is nondecreasing in $\delta$, nonincreasing in $\sigma$ and $\check{r}(\delta, \delta) = 0$. If we extend $\check{r}$ to $\sigma > \delta$ by setting $\check{r}(\sigma; \delta) := \check{r}(\delta; \sigma)$, then, using the triangle inequality for $\rho_P$, it is easy to check that $\check{r}$ is a pseudometric. Clearly, $\check{r}(\sigma, \delta) \leq D(\delta)$. Moreover, it is not hard to imagine the situations when $\check{r}(0; \delta)$ is significantly smaller than $D(\delta)$ [say, $\check{r}(0; \delta) \to 0$ as $\delta \to 0$ whereas $D(\delta)$ is bounded away from 0]. Suppose, for instance, that $\mathcal{F} := \bigcup_j \mathcal{F}_j$, where $\mathcal{F}_j$ are classes of functions such that $\forall k, j : \min_{\mathcal{F}_j} Pf = \min_{\mathcal{F}_k} Pf$ (we assume that the minima are attained). Then it is easy to check that $\check{r}(0; \delta) \leq \sup_j D_P(\mathcal{F}_j; \delta)$. Of course, one can come up with examples of this sort in which $\check{r}(0, \delta) \to 0$ as $\delta \to 0$, but $D(\delta)$ is bounded away from 0.

It is not completely unnatural to expect that the function $\check{r}$ satisfies the condition of the following type:

(4.1) $$\check{r}(0; c_1 \delta) \leq c_2 \check{r}(0; \delta), \qquad \delta \in (0, 1]$$

for some constants $c_1, c_2 < 1$. Since $\check{r}(0; \delta) \leq \check{r}(0; c_1 \delta) + \check{r}(c_1 \delta, \delta)$, we get for all $\sigma \leq c_1 \delta$

$$\check{r}(\sigma; \delta) \leq \check{r}(0; \delta) \leq (1 - c_2)^{-1} \check{r}(\sigma; \delta),$$

which means that the values of $\check{r}(\sigma; \delta)$ are within a constant one from another for all $\sigma$ that are not too close to $\delta$ ($\sigma \leq c_1 \delta$).

Let

$$\check{\psi}_n(\sigma, \delta) := \lim_{\varepsilon \to 0} \mathbb{E} \sup_{g \in \mathcal{F}(\sigma)} \sup_{f \in \mathcal{F}(\delta), \rho_P(f,g) \leq \check{r}(\sigma, \delta) + \varepsilon} |(P_n - P)(f - g)|$$

and

$$\check{U}_n(\sigma; \delta; t) := \check{\psi}_n(\sigma, \delta) + \sqrt{2 \frac{t}{n} (\check{r}^2(\sigma, \delta) + 2\check{\psi}_n(\sigma, \delta))} + \frac{t}{2n}.$$

Almost as before, we will need

$$\check{V}_n(\sigma; \delta; t) := \sup_{j : \delta_j \geq \delta} \frac{\check{U}_n(\sigma; \delta_j; t) + \sigma}{\delta_j}.$$

Finally, we define $\check{\delta}_n(\sigma; t) := \inf\{\delta : \check{V}_n(\sigma; \delta; t) \leq 1/2q\}$. Clearly, $\check{\delta}_n(\sigma; t)$ is the $\sharp, q$-transform of the function $\delta \mapsto \check{U}_n(\sigma; \delta; t) + \sigma$ computed at the point $1/2q$. We obtain the following version of Theorem 1.

THEOREM 4. *For all $\sigma \in (0, 1]$, all $t > 0$ and all $\delta \geq \check{\delta}_n(\sigma; t)$,*

$$\mathbb{P}\{\mathcal{E}(\hat{f}_n) \geq \delta\} \leq \log_q \frac{q}{\delta} \exp\{-t\}$$

*and*

$$\mathbb{P}\left\{\exists f \in \mathcal{F}: \mathcal{E}(f) \geq \delta \text{ and } \frac{\hat{\mathcal{E}}_n(f)}{\mathcal{E}(f)} \leq 1 - q\check{V}_n(\sigma;\delta;t)\right\} \leq \log_q \frac{q}{\delta}\exp\{-t\}.$$

Note that, unlike the inequalities of Theorem 1, we have here only a *one-sided* bound for the ratio $\frac{\hat{\mathcal{E}}_n(f)}{\mathcal{E}(f)}$. As a result, it is easy to show that, for all $\sigma \in (0,1]$ and all $t > 0$, there exists an event of probability at least $1 - \log_q \frac{q^2}{\check{\delta}_n(\sigma;t)} e^{-t}$ such that on this event $\forall \delta \geq \check{\delta}_n(\sigma,t)$ the inclusion $\hat{\mathcal{F}}_n(\delta) \subset \mathcal{F}(2\delta)$ holds, but not the other inclusion of Lemma 2. The following proposition shows that this difficulty is unavoidable and the set $\hat{\mathcal{F}}_n(\delta)$ does not include even $\mathcal{F}(0)$ for the values of $\delta$ of the order $\check{\delta}_n(\sigma;t)$, or even larger. Because of this reason, the estimation of the quantity $\check{r}(\sigma;\delta)$ based on the data is a much harder problem than the estimation of the diameter $D_P(\mathcal{F};\delta)$. The discussion of this problem goes beyond the scope of the paper.

PROPOSITION 2. *Let $S := \{0,1\}^{N+1}$ and $P$ be the uniform distribution on $\{0,1\}^{N+1}$. Let $\mathcal{F} := \{f_j : 1 \leq j \leq N+1\}$, where $f_j(x) = x_j$, $x = (x_1, \ldots, x_{N+1}) \in \{0,1\}^{N+1}$. Then the following statements hold:*

(i) *$\mathcal{E}_P(\hat{f}) = 0$;*
(ii) *with some $C > 0$, $\check{\delta}_n(\sigma;t) \leq Ct/n$;*
(iii) *with some $c > 0$, $\delta_n(t) \geq c((\log N/n)^{1/2} + (t/n)^{1/2})$;*
(iv) *for any $\varepsilon > 0$ there exists $N_0$ such that, for $N_0 \leq N \leq \sqrt{n}$ and for $\delta = 0.25(\log N/n)^{1/2}$, the inclusion $\mathcal{F}(0) \subset \hat{\mathcal{F}}_n(\delta)$ does not hold with probability at least $1 - \varepsilon$.*

**5. Model selection.** Consider a family of function classes $\{\mathcal{F}_k\}$ such that $\forall k, \mathcal{F}_k \subset \mathcal{F}$. In applications, the classes $\{\mathcal{F}_k\}$ are used to find an approximate solution of risk minimization problem on the bigger class $\mathcal{F}$ of functions of interest. Let $\hat{f}_k := \hat{f}_{n,k} := \operatorname{argmin}_{f \in \mathcal{F}_k} P_n f$ be the corresponding empirical risk minimizers (we assume for simplicity that they exist). The goal is to construct, based on $\{\hat{f}_{n,k}\}$, a function $\hat{f} \in \mathcal{F}$ for which the excess risk $\mathcal{E}_P(\mathcal{F}; \hat{f})$ is small. To formulate the problem more precisely, suppose that there exists an index $k(P)$ such that $\inf_{\mathcal{F}_{k(P)}} Pf = \inf_{\mathcal{F}} Pf$, that is, a risk minimizer over the large class $\mathcal{F}$ can be found in a smaller class $\mathcal{F}_{k(P)}$. Let $\tilde{\delta}_n(k)$ be an upper bound on the excess risk (with respect to the class $\mathcal{F}_k$) of $\hat{f}_{n,k}$ that provides the optimal (in a minimax sense), or just a desirable accuracy of the solution of empirical risk minimization problem on the class $\mathcal{F}_k$. If there were an oracle who could tell a statistician that $k(P) = k$ is the right index of the class to be used, then the risk minimization problem could be solved with the accuracy $\tilde{\delta}_n(k)$. The *model selection problem* deals with constructing a data-dependent index $\hat{k} = \hat{k}(X_1, \ldots, X_n)$ of the model such that the excess risk of $\hat{f} := \hat{f}_{n,\hat{k}}$ is within a constant from $\tilde{\delta}_n(k(P))$ with a high probability. More generally, in the case when the global minimum over $\mathcal{F}$ is not attained precisely in any of the classes $\mathcal{F}_k$, one can still hope to show that with a high probability

$$\mathcal{E}_P(\mathcal{F}; \hat{f}) \leq C \inf_k \Big[\inf_{\mathcal{F}_k} Pf - Pf_* + \tilde{\pi}_n(k)\Big],$$

where $f_* := \operatorname{argmin}_{f \in \mathcal{F}} Pf$ (its existence will be assumed in what follows), $\tilde{\pi}_n(k)$ are "ideal" distribution-dependent complexity penalties associated with risk minimization over $\mathcal{F}_k$ and $C$ is a constant (preferably, $C=1$ or at least close to 1). The inequalities that express such a property are often referred to as *oracle inequalities*.

Among the most popular approaches to model selection are *penalization methods*, in which $\hat{k}$ is defined as a solution of the following minimization problem:

(5.1) $$\hat{k} := \operatorname*{argmin}_{k \geq 1} \{P_n \hat{f}_k + \hat{\pi}(k)\},$$

where $\hat{\pi}(k)$ is a *complexity penalty* (generally, data dependent) associated with the class (the model) $\mathcal{F}_k$. In other words, instead of minimizing the empirical risk on the whole class $\mathcal{F}$ we now minimize a penalized empirical risk. We discuss below two penalization methods (one in spirit of [34], another one more in spirit of [36]) with the penalties based on data-dependent bounds on excess risk developed in previous sections. Penalization methods proved to be very useful in a variety of statistical problems, including nonparametric regression. However, there are substantial difficulties in implementing model selection techniques based on penalization in nonparametric classification problems. Up to our best knowledge, this approach has failed so far to produce adaptive classification rules with fast Tsybakov's-type convergence rates (an exception is the recent result by [45] that achieves this goal, but only in a very special and somewhat artificial framework). As an alternative, we discuss a general model selection technique based on comparing the minima of empirical risk for different models with certain data-dependent thresholds (defined in terms of excess risk confidence bounds of the previous sections) that allows one to recover Tsybakov's convergence rates in very general risk minimization problems, including classification (note that Tsybakov [44] also used a version of comparison method in a specialized classification framework).

To provide some motivation for the approaches discussed below, note that ideally one would want to find $\hat{k}$ by minimizing over $k$ the *global* excess risk $\mathcal{E}_P(\mathcal{F}; \hat{f}_{n,k})$ of the solutions. This is impossible without oracle's help, so one has to develop some data-dependent upper confidence bounds on the excess risk. The following trivial representation (that plays the role of "bias-variance decomposition")

$$\mathcal{E}_P(\mathcal{F}; \hat{f}_{n,k}) = \inf_{\mathcal{F}_k} Pf - Pf_* + \mathcal{E}_P(\mathcal{F}_k; \hat{f}_{n,k})$$

shows that part of the problem is to come up with data-dependent upper bounds on the *local* excess risk $\mathcal{E}_P(\mathcal{F}_k; \hat{f}_{n,k})$, which is precisely what we considered in the previous sections. Another part is to bound $\inf_{\mathcal{F}_k} Pf - Pf_*$ in terms of $\inf_{\mathcal{F}_k} P_n f - P_n f_*$, which is what we do in Lemma 4 below. Combining these two bounds provides an upper bound on the global excess risk that can be now minimized with respect to $k$ ($P_n f_*$ can be dropped since it does not depend on $k$). Another approach is to use the representation

$$\mathcal{E}_P(\mathcal{F}; \hat{f}_{n,k}) - \mathcal{E}_P(\mathcal{F}; \hat{f}_{n,l}) = \inf_{\mathcal{F}_k} Pf - \inf_{\mathcal{F}_l} Pf + \mathcal{E}_P(\mathcal{F}_k; \hat{f}_{n,k}) - \mathcal{E}_P(\mathcal{F}_l; \hat{f}_{n,l})$$

and data-dependent bounds on local excess risk to develop a model selection technique based on comparison of the difference between $\inf_{\mathcal{F}_k} P_n f$ and $\inf_{\mathcal{F}_l} P_n f$ with certain data-dependent thresholds (which is done in Section 5.3 below).

For $\mathcal{G} \subset \mathcal{F}$, the distribution-dependent complexity $\bar{\delta}_n(\mathcal{G}; t)$ is defined as in Section 3 $[\bar{\delta}_n(t) = \bar{U}_{n,t}^\sharp(1/2q^3)]$. Let $t_k \geq 0$ and let $\hat{\delta}_n(\mathcal{F}_k; t_k)$ and $\tilde{\delta}_n(\mathcal{F}_k; t_k)$ be, respectively, data-dependent and distribution-dependent complexities such that

(5.2) $\quad \forall k \quad \mathbb{P}\{\bar{\delta}_n(\mathcal{F}_k; t_k) \leq \hat{\delta}_n(\mathcal{F}_k; t_k) \leq \tilde{\delta}_n(\mathcal{F}_k; t_k)\} \geq 1 - p_k.$

In particular, one can use the version of these complexities constructed in Section 3, in which case $p_k := \log_q \frac{q^2 n}{t_k} e^{-t_k} + 4 \log_q \frac{qn}{t_k} e^{-t_k}$, by Theorem 3. We use these notations throughout the section.

5.1. *Penalization method*: *version* 1. Define the following penalties:

$$\hat{\pi}(k) := \hat{K}\left[\hat{\delta}_n(\mathcal{F}_k, t_k) + \sqrt{\frac{t_k}{n} \inf_{\mathcal{F}_k} P_n f} + \frac{t_k}{n}\right] \quad \text{and}$$

$$\tilde{\pi}(k) := \tilde{K}\left[\tilde{\delta}_n(\mathcal{F}_k, t_k) + \sqrt{\frac{t_k}{n} \inf_{\mathcal{F}_k} Pf} + \frac{t_k}{n}\right],$$

where $\hat{K}, \tilde{K}$ are sufficiently large numerical constants. Here $\tilde{\pi}(k)$ represents a "desirable accuracy" of risk minimization on the class $\mathcal{F}_k$. The index estimate $\hat{k}$ is defined according to standard penalization method (5.1) and we set $\hat{f} := \hat{f}_{n,\hat{k}}$.

THEOREM 5. *There exists a choice of $\hat{K}, \tilde{K}$ such that for any sequence $\{t_k\}$ of positive numbers,*

$$\mathbb{P}\Big\{P\hat{f} \geq \inf_{k\geq 1}\{P_n\hat{f}_{n,k} + \hat{\pi}(k)\}\Big\} \leq \sum_{k=1}^{\infty}\Big(p_k + \log_q \frac{q^3 n}{t_k}e^{-t_k}\Big)$$

*and*

$$\mathbb{P}\Big\{\mathcal{E}_P(\mathcal{F};\hat{f}) \geq \inf_{k\geq 1}\Big\{\inf_{f\in\mathcal{F}_k} Pf - \inf_{f\in\mathcal{F}} Pf + \tilde{\pi}(k)\Big\}\Big\} \leq \sum_{k=1}^{\infty}\Big(p_k + \log_q \frac{q^3 n}{t_k}e^{-t_k}\Big).$$

The first bound of the theorem is an upper confidence bound on the risk of $\hat{f}$ in terms of minimal penalized empirical risk. The second bound is an oracle inequality showing that the excess risk of the function $\hat{f}$ is nearly optimal (up to complexity penalty terms).

The proof relies on the following lemma, which might be of independent interest.

LEMMA 3. *Given a class $\mathcal{F}$ of measurable functions from $S$ into $[0,1]$, suppose that, for some $t > 0$ and $p \in (0,1)$, $\mathbb{P}\{\bar{\delta}_n(\mathcal{F};t) \leq \hat{\delta}_n(\mathcal{F};t)\} \geq 1 - p$. Then the following inequalities hold:*

$$\mathbb{P}\Bigg\{\Big|\inf_{\mathcal{F}} P_n f - \inf_{\mathcal{F}} Pf\Big| \geq 2\bar{\delta}_n(\mathcal{F};t) + \sqrt{\frac{2t}{n}\inf_{\mathcal{F}} Pf} + \frac{t}{n}\Bigg\} \leq \log_q \frac{q^3}{\bar{\delta}_n(t)}e^{-t}$$

*and*

$$\mathbb{P}\Bigg\{\Big|\inf_{\mathcal{F}} P_n f - \inf_{\mathcal{F}} Pf\Big| \geq 4\hat{\delta}_n(\mathcal{F};t) + 2\sqrt{\frac{2t}{n}\inf_{\mathcal{F}} P_n f} + \frac{8t}{n}\Bigg\} \leq p + \log_q \frac{q^3}{\bar{\delta}_n(t)}e^{-t}.$$

5.2. *Penalization method: version* 2. For this version of penalization method, the following assumption is crucial:

(5.3) $$\forall f \in \mathcal{F} \qquad Pf - Pf_* \geq \varphi(\sqrt{\mathrm{Var}_P(f - f_*)}),$$

where $\varphi$ is a convex nondecreasing function on $[0, +\infty)$ with $\varphi(0) = 0$. We also assume that $\varphi(uv) \leq \varphi(u)\varphi(v)$, $u, v \geq 0$. The function $\varphi$ is supposed to be known and is involved in the definition of the penalties. This is the case, for instance, in least squares regression where one can use $\varphi(u) = u^2/2$ (see Section 8). However, in classification problems $\varphi$ is typically unknown, but it has a significant impact on the convergence rates. Adapting to unknown function $\varphi$ is a challenge for model selection in classification setting.

Denote $\varphi^*(v) := \sup_{u\geq 0}[uv - \varphi(u)]$ the conjugate of $\varphi$. We have $uv \leq \varphi(u) + \varphi^*(v)$, $u, v \geq 0$. For a fixed $\varepsilon > 0$, define the penalties as follows:

$$\hat{\pi}(k) := A(\varepsilon)\hat{\delta}_n(\mathcal{F}_k;t_k) + \varphi^*\Bigg(\sqrt{\frac{2t_k}{\varepsilon n}}\Bigg) + \frac{t_k}{n}$$

and

$$\tilde{\pi}(k) := \frac{A(\varepsilon)}{1+\varphi(\sqrt{\varepsilon})}\tilde{\delta}_n(\mathcal{F}_k; t_k) + \frac{2}{1+\varphi(\sqrt{\varepsilon})}\varphi^*\left(\sqrt{\frac{2t_k}{\varepsilon n}}\right) + \frac{2}{(1+\varphi(\sqrt{\varepsilon}))}\frac{t_k}{n},$$

where $A(\varepsilon) := \frac{5}{2} - \varphi(\sqrt{\varepsilon})$. As before, $\hat{k}$ is defined by (5.1) and $\hat{f} := \hat{f}_{n,\hat{k}}$.

THEOREM 6. *For any sequence $\{t_k\}$ of positive numbers,*

$$\mathbb{P}\Big\{\mathcal{E}_P(\mathcal{F};\hat{f}) \geq C(\varepsilon) \inf_{k \geq 1}\Big\{\inf_{f \in \mathcal{F}_k} Pf - \inf_{f \in \mathcal{F}} Pf + \tilde{\pi}(k)\Big\}\Big\}$$

$$\leq \sum_{k=1}^{\infty}\left(p_k + 2\log_q \frac{q^2 n}{t_k}e^{-t_k}\right),$$

where $C(\varepsilon) := \frac{1+\varphi(\sqrt{\varepsilon})}{1-\varphi(\sqrt{\varepsilon})}$.

The following lemma is used in the proof.

LEMMA 4. *Let $\mathcal{G} \subset \mathcal{F}$. For all $t > 0$, there exists an event $E$ with probability at least $1 - \log_q \frac{q^3 n}{t}e^{-t}$ such that on this event*

$$(5.4) \quad \inf_{\mathcal{G}} P_n f - P_n f_* \leq (1+\varphi(\sqrt{\varepsilon}))\Big(\inf_{\mathcal{G}} Pf - Pf_*\Big) + \varphi^*\left(\sqrt{\frac{2t}{\varepsilon n}}\right) + \frac{t}{n}$$

*and*

$$(5.5) \quad \begin{aligned}\inf_{\mathcal{G}} Pf - Pf_* &\leq (1-\varphi(\sqrt{\varepsilon}))^{-1} \\ &\times \left[\inf_{\mathcal{G}} P_n f - P_n f_* + \frac{3}{2}\bar{\delta}_n(\mathcal{G}; t) + \varphi^*\left(\sqrt{\frac{2t}{\varepsilon n}}\right) + \frac{t}{n}\right].\end{aligned}$$

*In addition, if there exists $\bar{\delta}_n(\mathcal{G}; \varepsilon; t)$ such that*

$$\bar{\delta}_n(\mathcal{G}; t) \leq \varepsilon\Big(\inf_{\mathcal{G}} Pf - Pf_*\Big) + \bar{\delta}_n(\mathcal{G}; \varepsilon; t),$$

*then*

$$(5.6) \quad \begin{aligned}\inf_{\mathcal{G}} Pf - Pf_* &\leq \left(1 - \varphi(\sqrt{\varepsilon}) - \frac{3}{2}\varepsilon\right)^{-1} \\ &\times \left[\inf_{\mathcal{G}} P_n f - P_n f_* + \frac{3}{2}\bar{\delta}_n(\mathcal{G}; \varepsilon; t) + \varphi^*\left(\sqrt{\frac{2t}{\varepsilon n}}\right) + \frac{t}{n}\right].\end{aligned}$$

REMARKS. 1. It is easily seen from the proofs that the same inequality holds for arbitrary penalties $\hat{\pi}(k)$ and $\tilde{\pi}(k)$ such that with probability at least $1 - p_k$

$$\hat{\pi}(k) \geq A(\varepsilon)\bar{\delta}_n(\mathcal{F}_k; t_k) + \varphi^*\left(\sqrt{\frac{2t_k}{\varepsilon n}}\right) + \frac{t_k}{n}$$

and

$$\tilde{\pi}(k) \geq \frac{\hat{\pi}(k)}{1+\varphi(\sqrt{\varepsilon})} + \frac{\varphi^*(\sqrt{\frac{2t_k}{\varepsilon n}})}{1+\varphi(\sqrt{\varepsilon})} + \frac{t_k}{(1+\varphi(\sqrt{\varepsilon}))n}.$$

2. Suppose that the following condition holds:

$$\bar{\delta}_n(\mathcal{F}_k; t) \leq \varepsilon \Big(\inf_{\mathcal{F}_k} Pf - Pf_*\Big) + \bar{\delta}_n(\mathcal{F}_k; \varepsilon; t),$$

as is the case in Lemma 5 below. Suppose also that there exist $\hat{\delta}_n(\mathcal{F}_k; \varepsilon; t_k)$, $\tilde{\delta}_n(\mathcal{F}_k; \varepsilon; t_k)$ such that

$$\forall k \qquad \mathbb{P}\{\bar{\delta}_n(\mathcal{F}_k; \varepsilon; t_k) \leq \hat{\delta}_n(\mathcal{F}_k; \varepsilon; t_k) \leq \tilde{\delta}_n(\mathcal{F}_k; \varepsilon; t_k)\} \geq 1 - p_k.$$

Then, using the bound (5.6) of Lemma 4, one can easily modify Theorem 6 replacing in the definition of the penalties the quantities $\bar{\delta}_n(\mathcal{F}_k; t_k)$, $\hat{\delta}_n(\mathcal{F}_k; t_k)$, $\tilde{\delta}_n(\mathcal{F}_k; t_k)$, by $\bar{\delta}_n(\mathcal{F}_k; \varepsilon; t_k)$, $\hat{\delta}_n(\mathcal{F}_k; \varepsilon; t_k)$, $\tilde{\delta}_n(\mathcal{F}_k; \varepsilon; t_k)$ and also defining

$$A(\varepsilon) := \tfrac{3}{2} + (1 - \varphi(\sqrt{\varepsilon}) - \tfrac{3}{2}\varepsilon)/(1+\varepsilon) \quad \text{and}$$
$$C(\varepsilon) := (1 + \varphi(\sqrt{\varepsilon}))(1+\varepsilon)/(1 - \varphi(\sqrt{\varepsilon}) - \tfrac{3}{2}\varepsilon).$$

3. Note also that if $\bar{\delta}_n(\mathcal{F}_k; t_k)$ is replaced by $\check{\delta}_n(\mathcal{F}_k; t_k)$, defined as in Theorem 2, the result of Theorem 6 is also true, and, moreover, the logarithmic factor in the oracle inequality can be dropped: the expression in the right-hand side of the bound of Theorem 6 becomes $\sum_{k=1}^{\infty}(p_k + 4e^{-t_k})$.

4. The result also holds if condition (5.3) holds for each $k$ and for all $f \in \mathcal{F}_k$ with its own function $\varphi_k$ (but with the same function $f_*$) and the sequence of functions $\{\varphi_k\}$ is nonincreasing: $\forall k \ \varphi_k \geq \varphi_{k+1}$. In this case, one should use the function $\varphi_k$ in the definitions of $\hat{\pi}(k), \tilde{\pi}(k)$. $C(\varepsilon)$ is defined as before with $\varphi = \varphi_1$.

5.3. *Comparison method.* The version of comparison method presented here relies on the following assumption: $\mathcal{F}_1 \subset \mathcal{F}_2 \subset \cdots$. Denote

$$\bar{\delta}_n(k) := \max_{1 \leq j \leq k} \bar{\delta}_n(\mathcal{F}_j; t_j), \qquad \hat{\delta}_n(k) := \max_{1 \leq j \leq k} \hat{\delta}_n(\mathcal{F}_j; t_j),$$

$$\tilde{\delta}_n(k) := \max_{1 \leq j \leq k} \tilde{\delta}_n(\mathcal{F}_j; t_j)$$

and define with some numerical constants $\bar{c}, \hat{c}, \tilde{c}$ and with inf being $\infty$ if the set of $k$'s is empty:

$$k^* := k^*(P) := \inf\Big\{k : \forall l > k \inf_{\mathcal{F}_k} Pf = \inf_{\mathcal{F}_l} Pf\Big\},$$

$$\bar{k} := \bar{k}(P) := \inf\Big\{k : \forall l > k \inf_{\mathcal{F}_k} Pf - \inf_{\mathcal{F}_l} Pf \leq \bar{c}\bar{\delta}_n(l)\Big\},$$

$$\hat{k} := \inf\Big\{k : \forall l > k \inf_{\mathcal{F}_k} P_n f - \inf_{\mathcal{F}_l} P_n f \leq \hat{c}\hat{\delta}_n(l)\Big\},$$

$$\tilde{k} := \tilde{k}(P) := \inf\Big\{k : \forall l > k \inf_{\mathcal{F}_k} Pf - \inf_{\mathcal{F}_l} Pf \leq \tilde{c}\tilde{\delta}_n(l)\Big\}.$$

Finally, let $\hat{f} := \hat{f}_{n,\hat{k}}$ (if $\hat{k} = \infty$, $\hat{f}$ can be defined in an arbitrary way, say, $\hat{f} = \hat{f}_{n,1}$).

THEOREM 7. *There exists a choice of constants $\bar{c}, \hat{c}, \tilde{c}$ such that with some constant $C > 0$ for any sequence $\{t_k\}$, $t_k > 0$*

$$\mathbb{P}\Big\{P\hat{f} - \inf_k \inf_{\mathcal{F}_k} Pf \geq \inf_{k \geq \bar{k}(P)}\Big[\inf_{\mathcal{F}_k} Pf - \inf_k \inf_{\mathcal{F}_k} Pf + C\tilde{\delta}_n(k)\Big]\Big\}$$
$$\leq \sum_{k=1}^{\infty}\Big(p_k + \log_q \frac{q^2 n}{t_k} e^{-t_k}\Big).$$

*In particular, if $k^*(P) < \infty$, then*

$$\mathbb{P}\Big\{P\hat{f} - \inf_k \inf_{\mathcal{F}_k} Pf \geq C\tilde{\delta}_n(k^*(P))\Big\} \leq \sum_{k=1}^{\infty}\Big(p_k + \log_q \frac{q^2 n}{t_k} e^{-t_k}\Big).$$

REMARKS. 1. If $\bar{k}(P) = \infty$, assume that the infimum over $k \geq \bar{k}(P)$ is equal to 1, which makes the first bound trivial. If $\bar{k}(P) < \infty$, it follows from the proof that so is $\hat{k}$ (with an exception of the event whose probability is controlled in the theorem).

2. If $\bar{\delta}_n(\mathcal{F}_k; t_k)$ is replaced by $\check{\delta}_n(\mathcal{F}_k; t_k)$ (as defined in Theorem 2), then the logarithmic factor in the oracle inequality can be dropped and the expression in the right-hand side of the bounds becomes $\sum_{k=1}^{\infty}(p_k + 2e^{-t_k})$.

**6. Connection to several recent results.** In this section, we discuss the connection of our main results to some other recent work on model selection in risk minimization problems, including [34, 36, 44].

6.1. *Tsybakov.* Our first example is motivated by the recent work of Tsybakov [44] (see also the earlier paper by Mammen and Tsybakov [35]), on fast convergence rates in classification. Let $\rho_P^2(f,g) := P(f-g)^2$. Define the expected continuity modulus $\omega_n(\mathcal{F}; \delta)$ as in Section 3. For $\rho \in (0,1)$, $\kappa \geq 1$ and $C > 0$, let $\mathcal{P}_{\rho,\kappa,C}(\mathcal{F})$ denote the class of probability measures $P$ such that the following two conditions hold:

(i) $\omega_n(\mathcal{F}; \delta) \leq C\delta^{1-\rho} n^{-1/2}$;
(ii) $D_P(\mathcal{F}; \delta) \leq C\delta^{\frac{1}{2\kappa}}$.

THEOREM 8. *Under conditions* (i) *and* (ii), $\sup_{P \in \mathcal{P}_{\rho,\kappa,C}(\mathcal{F})} \mathbb{E}\mathcal{E}_P(\mathcal{F}; \hat{f}_n) = O(n^{-\frac{\kappa}{2\kappa+\rho-1}})$.

This result generalizes Theorem 1 in [44]. Namely, using the standard Dudley's entropy integral bound on the expected continuity modulus of the empirical process under the condition that the $L_2(P)$-entropy with bracketing of the class $\mathcal{F}$ grows as $O(\varepsilon^{-2\rho})$ (see, e.g., [47], Theorem 2.14.2) yields condition (i). If

(6.1) $\quad f_* := f_{*,P} := \underset{f \in \mathcal{F}}{\operatorname{argmin}} Pf \quad \text{and} \quad Pf - Pf_* \geq c_0 \rho_P^{2\kappa}(f, f_*),$

then also condition (ii) is satisfied. The conditions above, being translated to the case of classes of sets (which was the case considered by Tsybakov whose paper dealt with the binary classification problem), are precisely the assumptions (A1) and (A2) in Tsybakov [44] and the rate of convergence $(n^{-\frac{\kappa}{2\kappa+\rho-1}})$ is the one obtained by Tsybakov. Of course, condition (i) will be also satisfied under many other assumptions common in empirical processes theory; for example, it can be expressed in terms of random entropies of the class. Also, the diameter $D_P(\mathcal{F}; \delta)$ in condition (ii) can be replaced by a more subtle geometric characteristic $\check{r}(0; \delta) = \check{r}_P(\mathcal{F}; 0, \delta)$ defined in Section 4. In other words, condition (6.1) can be replaced by the following:

(6.2) $\quad \forall f \in \mathcal{F} \; \exists f_* \in \underset{f \in \mathcal{F}}{\operatorname{argmin}} Pf = \mathcal{F}(0): \qquad Pf - Pf_* \geq c_0 \rho_P^{2\kappa}(f, f_*),$

including the case when the risk $Pf$ has multiple minima on $\mathcal{F}$. Theorem 8 holds in this case with only minor changes in the proof.

Next we turn to model selection.

THEOREM 9. *Consider a family* $\{(\mathcal{F}_j, \mathcal{P}_j)\}_{1 \leq j \leq N}$, *such that* $\mathcal{F}_j \subset \mathcal{F}$, $\mathcal{P}_j := \mathcal{P}_{\rho_j,\kappa_j,C}(\mathcal{F}_j)$ *and for all* $P \in \mathcal{P}_j$ *we have* $f_{*,P} \in \mathcal{F}_j$. *Moreover, assume that* $\mathcal{F}_1 \subset \mathcal{F}_2 \subset \cdots \subset \mathcal{F}_N$, *that for all* $P \in \mathcal{P}_j$, $k^*(P) = j$ *(with* $k^*(P)$ *defined in Section* 5.3*) and that the numbers* $\beta_j := \kappa_j/(2\kappa_j + \rho_j - 1)$ *satisfy the condition* $\beta_1 \geq \beta_2 \geq \cdots \geq \beta_N$. *Define* $\hat{k}$ *and* $\hat{f}$ *as in Theorem* 7 *(with* $t_k := \log N + 3 \log n$, $k = 1, \ldots, n$). *Then*

$$\max_{1 \leq j \leq N} \sup_{P \in \mathcal{P}_j} n^{\beta_j} \mathbb{E}(P\hat{f} - Pf_*) = O(1) \qquad \textit{as } n \to \infty.$$

Note that the result is also true if $N = N_n$, where $N_n$ grows not too fast, say, so that for all $\delta > 0$, $\log N_n = o(n^\delta)$ as $n \to \infty$. This should be compared with Theorem 3 in [44] where another method of constructing an adaptive empirical risk minimizer was suggested in a more special classification framework and it was proved that the optimal convergence rate is attained at this estimate up to a logarithmic factor. Our Theorem 9 extends these types of result to a more general framework of abstract empirical risk minimization and refines them by removing the logarithmic factor.

6.2. *Lugosi and Wegkamp.* Next we turn to the results of a recent paper of Lugosi and Wegkamp [34]. Suppose that $\mathcal{F}$ is a class of measurable functions on $S$ taking values in $\{0,1\}$ (binary functions). As in Section 2, Example 6, $\Delta^{\mathcal{F}}(X_1,\ldots,X_n)$ denotes the shattering number of the class $\mathcal{F}$ on the sample $(X_1,\ldots,X_n)$.

Given a sequence $\{\mathcal{F}_k\}$, $\mathcal{F}_k \subset \mathcal{F}$, of classes of binary functions, define the penalties

$$\hat{\pi}(k) := \hat{K}\left[\sqrt{\inf_{f\in\mathcal{F}_k} P_n f \frac{\log \Delta^{\mathcal{F}_k}(X_1,\ldots,X_n)+t_k}{n}} + \frac{\log \Delta^{\mathcal{F}_k}(X_1,\ldots,X_n)+t_k}{n}\right]$$

and

$$\tilde{\pi}(k) := \tilde{K}\left[\sqrt{\inf_{f\in\mathcal{F}_k} Pf \frac{\mathbb{E}\log \Delta^{\mathcal{F}_k}(X_1,\ldots,X_n)+t_k}{n}} + \frac{\mathbb{E}\log \Delta^{\mathcal{F}_k}(X_1,\ldots,X_n)+t_k}{n}\right],$$

and let $\hat{k}$ solve the penalized empirical risk minimization problem (5.1), $\hat{f} := \hat{f}_{n,\hat{k}}$.

THEOREM 10. *There exists a choice of $\hat{K}, \tilde{K}$ such that for all $t_k > 0$,*

$$\mathbb{P}\Big\{\mathcal{E}_P(\mathcal{F};\hat{f}) \geq \inf_{k\geq 1}\Big\{\inf_{f\in\mathcal{F}_k} Pf - \inf_{f\in\mathcal{F}} Pf + \tilde{\pi}(k)\Big\}\Big\} \leq 2\sum_{k=1}^{\infty} \log_q \frac{q^4 n}{t_k} e^{-t_k}.$$

The development of penalization techniques that lead to these types of oracle inequalities was one of the major goals of the paper of Lugosi and Wegkamp [34]. A little bit sharper results obtained in this paper (involving the shattering numbers or Rademacher complexities of the classes $\hat{\mathcal{F}}_k(\delta_k)$ for suitably chosen $\delta_k$ instead of the global shattering numbers) can be also recovered from Theorem 7 relatively easily (using Lemma 2).

6.3. *Massart.* We consider now some recent results of Massart [36] that we formulate in a somewhat different form. Suppose that $\mathcal{F}$ is a class of measurable functions from $S$ into $[0,1]$ and $f_*:S \mapsto [0,1]$ is a measurable function such that with some numerical constant $D > 0$

(6.3) $\quad D(Pf - Pf_*) \geq \rho_P^2(f, f_*) \geq P(f - f_*)^2 - (P(f - f_*))^2,$

where $\rho_P$ is a (pseudo)metric. We will assume, for simplicity, that the infimum of $Pf$ over $\mathcal{F}$ is attained at a function $\bar{f} \in \mathcal{F}$ (the result can be easily modified if this is not the case). Recall the definition of $\theta_n(\delta)$ in Section 2. The following lemma will be crucial.

LEMMA 5. *There exists a large enough numerical constant $K > 0$ such that for all $\varepsilon \in (0,1]$ and for all $t > 0$*

$$\bar{\delta}_n(\mathcal{F}; t) \leq \varepsilon \Big(\inf_{\mathcal{F}} Pf - Pf_*\Big) + \frac{1}{D}\theta_n^\sharp\Big(\frac{\varepsilon}{KD}\Big) + \frac{KD}{\varepsilon}\frac{t}{n}.$$

It immediately follows from the lemma and Theorem 1 that

$$\mathbb{P}\Big\{P\hat{f} - Pf_* \geq (1+\varepsilon)\Big(\inf_{\mathcal{F}} Pf - Pf_*\Big) + \frac{1}{D}\theta_n^\sharp\Big(\frac{\varepsilon}{KD}\Big) + \frac{KD}{\varepsilon}\frac{t}{n}\Big\} \leq \log_q \frac{qn}{t} e^{-t}$$

(and, due to Theorem 2, a version without the logarithmic factor holds with $\theta_n$ replaced by an upper bound $\check{\theta}_n$ of strictly concave type).

Now suppose that $\{\mathcal{F}_j\}$ is a sequence of function classes such that condition (6.3) holds for each class $\mathcal{F}_j$ with some constant $D_j \geq 1$ (and with the same $f_*$). Assume also that the sequence $\{D_j\}$ is nondecreasing. We denote $\bar{\delta}_n(\varepsilon; j) := D_j^{-1}\theta_n^\sharp(\varepsilon/KD_j)$ and suppose that for any $j$ there exist a data-dependent quantity $\hat{\delta}_n(\varepsilon; j)$ and a distribution-dependent quantity $\tilde{\delta}_n(\varepsilon; j)$ such that $\forall j, \mathbb{P}\{\bar{\delta}_n(\varepsilon; j) \leq \hat{\delta}_n(\varepsilon; j) \leq \tilde{\delta}_n(\varepsilon; j)\} \geq 1 - p_j$. Now we define the penalties as follows:

$$\hat{\pi}(\varepsilon; j) := 3\hat{\delta}_n(\varepsilon; j) + \frac{\hat{K}D_j t_j}{\varepsilon n} \quad \text{and} \quad \tilde{\pi}(\varepsilon; j) := 3\tilde{\delta}_n(\varepsilon; j) + \frac{\tilde{K}D_j t_j}{\varepsilon n}$$

with some numerical constants $\hat{K}, \tilde{K}$. Define $\hat{k}$ according to (5.1), $\hat{f} := \hat{f}_{\hat{k}}$.

The next result follows from Lemma 5 and Theorem 6.

THEOREM 11. *There exist numerical constants $\hat{K}, \tilde{K}$ such that for any sequence $\{t_k\}$ of positive numbers,*

$$\mathbb{P}\Big\{P\hat{f} - Pf_* \geq \frac{1+\varepsilon}{1-\varepsilon} \inf_{k \geq 1}\Big\{\inf_{f \in \mathcal{F}_k} Pf - Pf_* + \tilde{\pi}(\varepsilon; k)\Big\}\Big\}$$
$$\leq \sum_{k=1}^\infty \Big(p_k + 2\log_q \frac{q^2 n}{t_k} e^{-t_k}\Big).$$

*If, in addition, $\forall j, \forall \delta > 0 : \theta_n(\mathcal{F}_j; \delta) \leq \check{\theta}_n(\mathcal{F}_j; \delta)$, where $\check{\theta}_n(\mathcal{F}_j; \cdot) = \check{\theta}_{n,\mathcal{F}_j}(\cdot)$ is a function of strictly concave type, then one can replace $\bar{\delta}_n(\varepsilon; j)$ by $\check{\delta}_n(\varepsilon; j) := D_j^{-1}\check{\theta}_{n,\mathcal{F}_j}^\sharp(\varepsilon/KD_j)$, the right-hand side of the bound being in this case $\sum_{k=1}^\infty (p_k + 4e^{-t_k})$.*

This result has a number of applications. In a sense, most of the important complexity penalties used in learning theory can be derived as its consequence. For example (pointed out already in [36]), if $\mathcal{F}_k$ are classes of binary functions and

$$\hat{\pi}(k) := \frac{6 \log \Delta^{\mathcal{F}_k}(X_1, \ldots, X_n) + K t_k}{n},$$

one can use Theorem 11, the bounds of Example 6, Section 2 and the deviation inequalities for shattering numbers [12] to get very easily the following oracle inequality:

$$\mathbb{P}\bigg\{ P\hat{f} - Pf_* \geq C \inf_{k \geq 1} \bigg\{ \inf_{f \in \mathcal{F}_k} Pf - Pf_* + \frac{\mathbb{E} \log \Delta^{\mathcal{F}_k}(X_1, \ldots, X_n) + t_k}{n} \bigg\} \bigg\}$$
$$\leq 5 \sum_{k=1}^{\infty} e^{-t_k},$$

with some constant $C > 1$. One can also combine Theorem 11 with Lemma 1 to obtain oracle inequalities for penalization method based on localized Rademacher complexities (defined in terms of continuity modulus of Rademacher process).

**7. Loss functions and empirical risk minimization.** Let $T$ be a measurable space with $\sigma$-algebra $\mathcal{T}$, and let $(X, Y)$ be a random couple in $S \times T$ with joint distribution $P$. The distribution of $X$ will be denoted $\Pi$. Consider a sample $(X_1, Y_1), \ldots, (X_n, Y_n)$ of independent copies of $(X, Y)$ and let $P_n$ be the empirical distribution in $S \times T$ based on this sample, while $\Pi_n$ will denote the empirical distribution in $S$ based on the sample $(X_1, \ldots, X_n)$. Let $\ell: T \times \mathbb{R} \mapsto \mathbb{R}_+$ be a loss function. Given a class $\mathcal{G}$ of measurable functions from $S$ into $\mathbb{R}$, consider the following risk minimization problem:

$$\mathbb{E}\ell(Y, g(X)) \to \min, \qquad g \in \mathcal{G}.$$

If we denote $(\ell \bullet g)(x, y) := \ell(y; g(x))$, then we can rewrite this problem as $P(\ell \bullet g) \to \min$, $g \in \mathcal{G}$, or

$$Pf \to \min, \qquad f \in \mathcal{F} := \ell \bullet \mathcal{G} := \{\ell \bullet g : g \in \mathcal{G}\},$$

so we are dealing with problem (1.1) for a class $\mathcal{F}$ of special structure (the "loss class") and the results of previous sections can be specialized in this case.

Let $\mu_x$ denote a version of conditional distribution of $Y$ given $X = x$. Then the following representation of the risk holds under some mild regularity assumptions:

$$P(\ell \bullet g) = \int_S \int_T \ell(y; g(x)) \mu_x(dy) \Pi(dx).$$

Given a probability measure $\mu$ on $(T, \mathcal{T})$, let $u_\mu \in \operatorname{argmin}_{u \in \mathbb{R}} \int_T \ell(y; u)\mu(dy)$. If

$$g_*(x) := u_{\mu_x} = \underset{u \in \mathbb{R}}{\operatorname{argmin}} \int_T \ell(y; u)\mu_x(dy),$$

then we have (assuming that the function $g_*$ is well defined and measurable) $\forall g, P(\ell \bullet g) \geq P(\ell \bullet g_*)$, so $g_*$ is a *global* minimal point of $P(\ell \bullet g)$.

The corresponding empirical risk minimization problem is

$$P_n(\ell \bullet g) = n^{-1} \sum_{j=1}^n \ell(Y_j; g(X_j)) \to \min, \qquad g \in \mathcal{G},$$

and $\hat{g}_n$ will denote its solution (we assume its existence for simplicity). The following assumption on the loss function $\ell$ is very useful in the analysis of this problem. Suppose there exists a function $D(u, \mu) \geq 0$ such that for any measure $\mu = \mu_x$, $x \in S$

(7.1) $\quad \int_T (\ell(y, u) - \ell(y, u_\mu))^2 \mu(dy) \leq D(u, \mu) \int_T (\ell(y, u) - \ell(y, u_\mu))\mu(dy).$

In the case when the functions in the class $\mathcal{G}$ take their values in the interval $[-M/2, M/2]$ and $D(u, \mu_x)$, $|u| \leq M/2$, $x \in S$ is uniformly bounded by a constant $D > 0$, it immediately follows from (7.1) [by plugging in $u = g(x)$, $\mu = \mu_x$ and integrating with respect to $\Pi(dx)$] that for all $g \in \mathcal{G}$

(7.2) $\qquad\qquad P(\ell \bullet g - \ell \bullet g_*)^2 \leq DP(\ell \bullet g - \ell \bullet g_*).$

As a result, if $g_* \in \mathcal{G}$, then the $L_2(P)$-diameter of the $\delta$-minimal set of $\mathcal{F}$, $D(\mathcal{F}; \delta) \leq 2(D\delta)^{1/2}$. Moreover, even if $g_* \notin \mathcal{G}$, the condition (6.3) still holds for the loss class $\mathcal{F}$ with $f_* = \ell \bullet g_*$, opening the way for Massart's penalization method in these types of problems. The idea to control variance in terms of expectation has been extensively used in [36] (and even in earlier work of Birgé and Massart) and in learning theory literature [5, 6, 7, 8, 10, 37].

The analysis of risk minimization problems (in particular, proving the existence of $g_*$, checking condition (7.1), etc.) becomes much simpler under the convexity of the loss, that is, when for all $y \in T$, $\ell(y, \cdot)$ is a convex function. The problems of this type are called *convex risk minimization*. Both the least squares regression and $L_1$-regression as well as some of the methods of large margin classification (such as boosting) can be viewed as versions of convex risk minimization.

Assuming again that the functions in $\mathcal{G}$ take values in $[-M/2, M/2]$, we will introduce some even stricter assumptions on the loss function $\ell$. Namely, assume that $\ell$ satisfies the Lipschitz condition with some $L > 0$:

(7.3) $\quad \forall y \in T, \ \forall u, v \in [-M/2, M/2] \qquad |\ell(y, u) - \ell(y, v)| \leq L|u - v|$

and also that the following assumption on convexity modulus of $\ell$ holds with some $\Lambda > 0$:

(7.4) $\quad \forall y \in T, \ \forall u, \ v \in [-M/2, M/2] \quad \dfrac{\ell(y,u) + \ell(y,v)}{2} - \ell\left(y; \dfrac{u+v}{2}\right) \geq \Lambda |u-v|^2.$

Note that if $g_*$ is bounded by $M/2$, conditions (7.3) and (7.4) imply (7.1) with $D(u,\mu) \leq \dfrac{L^2}{2\Lambda}$. To see this, it is enough to use (7.4) with $v = u_\mu$, $\mu = \mu_x$ and integrate it with respect to $\mu$ to get for $L(u) := \int_T \ell(y,u)\mu(dy)$ (the minimum of $L$ is at $u_\mu$):

$$\frac{L(u) - L(u_\mu)}{2} = \frac{L(u) + L(u_\mu)}{2} - L(u_\mu)$$
$$\geq \frac{L(u) + L(u_\mu)}{2} - L\left(\frac{u + u_\mu}{2}\right) \geq \Lambda |u - u_\mu|^2$$

and then to use the Lipschitz condition to get

$$\int_T |\ell(y,u) - \ell(y,u_\mu)|^2 \mu(dy) \leq L^2 |u - u_\mu|^2.$$

This nice and simple trick, based on strict convexity, has been used repeatedly in the theory (see, e.g., [6]). We will use it again in the proof of Lemma 6. Sometimes a more general version of condition (7.4) is needed. It can be formulated as follows:

(7.5) $\quad \forall y \in T, \ \forall u, v \in [-M/2, M/2] \quad \dfrac{\ell(y,u) + \ell(y,v)}{2} - \ell\left(y; \dfrac{u+v}{2}\right) \geq \psi(|u-v|^r),$

where $\psi$ is a convex nondecreasing function and $r \in (0,2]$. The following lemma will allow us to bound the local complexities of the loss class $\mathcal{F} = \ell \bullet \mathcal{G}$ in terms of local complexities of the class $\mathcal{G}$, which is often needed in applications. Let

$$\bar{W}_n(\delta; t) = \bar{W}_{n,t}(\delta) := \bar{W}_n(\mathcal{G}; \delta; t)$$
$$:= C\left[L\theta_n(\mathcal{G}; \bar{g}; M^{2-r}\psi^{-1}(\delta/2)) + L\sqrt{\dfrac{M^{2-r}\psi^{-1}(\delta/2)(t+1)}{n}} + \dfrac{t}{n}\right],$$

where $C > 0$ is a numerical constant and $\theta_n$ is defined in Section 2.4.

LEMMA 6. *Suppose that $\mathcal{G}$ is a convex class of functions taking values in $[-M/2, M/2]$. Assume that the minimum of $P(\ell \bullet g)$ over $\mathcal{G}$ is attained at $\bar{g} \in \mathcal{G}$. Under the conditions (7.3) and (7.5), there is a choice of numerical constants $C$ and $\kappa_W$ such that $\forall \delta, t$, $\bar{U}_n(\mathcal{F}; \delta; t) \leq \bar{W}_n(\mathcal{G}; \delta; t)$ and $\bar{\delta}_n(\mathcal{F}; t) \leq \bar{\delta}_n^W(\mathcal{G}; t) := \bar{W}_{n,t}^{\sharp}(\kappa_W).$*

We are especially interested in the case when $\mathcal{G} := M \operatorname{conv}(\mathcal{H})$, where $\mathcal{H}$ is a base class of functions from $S$ into $[-1/2, 1/2]$ (see Example 5, Section 2.5). In this case, there are a number of powerful functional gradient descent-type algorithms (boosting algorithms) that allow one to implement convex empirical risk minimization over such classes. Assume that condition (2.1) holds for the class $\mathcal{H}$ with some $V > 0$. Define

$$\pi_n(M, L, \Lambda; t) := C\left[\Lambda M^{V/(V+1)}\left(\frac{L}{\Lambda} \vee 1\right)^{(V+2)/(V+1)} n^{-\frac{1}{2}\frac{V+2}{V+1}} + \frac{L^2}{\Lambda}\frac{t+1}{n}\right]$$

with some numerical constant $C$. The next result is essentially a slightly generalized version of a theorem due to Bartlett, Jordan and McAuliffe [6]. We will derive it as a corollary of our Theorem 2, using several nice observations of Bartlett, Jordan and McAuliffe [6] (contained in the proof of Lemma 6).

THEOREM 12. *Under the conditions* (7.3) *and* (7.4), $\bar{\delta}_n(\mathcal{F}; t) \leq \pi_n(M, L, \Lambda; t)$ *and as a result*

$$\mathbb{P}\left\{P(\ell \bullet \hat{g}_n) \geq \min_{g \in \mathcal{G}} P(\ell \bullet g) + \pi_n(M, L, \Lambda; t)\right\} \leq e^{-t}.$$

Because of the generality of the methods, the results can be easily extended to other examples of convex risk minimization problems. For instance, let $K$ be a symmetric nonnegatively definite kernel on $S \times S$ such that $|K(x,x)| \leq 1$ for all $x \in S$. As in Example 7, Section 2.5, $H_K$ is the reproducing kernel Hilbert space and $B_K$ is its unit ball. Let $\mathcal{G} := \mathcal{G}_M := \frac{M}{2} B_K$. This example is of importance in the theory of kernel machines. Clearly, $\mathcal{G}_M$ is a convex class of functions and, by elementary properties of reproducing kernel spaces, $\forall g \in \mathcal{G}_M, x \in S : |g(x)| \leq M/2$. We will use now slightly rescaled Mendelson's complexities of Example 8. It is easy to check (using Mendelson's inequalities of Example 8, Lemma 6 and the argument used at the beginning of the proof of Theorem 12) that

$$\bar{\delta}_n(\mathcal{F}; t) \leq \bar{\delta}_n^W(\mathcal{G}_M; t) \leq C\left[M^2 \Lambda \bar{\gamma}_n^\sharp\left(\frac{M\Lambda}{L}\right) + \frac{L^2}{\Lambda}\frac{t+1}{n}\right] =: \bar{\pi}_n(M, L, \Lambda, t).$$

With this new definition, the assertion of Theorem 12 still holds, and, moreover, based on the discussion in Example 7, one can replace in the bound the distribution-dependent Mendelson's complexity by its data-dependent version.

An alternative to the approach of Lemma 6, exploited, for instance, in the paper of Blanchard, Lugosi and Vayatis [10], is based on a straightforward comparison of $L_2(P_n)$-distances and the corresponding entropies for the classes $\mathcal{G}$ and $\mathcal{F} = \ell \bullet \mathcal{G}$ (which is easy under the Lipschitz assumption on $\ell$) and then bounding localized complexities of $\mathcal{F}$ using inequality (2.4).

It is not hard to combine the bounds of this type with model selection results of Section 5 to obtain various oracle inequalities for model selection in convex risk minimization problems. In particular, in the case of model selection for a sequence of function classes $\mathcal{G}_k := M_k \operatorname{conv}(\mathcal{H})$, where $\mathcal{H}$ is a VC-class, one would easily obtain a slight generalization of a recent result of Blanchard, Lugosi and Vayatis [10] on convergence rates of regularized boosting algorithm.

**8. Comments on regression and classification.** The general least squares regression is among statistical problems for which the penalization techniques have been very successful so far. In addition to already mentioned papers by Birgé and Massart [8], Barron, Birgé and Massart [3] and Massart [36], we refer the reader to a book by van de Geer [46], a book by Györfi, Kohler, Krzyzak and Walk [22] and papers by Baraud [2] and Kohler [25]. Our goal here is only to outline the connection of regression problems to a more general theory considered in the previous sections.

To simplify the matter, we consider only the case of least squares regression with bounded noise, that is, $T = [0, 1]$, $\ell(y, u) := (y - u)^2$. Thus, the regression problem is a convex risk minimization problem and it is well known and straightforward that in this case $g_*$ is the regression function: $g_*(x) := \mathbb{E}(Y | X = x)$. Given a class $\mathcal{G}$ of functions $g : S \mapsto [0, 1]$, a solution $\hat{g}_n$ of the empirical risk minimization problem (over the class $\mathcal{G}$) is a well-known least squares estimate of the regression function. The first problem of interest is to provide upper bounds on $\|\hat{g}_n - g_*\|_{L_2(\Pi)}$.

To relate this to the general framework of convex risk minimization, note that in this case $u_\mu := \operatorname{argmin}_u \int_0^1 (y - u)^2 \mu(dy) = \int_0^1 y\mu(dy)$ and by a very simple algebra

$$(\ell(y, u) - \ell(y, u_\mu))^2 = ((y - u)^2 - (y - u_\mu)^2)^2$$
$$= (u - u_\mu)^2 (2y - u - u_\mu)^2 \leq 4(u - u_\mu)^2$$

and

$$(8.1) \quad \int_0^1 (\ell(y, u) - \ell(y, u_\mu))\mu(dy) = \int_0^1 [(y - u)^2 - (y - u_\mu)^2]\mu(dy)$$
$$= (u - u_\mu)^2.$$

As a result, condition (7.1) holds with $D(u, \mu) \equiv 4$. Note also that identity (8.1) also implies (by integration) the formula $P(\ell \bullet g) - P(\ell \bullet g_*) = \|g - g_*\|^2$ that immediately reduces the study of $\|\hat{g}_n - g_*\|^2_{L_2(\Pi)}$ to excess risk bounds.

These observations allow one to simplify the arguments used in the previous section and to obtain the following result, using Theorem 1 and Lemma 5, more precisely; see the bound right after this lemma. In the case when the class $\mathcal{G}$ is convex, there is a way to improve the bound of the lemma. The key observation is that under the convexity assumption for all $g \in \mathcal{G}$, $\|g - \bar{g}\|^2_{L_2(\Pi)} \leq \|g - g_*\|^2_{L_2(\Pi)} - \|\bar{g} - g_*\|^2_{L_2(\Pi)}$ (see, e.g., [1], Lemma 20.9), which is a simplification and a specialization of the convexity inequalities used in the proof of Lemma 6.

THEOREM 13. *Let $\theta_n(\delta) := \theta_n(\mathcal{G}; \delta) := \theta_{n,\mathcal{G}}(\delta)$. There exists a constant $K$ such that for all $\varepsilon \in (0, 1]$*

$$\mathbb{P}\left\{\|\hat{g}_n - g_*\|^2_{L_2(\Pi)} \geq (1+\varepsilon)\inf_{h \in \mathcal{G}}\|h - g_*\|^2_{L_2(\Pi)} + K\left(\theta_n^\sharp\left(\frac{\varepsilon}{K}\right) + \frac{t+1}{\varepsilon n}\right)\right\}$$
$$\leq \log_q \frac{qn}{t} e^{-t}.$$

*If $\mathcal{G}$ is convex, then*

$$\mathbb{P}\left\{\|\hat{g}_n - g_*\|^2 \geq \inf_{g \in \mathcal{G}}\|g - g_*\|^2 + K\left(\theta_n^\sharp\left(\frac{1}{K}\right) + \frac{t+1}{n}\right)\right\} \leq \log_q \frac{qn}{t} e^{-t}.$$

*Moreover, if $\theta_n$ can be upper bounded by a function $\check{\theta}_n$ which is of strictly concave type, then one can replace $\theta_n$ by $\check{\theta}_n$ and drop the logarithmic factor in the bound.*

The significance of the above inequalities is related to the fact that in many particular cases of regression problem they allow one to recover asymptotically correct convergence rates. This follows from computations of local Rademacher complexities in particular examples, given in Section 2.5.

In the model selection framework, it is assumed that there exists a sequence $\mathcal{G}_k$ of classes of functions (models) available for least squares regression estimation. Let $\hat{g}_{n,k}$ denote a least squares estimate in the class $\mathcal{G}_k$. Given data-dependent complexity penalties $\hat{\pi}_n(k)$ associated with classes $\mathcal{G}_k$, we define the penalized least squares estimator as follows:

$$\hat{k} := \operatorname{argmin}\left[n^{-1}\sum_{j=1}^n (Y_j - \hat{g}_{n,k}(X_j))^2 + \hat{\pi}(k)\right], \qquad \hat{g}_n := \hat{g}_{n,\hat{k}}.$$

It is very natural to use penalization techniques of Theorems 6 and 11 to design complexity penalties and to establish oracle inequalities for the corresponding penalized least squares estimators.

EXAMPLE 1 (Dimension-based penalization). Suppose that for each $k$, $\mathcal{G}_k$ is a subset of a finite-dimensional subspace of $L_2(\Pi)$ of dimension $d_k$ and define $\hat{\pi}(k) := \hat{K}\frac{d_k + t_k + 1}{n}$ where $\hat{K}$ is some numerical constant (see Example 1 of Section 2.5). The following oracle inequality holds with some constant $C > 0$:

$$\mathbb{P}\left\{\|\hat{g}_n - g_*\|^2_{L_2(\Pi)} \geq C \inf_{k \geq 1}\left\{\inf_{g \in \mathcal{G}_k}\|g - g_*\|^2_{L_2(\Pi)} + \frac{d_k + t_k + 1}{n}\right\}\right\} \leq 4 \sum_{k=1}^{\infty} e^{-t_k}.$$

EXAMPLE 2 (Kernel selection with Mendelson's complexities). In this example, one is given a sequence $\{K_j\}$ of symmetric nonnegatively definite kernels on $S \times S$, $\mathcal{G}_j$ being the unit ball in the reproducing kernel Hilbert space $H_{K_j}$ (see Example 7 of Section 2.5). For each $j$, one can define empirical Mendelson's complexity and true Mendelson's complexity of the class $\mathcal{G}_j$, as in Section 2.5. We use the notations $\bar{\gamma}_{n,j}(\cdot) = \bar{\gamma}_n(\mathcal{G}_j;\cdot)$ and $\hat{\gamma}_{n,j}(\cdot) = \hat{\gamma}_n(\mathcal{G}_j;\cdot)$ and define $\hat{\pi}(j) := \hat{K}(\hat{\gamma}^\sharp_{n,j}(1) + \frac{t_j+1}{n})$, where $\hat{K}$ is a numerical constant. Then, the following oracle inequality holds:

$$\mathbb{P}\left\{\|\hat{g}_n - g_*\|^2_{L_2(\Pi)} \geq C \inf_{k \geq 1}\left\{\inf_{g \in \mathcal{G}_k}\|g - g_*\|^2_{L_2(\Pi)} + \left(\bar{\gamma}^\sharp_{n,k}(1) + \frac{t_k + 1}{n}\right)\right\}\right\}$$
$$\leq 4 \sum_{k=1}^{\infty} \log_q \frac{q^2 n}{t_k} e^{-t_k}.$$

EXAMPLE 3 (Penalization based on Rademacher complexities). One can also use localized Rademacher complexities, defined in Section 2.4 (see Lemma 1), as general penalties for model selection in regression problems. Namely, given a sequence of classes $\mathcal{G}_k$, we set

$$\hat{\pi}(k) := \hat{K}\left(\hat{\omega}^\sharp_{n,k}\left(\frac{1}{\hat{K}}\right) + \frac{t_k + 1}{n}\right) \quad \text{and} \quad \tilde{\pi}(k) := \tilde{K}\left(\bar{\omega}^\sharp_{n,k}\left(\frac{1}{\tilde{K}}\right) + \frac{t_k + 1}{n}\right)$$

with some (large enough) numerical constants $\hat{K}, \tilde{K}$. Here $\bar{\omega}_{n,k}(\cdot) = \bar{\omega}_n(\mathcal{G}_k;\cdot)$ and $\hat{\omega}_{n,k}(\cdot) = \hat{\omega}_n(\mathcal{G}_k;\cdot)$. Then we have (for a penalized least squares estimator $\hat{g}_n$) with some constant $C$

$$\mathbb{P}\left\{\|\hat{g}_n - g_*\|^2_{L_2(\Pi)} \geq C \inf_{k \geq 1}\left\{\inf_{g \in \mathcal{G}_k}\|g - g_*\|^2_{L_2(\Pi)} + \tilde{\pi}(k)\right\}\right\} \leq 4 \sum_{k=1}^{\infty} \log_q \frac{q^2 n}{t_k} e^{-t_k}.$$

We turn now to binary classification problems. In this case, $T := \{-1, 1\}$ and the loss function is chosen as $\ell(y, u) := I(y \neq u)$. The variable $Y$ is interpreted as an unobservable label associated with an observable instance $X$. Binary measurable functions $g: S \mapsto \{-1, 1\}$ are called *classifiers*. The goal of classification is to find a classifier that minimizes the *generalization error* (the probability of misclassification)

$$\mathbb{P}\{Y \neq g(X)\} = P\{(x, y): y \neq g(x)\} = P(\ell \bullet g),$$

so the classification problem becomes a version of a risk minimization problem with a binary loss function. Its solution always exists and is given by the following classifier (Bayes classifier): $g_*(x) := g_{*,P}(x) = I(\eta(x) \geq 0)$, where $\eta(x) := \mathbb{E}(Y|X=x)$ is the regression function (see [15]). However, the distribution $P$ of $(X,Y)$ and the regression function $\eta$ are unknown and the Bayes classifier is to be estimated based on the training data $(X_1, Y_1), \ldots, (X_n, Y_n)$ consisting of $n$ i.i.d. copies of $(X,Y)$. This is done by minimizing the so-called *training error*

$$n^{-1} \sum_{j=1}^{n} I(Y_j \neq g(X_j)) = P_n\{(x,y) : y \neq g(x)\} = P_n(\ell \bullet g)$$

over a suitable class of $\mathcal{G}$ of binary classifiers, which is equivalent to empirical risk minimization over the loss class $\mathcal{F} = \ell \bullet \mathcal{G}$, and all the theory developed in the previous sections applies to classification problems.

It is straightforward to check that condition (7.1) holds for binary loss $\ell$ with $D(u, \mu_x) = \frac{1}{|\eta(x)|}$ (moreover, the inequality in this case becomes an equality). If for some $C > 0$, $\alpha > 0$

$$\forall t > 0: \qquad \Pi\{x : 0 < |\eta(x)| \leq t\} \leq Ct^{\alpha},$$

then it easily follows that

(8.2) $$P(\ell \bullet g) - P(\ell \bullet g_*) \geq c_0 \rho_P^{2\kappa}(\ell \bullet g, \ell \bullet g_*),$$

where $\rho_P(\ell \bullet g_1, \ell \bullet g_2) := \Pi^{1/2}\{x : g_1(x) \neq g_2(x)\} = \Pi^{1/2}(g_1 - g_2)^2$, and $\kappa = \frac{1+\alpha}{\alpha}$ (see [44]). To get $\kappa = 1$, one can assume that for some $t_0 > 0$, $\Pi\{x : 0 < |\eta(x)| \leq t_0\} = 0$. Roughly, the assumptions of this type describe the degree of separation of two classes in classification problem, or the level of the "noise" in the labels ("low noise assumption"). Now one can use Theorem 8 of Section 6.1 to get the convergence rates in classification obtained first by Mammen and Tsybakov [35] and Tsybakov [44]. Namely, if $\mathcal{P}$ denotes a class of probability distributions on $S \times \{-1, 1\}$ and $\mathcal{G}$ is a class of binary classifiers such that, for all $P \in \mathcal{P}$, $g_{*,P} \in \mathcal{G}$, condition (8.2) holds (with the same $\kappa$ and $c_0$) and the $L_2(\Pi)$ bracketing entropy of the class $\mathcal{G}$ is of the order $O(\varepsilon^{-2\rho})$ as $\varepsilon \to 0$ uniformly in $P \in \mathcal{P}$ for some $\rho \in (0,1)$, then for a classifier $\hat{g}_n$ that minimizes the training error over $\mathcal{G}$ we have

$$\sup_{P \in \mathcal{P}} [P\{(x,y) : y \neq \hat{g}_n(x)\} - P\{(x,y) : y \neq g_{*,P}(x)\}] = O(n^{-\frac{\kappa}{2\kappa+\rho-1}}).$$

This was the result originally proved by Mammen and Tsybakov [35]. They also showed the convergence rate to be optimal in a minimax sense [35, 44]. As a consequence of Theorem 9, it is also easy to get an improvement of the model selection result of Tsybakov [44] (see Theorem 3 there) in the sense that our version of adaptation gives the precise convergence rates (Tsybakov's bounds involve an extra logarithmic factor).

Unfortunately, minimization of the training error over huge classes of binary functions (with entropy growing as $\varepsilon^{-\rho}$) is most often a computationally intractable problem. In so-called *large margin* classification algorithms (such as boosting and many algorithms for kernel machines) this difficulty is avoided by replacing the binary loss by a smooth (often, convex) loss function that dominates the binary loss, and using a version of functional gradient descent to minimize the corresponding empirical risk. In this setting, it is common to use real-valued functions $g$ as classifiers. At the end, $\mathrm{sign}(g(x))$ is computed to predict the label of an instance $x$. Let $\phi$ be a nonnegative convex function such that $\phi(u) \geq I(u \leq 0)$. We set $\ell(y, u) := \phi(yu)$ and look at a convex risk minimization problem $P(\ell \bullet g) \to \min$ and its empirical version $P_n(\ell \bullet g) \to \min$. Recently, Bartlett, Jordan and McAuliffe [6] and Blanchard, Lugosi and Vayatis [10] obtained reasonably good convergence rates for these types of algorithms. Their analysis is, essentially, a special version of somewhat more general analysis of convex risk minimization problems given in the previous sections.

## 9. Main Proofs.

PROOF OF PROPOSITION 1. For the first part, note that

$$\sum_{j\,:\,\delta_j \geq \delta} \frac{\psi(\delta_j)}{\delta_j} = \sum_{j\,:\,\delta_j \geq \delta} \frac{\psi(\delta_j)}{\delta_j^\gamma \delta_j^{1-\gamma}} \leq \frac{\psi(\delta)}{\delta^\gamma} \sum_{j\,:\,\delta_j \geq \delta} \frac{1}{\delta_j^{1-\gamma}}$$

$$= \frac{\psi(\delta)}{\delta} \sum_{j\,:\,\delta_j \geq \delta} \left(\frac{\delta}{\delta_j}\right)^{1-\gamma} \leq \frac{\psi(\delta)}{\delta} \sum_{j \geq 0} q^{-j(1-\gamma)} = c_{\gamma,q} \frac{\psi(\delta)}{\delta}.$$

To prove the second part, note that by induction $\bar{\delta}_k$ is nonincreasing and takes values in $[\bar{\delta}, 1]$. Denote $d_k := \bar{\delta}_k - \bar{\delta}$. We have

$$d_{k+1} = \bar{\delta}_{k+1} - \bar{\delta} \leq \psi(\bar{\delta}_k) - \psi(\bar{\delta}) = \frac{\psi(\bar{\delta}_k)}{\bar{\delta}_k^\gamma} \bar{\delta}_k^\gamma - \frac{\psi(\bar{\delta})}{\bar{\delta}^\gamma} \bar{\delta}^\gamma,$$

and since $\psi$ is of strictly concave type with exponent $\gamma$ and $\bar{\delta}_k \geq \bar{\delta}$, we get

$$d_{k+1} \leq \frac{\psi(\bar{\delta})}{\bar{\delta}^\gamma} (\bar{\delta}_k^\gamma - \bar{\delta}^\gamma) \leq \frac{\psi(\bar{\delta})}{\bar{\delta}} \bar{\delta}^{1-\gamma} (\bar{\delta}_k - \bar{\delta})^\gamma = \bar{\delta}^{1-\gamma} d_k^\gamma.$$

The result now follows by induction. $\square$

PROOF OF LEMMA 1. The first bound trivially follows from symmetrization inequality $\theta_n(\delta) \leq 2\bar{\omega}_n(\delta)$ and the definition of $\sharp$-transform. Let $\delta_j := q^{-j}$. In what follows $\delta = \delta_i$ for some $i$. To prove the second bound, define

$$E(\delta) := \left\{ \bar{\omega}_n(\delta) \leq \sup_{P(f-g)^2 \leq \delta} |R_n(f-g)| + \sqrt{2\frac{t}{n}(\delta + 2\bar{\omega}_n(\delta))} + \frac{8t}{3n} \right\}$$

$$\cap \left\{ \sup_{P(f-g)^2 \leq \delta} |(P_n - P)((f-g)^2)| \leq \mathbb{E} \sup_{P(f-g)^2 \leq \delta} |(P_n - P)((f-g)^2)| \right.$$
$$\left. + \sqrt{2\frac{t}{n}\left(\delta + 2\mathbb{E} \sup_{P(f-g)^2 \leq \delta} |(P_n - P)((f-g)^2)|\right)} + \frac{t}{3n} \right\}.$$

It follows from Talagrand's concentration inequalities that $\mathbb{P}(E(\delta)) \geq 1 - 2e^{-t}$. By symmetrization and contraction inequalities,
$$\mathbb{E} \sup_{P(f-g)^2 \leq \delta} |(P_n - P)((f-g)^2)| \leq 2\mathbb{E} \sup_{P(f-g)^2 \leq \delta} |R_n((f-g)^2)| \leq 8\bar{\omega}_n(\delta).$$

Therefore, on the event $E(\delta)$,
$$P(f-g)^2 \leq \delta \implies P_n(f-g)^2 \leq \delta + 8\bar{\omega}_n(\delta) + 2\sqrt{\frac{t}{2n}\delta} + 2\sqrt{\frac{t}{n}8\bar{\omega}_n(\delta)} + \frac{t}{3n},$$

and using the inequality $2ab \leq a^2 + b^2$ the right-hand side can be further bounded by $2\delta + 16\bar{\omega}_n(\delta) + \frac{2t}{n}$. Assuming that $\delta \geq q^{-1}\bar{\omega}_n^{\sharp,q}(\varepsilon) \geq \frac{t}{n}$, and using the monotonicity of $\bar{\omega}_n$, we get

$$\bar{\omega}_n(\delta) \leq \delta \sup_{\delta_j \geq q^{-1}\bar{\omega}_n^{\sharp,q}(\varepsilon)} \frac{\bar{\omega}_n(\delta_j)}{\delta_j} \leq \delta \sup_{\delta_j \geq q^{-1}\bar{\omega}_n^{\sharp,q}(\varepsilon)} \frac{\bar{\omega}_n(q\delta_j)}{\delta_j}$$
$$\leq q\delta \sup_{\delta_j \geq \bar{\omega}_n^{\sharp,q}(\varepsilon)} \frac{\bar{\omega}_n(\delta_j)}{\delta_j} \leq q\varepsilon\delta.$$

Therefore, for $\varepsilon \in (0, 1]$ and $\delta \geq q^{-1}\bar{\omega}_n^{\sharp,q}(\varepsilon) \geq t/n$, on the event $E(\delta)$
$$P(f-g)^2 \leq \delta \implies P_n(f-g)^2 \leq 2\delta + 16\bar{\omega}_n(\delta) + \frac{2t}{n} \leq (4 + 16q)\delta.$$

Also, on the same event and under the same conditions,
$$\bar{\omega}_n(\delta) \leq \sup_{P(f-g)^2 \leq \delta} |R_n(f-g)| + \sqrt{2\frac{t}{n}(\delta + 2\bar{\omega}_n(\delta))} + \frac{8t}{3n}$$
$$\leq \sup_{P_n(f-g)^2 \leq (4+16q)\delta} |R_n(f-g)| + \sqrt{2\frac{t}{n}\delta} + 2\sqrt{\frac{\bar{\omega}_n(\delta)}{2}\frac{2t}{n}} + \frac{8t}{3n}$$
$$\leq \sup_{P_n(f-g)^2 \leq (4+16q)\delta} |R_n(f-g)| + \sqrt{2\frac{t}{n}\delta} + \frac{8t}{3n} + \frac{2t}{n} + \frac{\bar{\omega}_n(\delta)}{2},$$

where we again used the inequality $2ab \leq a^2 + b^2$. Therefore, on the event $E(\delta)$
$$\bar{\omega}_n(\delta) \leq 2 \sup_{P_n(f-g)^2 \leq (4+16q)\delta} |R_n(f-g)| + 2\sqrt{2}\sqrt{\frac{t}{n}\delta} + \frac{10t}{n}$$
$$= 2\hat{\omega}_n((4+16q)\delta) + 2\sqrt{2}\sqrt{\frac{t}{n}\delta} + \frac{10t}{n} =: \psi(\delta)$$

as soon as $\delta \geq q^{-1}\bar{\omega}_n^{\sharp,q}(\varepsilon) \geq \frac{t}{n}$.

Note that if $q^{-1}\bar{\omega}_n^{\sharp,q}(\varepsilon) < \frac{t}{n}$, then the second bound of the lemma is trivially satisfied. Otherwise, denote

$$E := \bigcap_{j\,:\,\delta_j \geq q^{-1}\bar{\omega}_n^{\sharp,q}(\varepsilon) \geq \frac{t}{n}} E(\delta_j).$$

Clearly, $\mathbb{P}(E) \geq 1 - 2\log_q \frac{qn}{t} e^{-t}$, and, on the event $E$, we have $\bar{\omega}_n(\delta_j) \leq \psi(\delta_j)$ for all $\delta_j \geq q^{-1}\bar{\omega}_n^{\sharp,q}(\varepsilon)$, which implies that (see Property 2' in Section 2.3) $\bar{\omega}_n^{\sharp,q}(\varepsilon) \leq \psi^{\sharp,q}(\varepsilon)$. Using the properties of $\sharp$-transform, this yields by a simple computation that

$$\bar{\omega}_n^\sharp(\varepsilon) \leq C\left(\hat{\omega}_n^\sharp(c\varepsilon) + \frac{t}{n\varepsilon^2}\right)$$

with some constants $C, c$ depending only on $q$.

To prove the third bound, we introduce the following event: $F := \bigcap_{\delta_j \geq \frac{t}{n}} F(\delta_j)$, where

$$F(\delta) := \left\{\sup_{P(f-g)^2 \leq c_q\delta} |R_n(f-g)| \leq \bar{\omega}_n(c_q\delta) + \sqrt{2\frac{t}{n}(c_q\delta + 2\bar{\omega}_n(c_q\delta))} + \frac{t}{3n}\right\}$$

$$\cap \left\{\sup_{P(f-g)^2 \leq \delta} |(P_n-P)((f-g)^2)| \leq \mathbb{E}\sup_{P(f-g)^2 \leq \delta} |(P_n-P)((f-g)^2)|\right.$$

$$\left. + \sqrt{2\frac{t}{n}\left(\delta + 2\mathbb{E}\sup_{P(f-g)^2 \leq \delta} |(P_n-P)((f-g)^2)|\right)} + \frac{t}{3n}\right\}$$

with a constant $c_q$ depending only on $q$ to be chosen later on. It follows from Talagrand's concentration inequalities that $\mathbb{P}(F) \geq 1 - 2\log_q \frac{qn}{t} e^{-t}$. Let $\delta = \delta_i$ for some $i$ and $\delta_i \geq \frac{t}{n}$. On the event $F$ the following implication holds:

$$P_n(f-g)^2 \leq \delta \quad \text{and} \quad P(f-g)^2 \in (\delta_{j+1}, \delta_j]$$

$$\implies \frac{\delta_j}{q} = \delta_{j+1} \leq P(f-g)^2 \leq \delta + \sup_{P(f-g)^2 \leq \delta_j} |(P_n-P)((f-g)^2)|$$

$$\leq \delta + 16\bar{\omega}_n(\delta_j) + \frac{\delta_j}{q^2} + \frac{(4/3 + q^2/2)t}{n},$$

where we used the same computation as in the previous part of the proof with minor modifications. If $\delta_j \geq \bar{\omega}_n^{\sharp,q}(\varepsilon)$, then $\bar{\omega}_n(\delta_j) \leq \varepsilon\delta_j$, and we can get

$$\delta_j(q^{-1} - q^{-2} - 16\varepsilon) \leq \delta + \frac{(4/3 + q^2/2)t}{n}.$$

If $\varepsilon < \frac{1}{32}(q^{-1} - q^{-2})$ (note that it is enough to prove the bound under this restriction and the general case would follow by changing the constants), then we get that

$$\delta_j \leq 2(q^{-1} - q^{-2})^{-1}\left(\delta + \frac{(4/3 + q^2/2)t}{n}\right).$$

What we proved so far can be formulated as follows: on the event $F$, for $\delta = \delta_i \geq \frac{t}{n}$,

$$P_n(f-g)^2 \leq \delta$$
$$\implies P(f-g)^2 \leq 2(q^{-1} - q^{-2})^{-1}\left(\delta + \frac{(4/3 + q^2/2)t}{n}\right) \vee \bar{\omega}_n^{\sharp,q}(\varepsilon),$$

which means that for $\delta \geq \bar{\omega}_n^{\sharp,q}(\varepsilon)$, $P_n(f-g)^2 \leq \delta \Rightarrow P(f-g)^2 \leq c_q \delta$ with a constant $c_q > 1$ depending only on $q$. This allows us to conclude that on the event $F$ for all $\delta = \delta_i \geq \bar{\omega}_n^{\sharp,q}(\varepsilon) \vee \frac{t}{n}$

$$\hat{\omega}_n(\delta) \leq \sup_{P(f-g)^2 \leq c_q \delta} |R_n(f-g)| \leq \bar{\omega}_n(c_q \delta) + \sqrt{2\frac{t}{n}(c_q \delta + 2\bar{\omega}_n(c_q \delta))} + \frac{t}{3n}$$

$$\leq 2\bar{\omega}_n(c_q \delta) + \sqrt{2c_q \delta \frac{t}{n}} + \frac{2t}{n} =: \psi(\delta).$$

Next we use the basic properties of the $\sharp$-transform to conclude the proof. Since $\psi(\delta) \geq \bar{\omega}_n(\delta) \vee \frac{t}{n}$, we get for all $\varepsilon \in (0,1]$, $\psi^{\sharp,q}(\varepsilon) \geq \bar{\omega}_n^{\sharp,q}(\varepsilon) \vee \frac{t}{n}$. Thus, for all $\delta \geq \psi^{\sharp,q}(\varepsilon)$, $\hat{\omega}_n(\delta) \leq \psi(\delta)$, implying that $\hat{\omega}_n^{\sharp,q}(\varepsilon) \leq \psi^{\sharp,q}(\varepsilon)$. Now it is easy to conclude that on the event $F$

$$\hat{\omega}_n^{\sharp}(\varepsilon) \leq C\left(\bar{\omega}_n^{\sharp}(c\varepsilon) + \frac{t}{n\varepsilon^2}\right)$$

with some constants $C, c$ depending only on $q$.

The proof for $\omega_{n,r}^{\sharp}$ is similar. $\square$

PROOF OF THEOREM 1. Let

$$E_{n,j}(t) := \left\{\sup_{f,g \in \mathcal{F}(\delta_j)} |(P_n - P)(f-g)| \leq U_n(\delta_j; t)\right\}.$$

By Talagrand's concentration inequality, $\mathbb{P}((E_{n,j}(t))^c) \leq e^{-t}$. Let $\delta_j \geq \delta$. Since on the event $E_{n,j}(t)$,

$$\hat{f}_n \in \mathcal{F}(\delta_{j+1}, \delta_j]$$
$$\implies \forall \varepsilon \in (0, \delta_{j+1}) \; \forall g \in \mathcal{F}(\varepsilon)$$
$$\delta_{j+1} < \mathcal{E}(\hat{f}_n) \leq P\hat{f}_n - Pg + \varepsilon$$

$$\leq P_n \hat{f}_n - P_n g + (P - P_n)(\hat{f}_n - g) + \varepsilon$$
$$\leq \hat{\mathcal{E}}_n(\hat{f}_n) + \sup_{f,g \in \mathcal{F}(\delta_j)} |(P_n - P)(f - g)| + \varepsilon$$
$$\leq U_n(\delta_j; t) + \varepsilon \leq V_n(\delta; t)\delta_j + \varepsilon$$
$$\implies V_n(\delta; t) \geq \frac{1}{q} > \frac{1}{2q}$$
$$\implies \delta \leq U_{n,t}^{\sharp,q}\left(\frac{1}{2q}\right) = \delta_n(t),$$

we can conclude that, for $\delta_j \geq \delta \geq \delta_n(t)$, $\{\hat{f}_n \in \mathcal{F}(\delta_{j+1}, \delta_j]\} \subset (E_{n,j}(t))^c$. Therefore, for $\delta \geq \delta_n(t)$, on the event $E_n(t) := \bigcap_{j : \delta_j \geq \delta} E_{n,j}(t)$ we have $\mathcal{E}(\hat{f}_n) \leq \delta$, implying that

$$\mathbb{P}\{\mathcal{E}(\hat{f}_n) > \delta\} \leq \sum_{j : \delta_j \geq \delta} \mathbb{P}((E_{n,j}(t))^c) \leq \log_q \frac{q}{\delta} e^{-t}.$$

Now, on the event $E_n(t)$, we have $\hat{f}_n \in \mathcal{F}(\delta)$ and for all $j$ such that $\delta_j \geq \delta$

$$f \in \mathcal{F}(\delta_{j+1}, \delta_j]$$
$$\implies \forall \varepsilon \in (0, \delta_j) \ \forall g \in \mathcal{F}(\varepsilon)$$
$$\mathcal{E}(f) \leq Pf - Pg + \varepsilon \leq P_n f - P_n g + (P - P_n)(f - g) + \varepsilon$$
$$\leq \hat{\mathcal{E}}_n(f) + U_n(\delta_j; t) + \varepsilon \leq \hat{\mathcal{E}}_n(f) + V_n(\delta; t)\delta_j + \varepsilon$$
$$\leq \hat{\mathcal{E}}_n(f) + qV_n(\delta; t)\mathcal{E}(f) + \varepsilon,$$

which means that on this event $\mathcal{E}(f) \geq \delta \Rightarrow \hat{\mathcal{E}}_n(f) \geq (1 - qV_n(\delta; t))\mathcal{E}(f)$. Similarly, we have on $E_n(t)$

$$f \in \mathcal{F}(\delta_{j+1}, \delta_j]$$
$$\implies \hat{\mathcal{E}}_n(f) = P_n f - P_n \hat{f}_n \leq Pf - P\hat{f}_n + (P_n - P)(f - \hat{f}_n)$$
$$\leq \mathcal{E}(f) + U_n(\delta_j; t) \leq \mathcal{E}(f) + V_n(\delta; t)\delta_j$$
$$\leq \mathcal{E}(f) + qV_n(\delta; t)\mathcal{E}(f) = (1 + qV_n(\delta; t))\mathcal{E}(f),$$

so that $\mathcal{E}(f) > \delta \Rightarrow \hat{\mathcal{E}}_n(f) \leq (1 + qV_n(\delta; t))\mathcal{E}(f)$. Since $\mathbb{P}((E_n(t))^c) \leq \log_q \frac{q}{\delta} e^{-t}$, the result follows. $\square$

PROOF OF LEMMA 2. Consider the following event:

$$E := \left\{ \forall f \in \mathcal{F} \text{ with } \mathcal{E}(f) \geq \delta_n(t) : \frac{1}{2} \leq \frac{\hat{\mathcal{E}}_n(f)}{\mathcal{E}(f)} \leq \frac{3}{2} \right\}.$$

It follows from Theorem 1 and the definition of $\delta_n(t)$ that $\mathbb{P}(E) \geq 1 - \log_q \frac{q}{\delta_n(t)} e^{-t}$. Consider also

$$F := \Big\{ \sup_{f,g \in \mathcal{F}(\delta_n(t))} |(P_n - P)(f - g)| \leq U_n(\delta_n(t); t) \Big\}.$$

It follows from the concentration inequality that $\mathbb{P}(F) \geq 1 - e^{-t}$. Therefore,

$$\mathbb{P}(E \cap F) \geq 1 - \log_q \frac{q^2}{\delta_n(t)} e^{-t}.$$

On the event $E$, we have

(9.1) $$\forall f \in \mathcal{F}: \quad \mathcal{E}(f) \leq 2\hat{\mathcal{E}}_n(f) \vee \delta_n(t),$$

which implies that for all $\delta \geq \delta_n(t)$, $\hat{\mathcal{F}}_n(\delta) \subset \mathcal{F}(2\delta)$. On the other hand, on the same event $E$, $\forall f \in \mathcal{F} : \mathcal{E}(f) \geq \delta_n(t) \Rightarrow \hat{\mathcal{E}}_n(f) \leq \frac{3}{2}\mathcal{E}(f)$.

On the event $F$,

$$\mathcal{E}(f) \leq \delta_n(t) \implies \hat{\mathcal{E}}_n(f) \leq \mathcal{E}(f) + \sup_{f,g \in \mathcal{F}(\delta_n(t))} |(P_n - P)(f - g)|$$

$$\leq \mathcal{E}(f) + U_n(\delta_n(t); t)$$

$$\leq \delta_n(t) + qV_n(\delta_n(t); t)\delta_n(t) \leq \frac{3}{2}\delta_n(t).$$

Thus, on the event $E \cap F$

(9.2) $$\forall f \in \mathcal{F}: \quad \hat{\mathcal{E}}_n(f) \leq \frac{3}{2}(\mathcal{E}(f) \vee \delta_n(t)),$$

which implies that $\forall \delta \geq \delta_n(t) : \mathcal{F}(\delta) \subset \hat{\mathcal{F}}_n(3\delta/2)$. □

PROOF OF THEOREM 2. It is similar to the proof of Theorem 1, but now our goal is to avoid using the concentration inequality many times (for each $\delta_j$) since this leads to a logarithmic factor. The trick was previously used in [36] and in the Ph.D. dissertation of Bousquet (see also [5]). Define

$$\mathcal{G}_\delta := \bigcup_{j : \delta_j \geq \delta} \frac{\delta}{\delta_j} \{f - g : f, g \in \mathcal{F}(\delta_j)\}.$$

Then the functions in $\mathcal{G}_\delta$ are bounded by 1 and

$$\sigma_P(\mathcal{G}_\delta) \leq \sup_{j : \delta_j \geq \delta} \frac{\delta}{\delta_j} \sup_{f,g \in \mathcal{F}(\delta_j)} \sigma_P(f - g) \leq \delta \sup_{j : \delta_j \geq \delta} \frac{\check{D}(\delta_j)}{\delta_j} \leq \check{D}(\delta),$$

since $\check{D}$ is of concave type. Also, since $\check{\phi}_n$ is of strictly concave type, Proposition 1 gives

$$\mathbb{E}\|P_n - P\|_{\mathcal{G}_\delta} = \mathbb{E}\sup_{j\,:\,\delta_j \geq \delta} \frac{\delta}{\delta_j} \sup_{f,g \in \mathcal{F}(\delta_j)} |(P_n - P)(f-g)|$$

$$\leq \sum_{j\,:\,\delta_j \geq \delta} \frac{\delta}{\delta_j} \mathbb{E} \sup_{f,g \in \mathcal{F}(\delta_j)} |(P_n - P)(f-g)|$$

$$\leq \delta \sum_{j\,:\,\delta_j \geq \delta} \frac{\check{\phi}_n(\delta_j)}{\delta_j} \leq c_{\gamma,q} \check{\phi}_n(\delta).$$

Now Talagrand's concentration inequality implies that there exists an event $E$ of probability $\mathbb{P}(E) \geq 1 - e^{-t}$ such that on this event $\|P_n - P\|_{\mathcal{G}_\delta} \leq \check{U}_n(\delta;t)$ (the constant $\check{K}$ in the definition of $\check{U}_n(\delta;t)$ should be chosen properly). Then on the event $E$

$$\forall j \text{ with } \delta_j \geq \delta\colon \qquad \sup_{f,g \in \mathcal{F}(\delta_j)} |(P_n - P)(f-g)| \leq \frac{\delta_j}{\delta} \check{U}_n(\delta;t) \leq \check{V}_n(\delta;t)\delta_j.$$

The rest repeats the proof of Theorem 1. $\square$

REMARK. There is also a way to prove a bound on $\mathcal{E}_P(\hat{f})$ based on the iterative localization method described in the Introduction and in the second statement of Proposition 1. Namely, one can assume that both $\check{\phi}_n$ and $\check{D}$ are of strictly concave type with exponent $\gamma \in (0,1)$. As a result, the function $\check{U}_{n,t}$ is also of strictly concave type with the same exponent. If now $\check{\delta}_n(t)$ denotes its fixed point, then by Proposition 1(ii), the number $N$ of iterations needed to achieve the bound $\bar{\delta}_N \leq 2\check{\delta}_n(t)$ is smaller than $\log\log_2((1-\check{\delta}_n(t))/\check{\delta}_n(t))/\log(1/\gamma) + 1$ in the case when $\check{\delta}_n(t) < 1/2$ and $N = 1$ otherwise. Thus, the argument described in the Introduction immediately shows that $\mathbb{P}\{\mathcal{E}_P(\hat{f}) \geq \check{\delta}_n(t)\} \leq Ne^{-t}$. This approach was first used in [27] (and later also in some of the arguments of [5]).

PROOF OF THEOREM 3. The proof consists of several steps. Throughout, $H$ will denote the event introduced in Lemma 2. According to this lemma, we have $\mathbb{P}(H) \geq 1 - \log_q \frac{q^2}{\delta_n(t)} e^{-t}$.

Step 1. *Bounding the Rademacher complexity.* Using Talagrand's concentration inequality, we get (for $\delta > 0$ and $t > 0$) on an event $F = F(\delta)$ with probability at least $1 - e^{-t}$

$$\mathbb{E}\sup_{f,g \in \mathcal{F}(\delta)} |R_n(f-g)| \leq \sup_{f,g \in \mathcal{F}(\delta)} |R_n(f-g)|$$
$$+ \sqrt{\frac{2t}{n}\Big(D^2(\delta) + 2\mathbb{E}\sup_{f,g \in \mathcal{F}(\delta)} |R_n(f-g)|\Big)} + \frac{8t}{3n},$$

which implies that

$$\mathbb{E} \sup_{f,g \in \mathcal{F}(\delta)} |R_n(f-g)| \leq \sup_{f,g \in \mathcal{F}(\delta)} |R_n(f-g)| + D(\delta)\sqrt{\frac{2t}{n}} + \frac{8t}{3n}$$

$$+ 2\sqrt{\frac{1}{2}\mathbb{E} \sup_{f,g \in \mathcal{F}(\delta)} |R_n(f-g)|\frac{2t}{n}}$$

$$\leq \sup_{f,g \in \mathcal{F}(\delta)} |R_n(f-g)| + D(\delta)\sqrt{\frac{2t}{n}} + \frac{8t}{3n}$$

$$+ \frac{1}{2}\mathbb{E} \sup_{f,g \in \mathcal{F}(\delta)} |R_n(f-g)| + \frac{2t}{n},$$

or

$$\mathbb{E} \sup_{f,g \in \mathcal{F}(\delta)} |R_n(f-g)| \leq 2 \sup_{f,g \in \mathcal{F}(\delta)} |R_n(f-g)| + 2\sqrt{2}D(\delta)\sqrt{\frac{t}{n}} + \frac{28t}{3n}.$$

This can be further bounded using Lemma 2. Namely, for all $\delta \geq \delta_n(t)$, we have on the event $H \cap F$ that

$$\mathbb{E} \sup_{f,g \in \mathcal{F}(\delta)} |R_n(f-g)| \leq 2 \sup_{f,g \in \hat{\mathcal{F}}_n(\frac{3}{2}\delta)} |R_n(f-g)| + 2\sqrt{2}D(\delta)\sqrt{\frac{t}{n}} + \frac{28t}{3n}.$$

*Step* 2. *Bounding the diameter* $D(\delta)$. Again, we apply Talagrand's concentration inequality to get on an event $G = G(\delta)$ with probability at least $1 - e^{-t}$

$$D^2(\delta) = \sup_{f,g \in \mathcal{F}(\delta)} P(f-g)^2$$

$$\leq \sup_{f,g \in \mathcal{F}(\delta)} P_n(f-g)^2 + \sup_{f,g \in \mathcal{F}(\delta)} |(P_n - P)((f-g)^2)|$$

$$\leq \sup_{f,g \in \mathcal{F}(\delta)} P_n(f-g)^2 + \mathbb{E} \sup_{f,g \in \mathcal{F}(\delta)} |(P_n - P)((f-g)^2)|$$

$$+ \sqrt{\frac{2t}{n}\left(D^2(\delta) + 2\mathbb{E} \sup_{f,g \in \mathcal{F}(\delta)} |(P_n - P)((f-g)^2)|\right)} + \frac{t}{3n},$$

where we also used that $\sup_{f,g \in \mathcal{F}(\delta)} \text{Var}_P((f-g)^2) \leq \sup_{f,g \in \mathcal{F}(\delta)} P(f-g)^2 = D^2(\delta)$, since the functions from $\mathcal{F}$ take their values in $[0,1]$. Using the symmetrization inequality and then the contraction inequality for Rademacher processes, we get

$$\mathbb{E} \sup_{f,g \in \mathcal{F}(\delta)} |(P_n - P)(f-g)^2| \leq 2\mathbb{E} \sup_{f,g \in \mathcal{F}(\delta)} |R_n((f-g)^2)|$$

$$\leq 8\mathbb{E} \sup_{f,g \in \mathcal{F}(\delta)} |R_n(f-g)|.$$

It follows from Lemma 2 that for all $\delta \geq \bar{\delta}_n(t)$ on the event $H$ we have

$$\sup_{f,g\in\mathcal{F}(\delta)} P_n(f-g)^2 \leq \sup_{f,g\in\hat{\mathcal{F}}_n(3\delta/2)} P_n(f-g)^2 = \hat{D}_n^2\left(\frac{3}{2}\delta\right).$$

Hence, on the event $H \cap G$

$$D^2(\delta) \leq \hat{D}_n^2\left(\frac{3}{2}\delta\right) + 8\mathbb{E}\sup_{f,g\in\mathcal{F}(\delta)} |R_n(f-g)| + D(\delta)\sqrt{\frac{2t}{n}}$$

$$+ 2\sqrt{\frac{8t}{n}\mathbb{E}\sup_{f,g\in\mathcal{F}(\delta)} |R_n(f-g)|} + \frac{t}{3n}$$

$$\leq \hat{D}_n^2\left(\frac{3}{2}\delta\right) + 9\mathbb{E}\sup_{f,g\in\mathcal{F}(\delta)} |R_n(f-g)| + D(\delta)\sqrt{\frac{2t}{n}} + \frac{9t}{n},$$

where we applied the inequality $2\sqrt{ab} \leq a+b$, $a,b \geq 0$. Next we use the resulting bound of Step 1 to get on $H \cap F \cap G$

$$D^2(\delta) \leq \hat{D}_n^2\left(\frac{3}{2}\delta\right) + 18\sup_{f,g\in\hat{\mathcal{F}}_n(3\delta/2)} |R_n(f-g)| + 19D(\delta)\sqrt{\frac{2t}{n}} + \frac{100t}{n}.$$

As before, we bound the term $19D(\delta)\sqrt{\frac{2t}{n}} = 2 \times 19\frac{D(\delta)}{\sqrt{2}}\sqrt{\frac{t}{n}}$ using the inequality $2ab \leq a^2 + b^2$ and this yields

$$D^2(\delta) \leq \frac{1}{2}D^2(\delta) + \hat{D}_n^2\left(\frac{3}{2}\delta\right) + 18\sup_{f,g\in\hat{\mathcal{F}}_n(3\delta/2)} |R_n(f-g)| + \frac{500t}{n}.$$

As a result, we get the following bound holding on the event $H \cap F \cap G$:

$$D^2(\delta) \leq 2\hat{D}_n^2\left(\frac{3}{2}\delta\right) + 36\sup_{f,g\in\hat{\mathcal{F}}_n(3\delta/2)} |R_n(f-g)| + \frac{1000t}{n},$$

which also implies

$$D(\delta) \leq \sqrt{2}\hat{D}_n\left(\frac{3}{2}\delta\right) + 6\sqrt{\sup_{f,g\in\hat{\mathcal{F}}_n(3\delta/2)} |R_n(f-g)|} + \frac{32t}{n}.$$

*Step 3. Bounding $\bar{U}_n$ in terms of $\hat{U}_n$.* We use the bound on $D(\delta)$ in terms of $\hat{D}_n(\frac{3}{2}\delta)$ (Step 2) to derive from the bound of Step 1 that

$$\mathbb{E}\sup_{f,g\in\mathcal{F}(\delta)} |R_n(f-g)| \leq 2\sup_{f,g\in\hat{\mathcal{F}}_n(3\delta/2)} |R_n(f-g)| + 4\hat{D}_n\left(\frac{3}{2}\delta\right)\sqrt{\frac{t}{n}}$$

$$+ 12\sqrt{2}\sqrt{\sup_{f,g\in\hat{\mathcal{F}}_n(3\delta/2)} |R_n(f-g)|}\sqrt{\frac{t}{n}} + \frac{100t}{n}$$

$$\leq 3 \sup_{f,g \in \hat{\mathcal{F}}_n(3\delta/2)} |R_n(f-g)| + 4\hat{D}_n\left(\frac{3}{2}\delta\right)\sqrt{\frac{t}{n}} + \frac{172t}{n},$$

which holds on the event $H \cap F \cap G$. By the symmetrization inequality, we also have

$$\mathbb{E} \sup_{f,g \in \mathcal{F}(\delta)} |(P_n - P)(f-g)| \leq 6 \sup_{f,g \in \hat{\mathcal{F}}_n(3\delta/2)} |R_n(f-g)| + 8\hat{D}_n\left(\frac{3}{2}\delta\right)\sqrt{\frac{t}{n}} + \frac{344t}{n},$$

which holds on the same event. Recalling the definition of $\bar{U}_n$ and $\hat{U}_n$, the last bound together with the bound of Step 2 shows that with a straightforward choice of numerical constants $\hat{K}, \hat{c}$ the following bound is true on the event $H \cap F \cap G$: $\bar{U}_n(\delta; t) \leq \hat{U}_n(\delta; t)$.

*Step* 4. *Bounding* $\hat{U}_n$ *in terms of* $\tilde{U}_n$. The derivation is similar to the previous one. First, by Lemma 2 and Talagrand's concentration inequality, for all $\delta \geq \delta_n(t)$,

$$\sup_{f,g \in \hat{\mathcal{F}}_n(\delta)} |R_n(f-g)| \leq \sup_{f,g \in \mathcal{F}(2\delta)} |R_n(f-g)| \leq \mathbb{E} \sup_{f,g \in \mathcal{F}(2\delta)} |R_n(f-g)|$$

$$+ \sqrt{2\frac{t}{n}\left(D^2(2\delta) + \mathbb{E}\sup_{f,g \in \mathcal{F}(2\delta)} |R_n(f-g)|\right)} + \frac{8t}{3n}$$

on the event $H \cap F'$, where $F' = F'(\delta)$ is such that $\mathbb{P}(F') \geq 1 - e^{-t}$. Next, using the desymmetrization inequality,

$$\mathbb{E} \sup_{f,g \in \mathcal{F}(2\delta)} |R_n(f-g)|$$

$$\leq \mathbb{E} \sup_{f,g \in \mathcal{F}(2\delta)} |R_n(f - g - P(f-g))| + \sup_{f,g \in \mathcal{F}(2\delta)} |P(f-g)|\mathbb{E}|R_n(1)|$$

$$\leq 2\mathbb{E} \sup_{f,g \in \mathcal{F}(2\delta)} |(P_n - P)(f-g)| + n^{-1/2} \sup_{f,g \in \mathcal{F}(2\delta)} P^{1/2}(f-g)^2$$

$$\leq 2\phi_n(2\delta) + n^{-1/2}D(2\delta).$$

Therefore, we get (by getting rid of $\phi_n$ under the square root)

$$\sup_{f,g \in \hat{\mathcal{F}}_n(\delta)} |R_n(f-g)| \leq 4\phi_n(2\delta) + D(2\delta)\left(\frac{1}{\sqrt{n}} + \sqrt{2}\sqrt{\frac{t}{n}}\right) + \frac{4t}{n}.$$

We turn now to bounding the empirical diameter $\hat{D}_n(\delta)$. Again, by Lemma 2 and Talagrand's concentration inequality, we have for all $\delta \geq \bar{\delta}_n(t)$ on the event $H \cap G'$, where $G' = G'(\delta)$ is such that $\mathbb{P}(G') \geq 1 - e^{-t}$,

$$\hat{D}_n^2(\delta) := \sup_{f,g \in \hat{\mathcal{F}}_n(\delta)} P_n(f-g)^2 \leq \sup_{f,g \in \mathcal{F}(2\delta)} P_n(f-g)^2$$

$$\leq \sup_{f,g \in \mathcal{F}(2\delta)} P(f-g)^2 + \sup_{f,g \in \mathcal{F}(2\delta)} |(P_n - P)((f-g)^2)|$$

$$\leq D^2(2\delta) + \mathbb{E} \sup_{f,g \in \mathcal{F}(2\delta)} |(P_n - P)((f-g)^2)|$$

$$+ \sqrt{2\frac{t}{n}\Big(D^2(2\delta) + 2\mathbb{E} \sup_{f,g \in \mathcal{F}(2\delta)} |(P_n - P)((f-g)^2)|\Big)} + \frac{t}{3n}.$$

As in Step 2, we use symmetrization and contraction inequalities to get

$$\mathbb{E} \sup_{f,g \in \mathcal{F}(2\delta)} |(P_n - P)((f-g)^2)| \leq 8\mathbb{E} \sup_{f,g \in \mathcal{F}(2\delta)} |R_n(f-g)|,$$

and then using the desymmetrization bound, as in Step 3, to get

$$\mathbb{E} \sup_{f,g \in \mathcal{F}(2\delta)} |(P_n - P)((f-g)^2)| \leq 16\phi_n(2\delta) + 8\frac{D(2\delta)}{\sqrt{n}}.$$

By a simple computation this implies that

$$\hat{D}_n^2(\delta) \leq D^2(2\delta) + 32\phi_n(2\delta) + D(2\delta)\left(\sqrt{\frac{2t}{n}} + \frac{16}{\sqrt{n}}\right) + \frac{2t}{n}.$$

The same algebra we already used in Step 3 yields the inequality $\hat{U}_n(\delta;t) \leq \tilde{U}_n(\delta;t)$ that holds on the event $H \cap F' \cap G'$ with properly chosen numerical constants $\tilde{K}, \tilde{c}$ in the definition of $\tilde{U}_n$.

*Step* 5. *Conclusion.* Using the inequalities of Steps 4 and 5 for $\delta = \delta_j \geq \delta_n(t)$ gives

$$\mathbb{P}(E) \geq 1 - \left(\log_q \frac{q^2}{\delta_n(t)} + 4\log_q \frac{q}{\delta_n(t)}\right)\exp\{-t\},$$

where

$$E := \{\forall \delta_j \geq \bar{\delta}_n(t) : \bar{U}_n(\delta_j;t) \leq \hat{U}_n(\delta_j;t) \leq \tilde{U}_n(\delta_j;t)\},$$

since

$$E \supset \bigcup_{j \,:\, \delta_j \geq \bar{\delta}_n(t)} (H \cap F(\delta_j) \cap G(\delta_j) \cap F'(\delta_j) \cap G'(\delta_j)).$$

Applying to $\psi(\delta) := \bar{U}_{n,t}(\delta)$ property $7'$ of the $\sharp, q$-transform, we get with $c = q^2$

$$q^2\delta_n(t) = q^2 U_{n,t}^{\sharp,q}\left(\frac{1}{2q}\right) \leq q^2 \bar{U}_{n,t}^{\sharp,q}\left(\frac{1}{2q}\right) \leq \bar{U}_{n,t}^{\sharp,q}\left(\frac{1}{2q^3}\right) = \bar{\delta}_n(t).$$

Therefore, using property $2'$ of the $\sharp, q$-transform, we get on the event $E$

$$\bar{\delta}_n(t) = \bar{U}_{n,t}^{\sharp,q}\left(\frac{1}{2q^3}\right) \leq \hat{\delta}_n(t) = \hat{U}_{n,t}^{\sharp,q}\left(\frac{1}{2q^3}\right)$$

and then, repeating the same argument for $\hat{\delta}_n(t)$, that

$$\hat{\delta}_n(t) = \hat{U}_{n,t}^{\sharp,q}\left(\frac{1}{2q^3}\right) \leq \tilde{\delta}_n(t) = \tilde{U}_{n,t}^{\sharp,q}\left(\frac{1}{2q^3}\right),$$

implying the result. $\square$

PROOF OF THEOREM 4. Denote

$$\breve{\psi}_n^\varepsilon(\sigma,\delta) := \mathbb{E} \sup_{g\in\mathcal{F}(\sigma)} \sup_{f\in\mathcal{F}(\delta),\rho_P(f,g)<\breve{r}(\sigma,\delta)+\varepsilon} |(P_n - P)(f-g)|.$$

Clearly, $\breve{\psi}_n^\varepsilon(\sigma,\delta) \downarrow \breve{\psi}_n(\sigma,\delta)$ as $\varepsilon \downarrow 0$. Define

$$\breve{U}_n^\varepsilon(\sigma;\delta;t) := \breve{\psi}_n^\varepsilon(\sigma,\delta) + \sqrt{2\frac{t}{n}((\breve{r}(\sigma,\delta)+\varepsilon)^2 + 2\breve{\psi}_n^\varepsilon(\sigma,\delta))} + \frac{t}{3n}.$$

We also have $\breve{U}_n^\varepsilon(\sigma;\delta;t) \downarrow \breve{U}_n(\sigma;\delta;t)$ as $\varepsilon \downarrow 0$. Let

$$E_{n,j}(t;\varepsilon) := \left\{\sup_{g\in\mathcal{F}(\sigma)} \sup_{f\in\mathcal{F}(\delta_j),\rho_P(f,g)<\breve{r}(\sigma,\delta_j)+\varepsilon} |(P_n-P)(f-g)| \leq \breve{U}_n^\varepsilon(\sigma,\delta_j;t)\right\}.$$

By Talagrand's concentration inequality, $\mathbb{P}((E_{n,j}(t;\varepsilon))^c) \leq e^{-t}$. Hence, for

$$E_n(t;\varepsilon) := \bigcap_{j\,:\,\delta_j \geq \delta} E_{n,j}(t;\varepsilon),$$

we have $\mathbb{P}((E_n(t;\varepsilon))^c) \leq \log_q \frac{q}{\delta} e^{-t}$. On the event $E_n(t;\varepsilon)$, for all $j$ such that $\delta_j \geq \delta$,

$$\begin{aligned} f \in \mathcal{F}(\delta_{j+1},\delta_j] &\implies \exists g \in \mathcal{F}(\sigma): \quad \rho_P(f,g) < \breve{r}(\sigma,\delta_j) + \varepsilon \\ &\implies \mathcal{E}(f) \leq Pf - Pg + \sigma \\ &\qquad \leq P_n f - P_n g + (P - P_n)(f-g) + \sigma \\ &\qquad \leq \hat{\mathcal{E}}_n(f) + \breve{U}_n^\varepsilon(\sigma,\delta_j;t) + \sigma. \end{aligned}$$

Therefore,

$$\mathbb{P}\{\exists j : \exists f \in \mathcal{F}(\delta_{j+1},\delta_j] : \delta_j \geq \delta,\ \mathcal{E}(f) > \hat{\mathcal{E}}_n(f) + \breve{U}_n^\varepsilon(\sigma,\delta_j;t) + \sigma\} \leq \log_q \frac{q}{\delta} e^{-t}.$$

Let

$$F := \{\exists f \in \mathcal{F} : \mathcal{E}(f) \geq \delta \text{ and } \hat{\mathcal{E}}_n(f) < (1 - q\breve{V}_n(\sigma,\delta;t))\mathcal{E}(f)\}.$$

Then,

$$\begin{aligned} F &\subset \{\exists j\ \exists f \in \mathcal{F}(\delta_{j+1},\delta_j] : \delta_j \geq \delta,\ \mathcal{E}(f) > \hat{\mathcal{E}}_n(f) + \breve{V}_n(\sigma,\delta;t)\delta_j\} \\ &\subset \{\exists j\ \exists f \in \mathcal{F}(\delta_{j+1},\delta_j] : \delta_j \geq \delta,\ \mathcal{E}(f) > \hat{\mathcal{E}}_n(f) + \breve{U}_n(\sigma,\delta_j;t) + \sigma\}. \end{aligned}$$

Because of the monotonicity of $\breve{U}_n^\varepsilon$ with respect to $\varepsilon$,

$$\mathbb{P}\{\exists j \ \exists f \in \mathcal{F}(\delta_{j+1}, \delta_j] : \delta_j \geq \delta, \ \mathcal{E}(f) > \hat{\mathcal{E}}_n(f) + \breve{U}_n(\sigma, \delta_j; t) + \sigma\}$$

$$= \lim_{\varepsilon \to 0} \mathbb{P}\{\exists j \ \exists f \in \mathcal{F}(\delta_{j+1}, \delta_j] : \delta_j \geq \delta, \ \mathcal{E}(f) > \hat{\mathcal{E}}_n(f) + \breve{U}_n^\varepsilon(\sigma, \delta_j; t) + \sigma\}$$

$$\leq \limsup_{\varepsilon \to 0} \mathbb{P}((E_n(t;\varepsilon))^c) \leq \log_q \frac{q}{\delta} e^{-t},$$

implying $\mathbb{P}(F) \leq \log_q \frac{q}{\delta} e^{-t}$. This proves the second bound of the theorem and it also implies the first bound since on the event $F^c$, $\mathcal{E}(\hat{f}_n) \leq \delta$; otherwise, we would have

$$0 = \hat{\mathcal{E}}_n(\hat{f}_n) \geq (1 - q\breve{V}_n(\sigma, \delta; t))\mathcal{E}(\hat{f}_n) \geq \delta/2,$$

a contradiction. □

PROOF OF PROPOSITION 2. We have $Pf = 1/2$ for all $f \in \mathcal{F}$ and as a result $\mathcal{F}(\delta) = \mathcal{F}$ for all $\delta \geq 0$. This implies $\forall 0 < \sigma \leq \delta : \breve{r}(\sigma; \delta) = 0$ and also $\breve{\psi}_n(\sigma; \delta) = 0$. Therefore, $\breve{\delta}_n(\sigma; t)$ is of the order $Ct/n$. Note also that $\forall k \neq j$: $P(f_k - f_j)^2 = 1/2$, so, $D_P(\mathcal{F}; \delta) = 1/2$. On the other hand,

$$\phi_n(\delta) = \mathbb{E} \sup_{f,g \in \mathcal{F}} |(P_n - P)(f - g)| = \mathbb{E} \max_{1 \leq k, j \leq N} |(P_n - P)(f_k - f_j)|,$$

which can be shown to be of the order $c(\log N/n)^{1/2}$. This easily yields the value of $\delta_n(t)$ of the order $c((\log N/n)^{1/2} + (t/n)^{1/2})$. The excess risk of $\hat{f}_n$ (and, as a matter of fact, of any $f \in \mathcal{F}$) is 0, so the bound $\delta_n(t)$ is not sharp at all. Next we show that (iv) also holds. To this end, note that

$$\mathbb{P}\{\mathcal{F}(0) \subset \hat{\mathcal{F}}_n(\delta)\} = \mathbb{P}\{\hat{\mathcal{F}}_n(\delta) = \mathcal{F}\}$$

$$= \mathbb{P}\left\{\forall j, 1 \leq j \leq N+1 : P_n f_j \leq \min_{1 \leq k \leq N+1} P_n f_k + \delta\right\}$$

$$\leq \mathbb{P}\{\forall j, 1 \leq j \leq N : P_n f_j \leq P_n f_{N+1} + \delta\}$$

$$= \mathbb{P}\{\forall j, 1 \leq j \leq N : \nu_{n,j} \leq \nu_n + \delta n\},$$

where $\nu_n, \nu_{n,j}, 1 \leq j \leq N$, are i.i.d. binomial random variables with parameters $n$ and $1/2$. Thus, we get

$$\mathbb{P}\{\mathcal{F}(0) \subset \hat{\mathcal{F}}_n(\delta)\} \leq \sum_{k=0}^{n} \mathbb{P}\{\nu_n = k\} \mathbb{P}\{\forall j, 1 \leq j \leq N : \nu_{n,j} \leq k + \delta n | \nu_n = k\}$$

$$= \sum_{k=0}^{n} \mathbb{P}\{\nu_n = k\} \prod_{j=1}^{N} \mathbb{P}\{\nu_{n,j} \leq k + \delta n\}$$

$$= \sum_{k=0}^{n} \mathbb{P}\{\nu_n = k\} \mathbb{P}^N\{\nu_n \leq k + \delta n\}$$

$$\leq \mathbb{P}\{\nu_n > \bar{k}\} + \mathbb{P}^N\{\nu_n \leq \bar{k} + \delta n\},$$

where $0 \leq \bar{k} \leq n$. Let $\bar{k} = \frac{n}{2} + n\delta$. Then, using Bernstein's inequality, we get

$$\mathbb{P}\{\nu_n > \bar{k}\} \leq \exp\left\{-\frac{n\delta^2}{4}\right\} = (\log N)^{-2^{-6}}.$$

On the other hand, using normal approximation of binomial distribution we get ($\Phi$ denoting the standard normal distribution function)

$$\mathbb{P}\{\nu_n \leq \bar{k} + \delta n\} \leq \Phi(4\delta\sqrt{n}) + n^{-1/2} = \Phi(\sqrt{\log N}) + n^{-1/2}.$$

Under the condition $N_0 \leq N \leq \sqrt{n}$ this easily gives (for a large enough $N_0$) $\mathbb{P}\{\mathcal{F}(0) \subset \hat{\mathcal{F}}_n(\delta)\} \leq \varepsilon$, which implies the claim. $\square$

PROOF OF LEMMA 3. First note that by Theorem 1 the event $\{\mathcal{E}(\hat{f}_n) \leq \bar{\delta}_n(t)\}$ holds with probability at least $1 - \log_q \frac{q}{\bar{\delta}_n(t)} e^{-t}$. On this event, we have for all $g \in \mathcal{F}(\varepsilon)$ with $\varepsilon < \bar{\delta}_n(t)$

$$\begin{aligned}
\left|\inf_{\mathcal{F}} P_n f - \inf_{\mathcal{F}} P f\right| &= \left|P_n \hat{f}_n - \inf_{\mathcal{F}} P f\right| \\
&\leq P\hat{f}_n - \inf_{\mathcal{F}} P f + |(P_n - P)(\hat{f}_n - g)| + |(P_n - P)(g)| \\
&\leq \bar{\delta}_n(t) + \sup_{f,g \in \mathcal{F}(\bar{\delta}_n(t))} |(P_n - P)(f - g)| + |(P_n - P)(g)|.
\end{aligned}$$
(9.3)

By Talagrand's inequality with probability at least $1 - e^{-t}$

$$\sup_{f,g \in \mathcal{F}(\bar{\delta}_n(t))} |(P_n - P)(f - g)| \leq \bar{U}_n(\bar{\delta}_n(t); t) \leq q\bar{V}_n(\bar{\delta}_n(t); t)\bar{\delta}_n(t) \leq \bar{\delta}_n(t).$$
(9.4)

On the other hand, by Bernstein's inequality, also with probability at least $1 - e^{-t}$

$$(9.5) \quad |(P_n - P)(g)| \leq \sqrt{2\frac{t}{n} \operatorname{Var}_P g} + \frac{2t}{3n} \leq \sqrt{2\frac{t}{n}\left(\inf_{\mathcal{F}} Pf + \varepsilon\right)} + \frac{2t}{3n},$$

since $g$ takes values in $[0,1]$, $g \in \mathcal{F}(\varepsilon)$, and hence $\operatorname{Var}_P g \leq Pg^2 \leq Pg \leq \inf_{\mathcal{F}} Pf + \varepsilon$. It follows from (9.3), (9.4) and (9.5) that on some event $E(\varepsilon)$ with probability at least $1 - \log_q \frac{q^3}{\bar{\delta}_n(t)} e^{-t}$ the following inequality holds:

$$(9.6) \quad \left|\inf_{\mathcal{F}} P_n f - \inf_{\mathcal{F}} P f\right| \leq 2\bar{\delta}_n(t) + \sqrt{2\frac{t}{n}\left(\inf_{\mathcal{F}} Pf + \varepsilon\right)} + \frac{t}{n}.$$

Since the events $E(\varepsilon)$ are monotone in $\varepsilon$, one can let $\varepsilon \to 0$ which yields the first bound of the lemma.

To prove the second bound, note that on the same event on which (9.6) with $\varepsilon = 0$ holds we also have

$$\left|\inf_{\mathcal{F}} P_n f - \inf_{\mathcal{F}} Pf\right| \leq \sqrt{2\frac{t}{n}\left|\inf_{\mathcal{F}} P_n f - \inf_{\mathcal{F}} Pf\right|} + 2\bar{\delta}_n(t) + \sqrt{2\frac{t}{n}\inf_{\mathcal{F}} P_n f} + \frac{t}{n}.$$
(9.7)

We either have

$$\left|\inf_{\mathcal{F}} P_n f - \inf_{\mathcal{F}} Pf\right| \leq \frac{8t}{n} \quad \text{or} \quad \frac{2t}{n} \leq \frac{|\inf_{\mathcal{F}} P_n f - \inf_{\mathcal{F}} Pf|}{4},$$

and in the last case (9.7) implies that

$$\left|\inf_{\mathcal{F}} P_n f - \inf_{\mathcal{F}} Pf\right| \leq 4\bar{\delta}_n(t) + 2\sqrt{2\frac{t}{n}\inf_{\mathcal{F}} P_n f} + \frac{2t}{n}.$$

We can use now the condition of the lemma to replace $\bar{\delta}_n(t)$ by $\hat{\delta}_n(t)$ and to get that with probability at least $1 - p - \log_q \frac{q^3}{\bar{\delta}_n(t)} e^{-t}$ the following bound holds:

$$\left|\inf_{\mathcal{F}} P_n f - \inf_{\mathcal{F}} Pf\right| \leq 4\hat{\delta}_n(t) + 2\sqrt{2\frac{t}{n}\inf_{\mathcal{F}} P_n f} + \frac{8t}{n}. \qquad \Box$$

PROOF OF THEOREM 5. We will use the following consequence of Theorem 1 and of Lemma 3 (and its proof): there exists an event $E$ of probability at least

$$1 - \sum_{k=1}^{\infty}\left(p_k + \log_q \frac{q^3 n}{t_k} e^{-t_k}\right)$$

such that on the event $E$, $\forall k \geq 1$:

$$P\hat{f}_k - \inf_{f \in \mathcal{F}_k} Pf \leq \bar{\delta}_n(\mathcal{F}_k; t_k) \leq \hat{\delta}_n(\mathcal{F}_k; t_k) \leq \tilde{\delta}_n(\mathcal{F}_k; t_k)$$

and

$$\left|\inf_{\mathcal{F}_k} P_n f - \inf_{\mathcal{F}_k} Pf\right| \leq 2\bar{\delta}_n(\mathcal{F}_k; t_k) + \sqrt{\frac{2t_k}{n}\inf_{\mathcal{F}_k} Pf} + \frac{t_k}{n},$$

$$\left|\inf_{\mathcal{F}_k} P_n f - \inf_{\mathcal{F}_k} Pf\right| \leq 4\hat{\delta}_n(\mathcal{F}_k; t_k) + 2\sqrt{\frac{2t_k}{n}\inf_{\mathcal{F}_k} P_n f} + \frac{8t_k}{n}.$$

Note also that the events involved in the proof of Lemma 3 are the same that are involved in the bound of Theorem 1; because of this reason, we do not have to add probabilities here. On the event $E$, we have

$$P\hat{f} = P\hat{f}_{\hat{k}} \leq \inf_{\mathcal{F}_{\hat{k}}} Pf + \bar{\delta}_n(\mathcal{F}_{\hat{k}}; t_{\hat{k}})$$

$$\leq \inf_{\mathcal{F}_{\hat{k}}} P_n f + 5\hat{\delta}_n(\mathcal{F}_{\hat{k}}; t_{\hat{k}}) + 2\sqrt{\frac{2t_{\hat{k}}}{n} \inf_{\mathcal{F}_{\hat{k}}} P_n f} + \frac{8t_{\hat{k}}}{n}$$

$$\leq \inf_{\mathcal{F}_{\hat{k}}} P_n f + \hat{\pi}(\hat{k}) = \inf_k \left[ \inf_{\mathcal{F}_k} P_n f + \hat{\pi}(k) \right],$$

provided that the constant $\hat{K}$ in the definition of $\hat{\pi}$ was chosen properly. This proves the first bound of the theorem.

To prove the second bound, note that since

$$\sqrt{\frac{t_k}{n} \inf_{\mathcal{F}_k} P_n f} \leq \sqrt{\frac{t_k}{n} \inf_{\mathcal{F}_k} Pf} + \sqrt{\frac{t_k}{n} \left| \inf_{\mathcal{F}_k} P_n f - \inf_{\mathcal{F}_k} Pf \right|}$$

$$\leq \sqrt{\frac{t_k}{n} \inf_{\mathcal{F}_k} Pf} + \frac{t_k}{2n} + \frac{1}{2} \left| \inf_{\mathcal{F}_k} P_n f - \inf_{\mathcal{F}_k} Pf \right|,$$

we also have on the event $E$ for all $k$

$$\hat{\pi}(k) = \hat{K} \left[ \hat{\delta}_n(\mathcal{F}_k; t_k) + \sqrt{\frac{t_k}{n} \inf_{\mathcal{F}_k} P_n f} + \frac{t_k}{n} \right]$$

$$\leq \frac{\tilde{K}}{2} \left[ \tilde{\delta}_n(\mathcal{F}_k; t_k) + \sqrt{\frac{t_k}{n} \inf_{\mathcal{F}_k} Pf} + \frac{t_k}{n} \right] = \tilde{\pi}(k)/2$$

and

$$\left| \inf_{\mathcal{F}_k} P_n f - \inf_{\mathcal{F}_k} Pf \right| \leq 2\bar{\delta}_n(\mathcal{F}_k; t_k) + \sqrt{\frac{2t_k}{n} \inf_{\mathcal{F}_k} Pf} + \frac{t_k}{3n}$$

$$\leq \frac{\tilde{K}}{2} \left[ \tilde{\delta}_n(\mathcal{F}_k; t_k) + \sqrt{\frac{t_k}{n} \inf_{\mathcal{F}_k} Pf} + \frac{t_k}{n} \right] = \tilde{\pi}(k)/2,$$

provided that the constant $\tilde{K}$ in the definition of $\tilde{\pi}(k)$ was chosen to be large enough. This yields on the event $E$

$$P\hat{f} \leq \inf_k \left[ \inf_{\mathcal{F}_k} P_n f + \hat{\pi}(k) \right] \leq \inf_k \left[ \inf_{\mathcal{F}_k} Pf + \tilde{\pi}(k) \right],$$

proving the second bound. □

PROOF OF LEMMA 4. We assume, for simplicity, that $Pf$ attains its minimum over $\mathcal{G}$ at some $\bar{f} \in \mathcal{G}$ (the proof can be easily modified if the minimum is not attained). Let $E$ be the event such that the following inequalities hold:

$$|(P_n - P)(\bar{f} - f_*)| \leq \sqrt{\frac{2t}{n} \mathrm{Var}_P(\bar{f} - f_*)} + \frac{t}{n} \quad \text{and}$$

$$\forall f \in \mathcal{G}: \qquad \hat{\mathcal{E}}_n(\mathcal{G};f) \leq \frac{3}{2}(\mathcal{E}_P(\mathcal{G};f) \vee \bar{\delta}_n(\mathcal{G};t)).$$

The first of these inequalities holds with probability at least $1 - e^{-t}$ by Bernstein's inequality; the second inequality takes place with probability at least $1 - \log_q \frac{q^2 n}{t} e^{-t}$ by (9.2) in the proof of Lemma 2. Hence, $\mathbb{P}(E) \geq 1 - \log_q \frac{q^3 n}{t} e^{-t}$. We also have $\operatorname{Var}_P^{1/2}(\bar{f} - f_*) \leq \varphi^{-1}(P\bar{f} - Pf_*)$ and hence, on the event $E$,

$$|(P - P_n)(\bar{f} - f_*)| \leq \varphi(\sqrt{\varepsilon}\varphi^{-1}(P\bar{f} - Pf_*)) + \varphi^*\left(\sqrt{\frac{2t}{\varepsilon n}}\right) + \frac{t}{n}$$

$$\leq \varphi(\sqrt{\varepsilon})(P\bar{f} - Pf_*) + \varphi^*\left(\sqrt{\frac{2t}{\varepsilon n}}\right) + \frac{t}{n},$$

implying

$$(9.8) \qquad P_n(\bar{f} - f_*) \leq (1 + \varphi(\sqrt{\varepsilon}))P(\bar{f} - f_*) + \varphi^*\left(\sqrt{\frac{2t}{\varepsilon n}}\right) + \frac{t}{n}$$

and

$$(9.9) \qquad P(\bar{f} - f_*) \leq (1 - \varphi(\sqrt{\varepsilon}))^{-1}\left[P_n(\bar{f} - f_*) + \varphi^*\left(\sqrt{\frac{2t}{\varepsilon n}}\right) + \frac{t}{n}\right].$$

Equation (9.8) immediately yields the first bound of the lemma. Since on the event $E$

$$P_n(\bar{f} - f_*) = P_n\bar{f} - \inf_{\mathcal{G}} P_n f + \inf_{\mathcal{G}} P_n f - P_n f_* = \hat{\mathcal{E}}_n(\mathcal{G};\bar{f}) + \inf_{\mathcal{G}} P_n f - P_n f_*$$

$$\leq \inf_{\mathcal{G}} P_n f - P_n f_* + \frac{3}{2}(\mathcal{E}_P(\mathcal{G};\bar{f}) \vee \bar{\delta}_n(\mathcal{G};t)),$$

and since $\mathcal{E}_P(\mathcal{G};\bar{f}) = 0$, we get

$$P_n(\bar{f} - f_*) \leq \inf_{\mathcal{G}} P_n f - P_n f_* + \frac{3}{2}\bar{\delta}_n(\mathcal{G};t).$$

Along with (9.9), this implies

$$\inf_{\mathcal{G}} Pf - Pf_* = P(\bar{f} - f_*) \leq (1 - \varphi(\sqrt{\varepsilon}))^{-1}\left[\inf_{\mathcal{G}} P_n f - P_n f_* + \frac{3}{2}\bar{\delta}_n(\mathcal{G};t)\right.$$

$$\left. + \varphi^*\left(\sqrt{\frac{2t}{\varepsilon n}}\right) + \frac{t}{n}\right],$$

which is the second bound of the lemma.

Finally, to prove the third bound plug into (5.5) the bound on $\bar{\delta}_n(\mathcal{G};t)$ and solve the resulting inequality with respect to $\inf_{\mathcal{G}} Pf - Pf_*$. □

PROOF OF THEOREM 6. Let $E_k$ be the event defined in Lemma 4 for $\mathcal{G} = \mathcal{F}_k$ and $t = t_k$. Let $E$ be the event such that the following inequalities and events $E_k$ hold for all $k$:

$$\mathcal{E}_P(\mathcal{F}_k; \hat{f}_k) = P\hat{f}_k - \inf_{\mathcal{F}_k} Pf \leq \bar{\delta}_n(\mathcal{F}_k; t_k)$$

and $\bar{\delta}_n(\mathcal{F}_k; t_k) \leq \hat{\bar{\delta}}_n(\mathcal{F}_k; t_k) \leq \tilde{\delta}_n(\mathcal{F}_k; t_k)$. The first of the inequalities holds with probability at least $1 - \log_q \frac{qn}{t_k} e^{-t_k}$ either by Theorem 1 or by Theorem 4; the second one holds with probability at least $1 - p_k$ by assumptions. Therefore, using Lemma 4,

$$\mathbb{P}(E) \geq 1 - \sum_{k=1}^{\infty} \left( p_k + 2\log_q \frac{q^2 n}{t_k} e^{-t_k} \right).$$

On the event $E$, using first bound (5.5) and then (5.4) of Lemma 4, we get

$$\mathcal{E}_P(\mathcal{F}; \hat{f}) = P\hat{f} - \inf_{\mathcal{F}} Pf = P\hat{f}_{\hat{k}} - Pf_* = P\hat{f}_{\hat{k}} - \inf_{\mathcal{F}_{\hat{k}}} Pf + \inf_{\mathcal{F}_{\hat{k}}} Pf - Pf_*$$

$$\leq \bar{\delta}_n(\mathcal{F}_{\hat{k}}; t_{\hat{k}}) + \inf_{\mathcal{F}_{\hat{k}}} Pf - Pf_*$$

$$\leq (1 - \varphi(\sqrt{\varepsilon}))^{-1} \left[ (1 - \varphi(\sqrt{\varepsilon}))\bar{\delta}_n(\mathcal{F}_{\hat{k}}; t_{\hat{k}}) + \inf_{\mathcal{F}_{\hat{k}}} P_n f - P_n f_* \right.$$

$$\left. + \frac{3}{2}\bar{\delta}_n(\mathcal{F}_{\hat{k}}; t_{\hat{k}}) + \varphi^*\left(\sqrt{\frac{2t_{\hat{k}}}{\varepsilon n}}\right) + \frac{t_{\hat{k}}}{n} \right]$$

$$\leq (1 - \varphi(\sqrt{\varepsilon}))^{-1} \left\{ \inf_k \left[ \inf_{\mathcal{F}_k} P_n f + (5/2 - \varphi(\sqrt{\varepsilon}))\hat{\bar{\delta}}_n(\mathcal{F}_k; t_k) \right. \right.$$

$$\left. \left. + \varphi^*\left(\sqrt{\frac{2t_k}{\varepsilon n}}\right) + \frac{t_k}{n} \right] - P_n f_* \right\}$$

$$= (1 - \varphi(\sqrt{\varepsilon}))^{-1} \left\{ \inf_k \left[ \inf_{\mathcal{F}_k} P_n f + \hat{\pi}(k) \right] - P_n f_* \right\}$$

$$\leq \frac{1 + \varphi(\sqrt{\varepsilon})}{1 - \varphi(\sqrt{\varepsilon})} \inf_k \left[ \inf_{\mathcal{F}_k} Pf - \inf_{\mathcal{F}} Pf + \frac{5/2 - \varphi(\sqrt{\varepsilon})}{1 + \varphi(\sqrt{\varepsilon})} \tilde{\delta}_n(\mathcal{F}_k; t_k) \right.$$

$$\left. + \frac{2}{1 + \varphi(\sqrt{\varepsilon})} \varphi^*\left(\sqrt{\frac{2t_k}{\varepsilon n}}\right) + \frac{2}{(1 + \varphi(\sqrt{\varepsilon}))} \frac{t_k}{n} \right]$$

$$= \inf_k \frac{1 + \varphi(\sqrt{\varepsilon})}{1 - \varphi(\sqrt{\varepsilon})} \left[ \inf_{\mathcal{F}_k} Pf - \inf_{\mathcal{F}} Pf + \tilde{\pi}(k) \right],$$

and the result follows. □

PROOF OF THEOREM 7. Let us define the event $E$ such that on this event $\forall l$ and $\forall k \leq l$

$$\inf_{f \in \mathcal{F}_k} \hat{\mathcal{E}}_n(\mathcal{F}_l, f) \leq 2\Big(\inf_{f \in \mathcal{F}_k} \mathcal{E}_P(\mathcal{F}_l, f) \vee \bar{\delta}_n(\mathcal{F}_l, t_l)\Big), \tag{9.10}$$

$$\inf_{f \in \mathcal{F}_k} \mathcal{E}_P(\mathcal{F}_l, f) \leq 2 \inf_{f \in \mathcal{F}_k} \hat{\mathcal{E}}_n(\mathcal{F}_l, f) \vee \bar{\delta}_n(\mathcal{F}_l, t_l), \tag{9.11}$$

and

$$\bar{\delta}_n(\mathcal{F}_l; t_l) \leq \hat{\bar{\delta}}_n(\mathcal{F}_l; t_l) \leq \tilde{\delta}_n(\mathcal{F}_l; t_l). \tag{9.12}$$

Then we have

$$\mathbb{P}(E) \geq 1 - \sum_{k=1}^{\infty}\bigg(p_k + \log_q \frac{q^2 n}{t_k} e^{-t_k}\bigg),$$

which is true because of the following reasons. First, for any $l$, we have with probability at least $1 - \log_q \frac{q^2}{\bar{\delta}_n(\mathcal{F}_l, t_l)} e^{-t_l}$ that for all $f \in \mathcal{F}_l$

$$\hat{\mathcal{E}}_n(\mathcal{F}_l, f) \leq 2(\mathcal{E}_P(\mathcal{F}_l, f) \vee \bar{\delta}_n(\mathcal{F}_l, t_l)) \quad \text{and} \quad \mathcal{E}_P(\mathcal{F}_l, f) \leq 2\hat{\mathcal{E}}_n(\mathcal{F}_l, f) \vee \bar{\delta}_n(\mathcal{F}_l, t_l)$$

[see the proof of Lemma 2, specifically, (9.1), (9.2)]. Then, by assumptions, for all $l$ with probability at least $1 - p_l$, $\bar{\delta}_n(\mathcal{F}_l; t_l) \leq \hat{\bar{\delta}}_n(\mathcal{F}_l; t_l) \leq \tilde{\delta}_n(\mathcal{F}_l; t_l)$. It remains to use the union bound to get the above lower bound on $\mathbb{P}(E)$.

Clearly, on the event $E$, $\forall l : \bar{\delta}_n(l) \leq \hat{\delta}_n(l) \leq \tilde{\delta}_n(l)$. We will show that on the same event $E$, $\tilde{k} \leq \hat{k} \leq \bar{k} \leq k^*$. The inequality $\bar{k} \leq k^*$ is obvious from the definitions. If $k < \hat{k}$, then there exists $l > k$ such that

$$\inf_{\mathcal{F}_k} \hat{\mathcal{E}}_n(\mathcal{F}_l, f) = \inf_{\mathcal{F}_k} P_n f - \inf_{\mathcal{F}_l} P_n f > \hat{c}\hat{\delta}_n(l).$$

We will use that, due to (9.10), on the event $E$

$$\inf_{\mathcal{F}_k} \hat{\mathcal{E}}_n(\mathcal{F}_l, f) \leq 2\Big(\inf_{\mathcal{F}_k} \mathcal{E}_P(\mathcal{F}_l, f) \vee \bar{\delta}_n(l)\Big).$$

Therefore (assuming that the constants $\hat{c}, \bar{c}$ have been chosen properly)

$$\inf_{\mathcal{F}_k} Pf - \inf_{\mathcal{F}_l} Pf = \inf_{\mathcal{F}_k} \mathcal{E}_P(\mathcal{F}_l, f) \geq \frac{\hat{c}}{2}\hat{\delta}_n(l) - \bar{\delta}_n(l) \geq \bigg(\frac{\hat{c}}{2} - 1\bigg)\bar{\delta}_n(l) \geq \bar{c}\bar{\delta}_n(l),$$

which implies that $k < \bar{k}$ and hence $\hat{k} \leq \bar{k}$. Similarly, if $k < \tilde{k}$, then there exists $l > k$ such that

$$\inf_{\mathcal{F}_k} \mathcal{E}_P(\mathcal{F}_l, f) = \inf_{\mathcal{F}_k} Pf - \inf_{\mathcal{F}_l} Pf > \tilde{c}\tilde{\delta}_n(l).$$

Due to (9.11), on the event $E$

$$\inf_{\mathcal{F}_k} \mathcal{E}_P(\mathcal{F}_l, f) \leq 2\inf_{\mathcal{F}_k} \hat{\mathcal{E}}_n(\mathcal{F}_l, f) \vee \bar{\delta}_n(l),$$

implying that

$$\inf_{\mathcal{F}_k} P_n f - \inf_{\mathcal{F}_l} P_n f = \inf_{\mathcal{F}_k} \hat{\mathcal{E}}_n(\mathcal{F}_l, f) \geq (\tilde{c}\tilde{\delta}_n(l) - \bar{\delta}_n(l))/2 \geq \left(\frac{\tilde{c}-1}{2}\right)\tilde{\delta}_n(l) > \hat{c}\hat{\delta}_n(l),$$

provided that the constants have been chosen properly. Therefore, $k < \hat{k}$ and hence $\tilde{k} \leq \hat{k}$.

Next we have on the event $E$ for all $k \geq \bar{k}$

$$P\hat{f} - \inf_j \inf_{\mathcal{F}_j} Pf = P\hat{f}_{\hat{k}} - \inf_{\mathcal{F}_k} Pf + \inf_{\mathcal{F}_k} Pf - \inf_j \inf_{\mathcal{F}_j} Pf$$

$$= P\hat{f}_{\hat{k}} - \inf_{\mathcal{F}_{\hat{k}}} Pf + \inf_{\mathcal{F}_{\hat{k}}} Pf - \inf_{\mathcal{F}_k} Pf + \inf_{\mathcal{F}_k} Pf - \inf_j \inf_{\mathcal{F}_j} Pf$$

$$\leq \bar{\delta}_n(\hat{k}) + \inf_{\mathcal{F}_{\hat{k}}} Pf - \inf_{\mathcal{F}_k} Pf + \inf_{\mathcal{F}_k} Pf - \inf_j \inf_{\mathcal{F}_j} Pf$$

$$\leq \bar{\delta}_n(\hat{k}) + \tilde{c}\tilde{\delta}_n(k) + \inf_{\mathcal{F}_k} Pf - \inf_j \inf_{\mathcal{F}_j} Pf$$

$$\leq \inf_{\mathcal{F}_k} Pf - \inf_j \inf_{\mathcal{F}_j} Pf + (\tilde{c}+1)\tilde{\delta}_n(k),$$

implying the first bound. The second bound follows immediately by plugging in $k = k^*$ (which is possible since $k^* \geq \bar{k}$) and observing that $\inf_{\mathcal{F}_{k^*}} Pf - \inf_j \inf_{\mathcal{F}_j} Pf = 0$. □

PROOF OF THEOREM 8. Since $\phi_n(\delta) \leq \omega_n(D(\delta))$, conditions (i) and (ii) imply that, for all $P \in \mathcal{P}_{\rho,\kappa,C}(\mathcal{F})$, $\phi_n(\delta) \leq Kn^{-1/2}\delta^{\frac{1-\rho}{2\kappa}}$. Then, by an easy computation,

$$\bar{\delta}_n(t) \leq K\left[\left(\frac{1}{n}\right)^{\frac{\kappa}{2\kappa+\rho-1}} \vee \left(\frac{t}{n}\right)^{\frac{\kappa}{2\kappa-1}} \vee \frac{t}{n}\right]$$

with some $K > 0$. It remains to recall that $\bar{\delta}_n(t) \geq \delta_n(t)$ and to use Theorem 1 with $t$ replaced by $t + \log\log_q n$ to get with some $K > 0$ for all $P \in \mathcal{P}_{\rho,\kappa,C}(\mathcal{F})$, the bound

$$\mathbb{P}\{n^{\frac{\kappa}{2\kappa+\rho-1}}\mathcal{E}(\hat{f}_n) \geq K(1+t)\} \leq e^{-t},$$

which implies the result. □

PROOF OF THEOREM 9. We use Theorem 7 to get for all $P$

$$\mathbb{P}\{P\hat{f} - Pf_* \geq K\tilde{\delta}_n(k^*(P))\} = O(n^{-2}).$$

Since for all $P \in \mathcal{P}_j$, $k^*(P) = j$, we have

$$\max_{1\leq j\leq N} \sup_{P\in\mathcal{P}_j} \mathbb{P}\{P\hat{f} - Pf_* \geq K\tilde{\delta}_n(j)\} = O(n^{-2}).$$

The same argument as in the proof of Theorem 8 shows that $\tilde{\delta}_n(j) \leq Kn^{-\beta_j}$. Therefore

$$\max_{1\leq j\leq N} \sup_{P\in\mathcal{P}_j} n^{\beta_j}\mathbb{E}(P\hat{f} - Pf_*) \leq \max_{1\leq j\leq N} n^{\beta_j} \sup_{P\in\mathcal{P}_j} \mathbb{P}\{P\hat{f} - Pf_* \geq Kn^{-\beta_j}\} + K$$

$$\leq K + O\Big(\max_{1\leq j\leq N} n^{\beta_j - 2}\Big) = O(1). \qquad \square$$

PROOF OF THEOREM 10. We first look at a single class $\mathcal{F}$ of binary functions. The following upper bounds hold:

$$D^2(\mathcal{F};\delta) = \sup_{f,g\in\mathcal{F}(\delta)} P(f-g)^2 \leq \sup_{f,g\in\mathcal{F}(\delta)} (Pf + Pg) \leq 2\Big(\inf_{f\in\mathcal{F}} Pf + \delta\Big)$$

and

$$(9.13) \quad \omega_n(\mathcal{F};\delta) \leq K\left[\delta\sqrt{\frac{\mathbb{E}\log \Delta^{\mathcal{F}}(X_1,\ldots,X_n)}{n}} + \frac{\mathbb{E}\log \Delta^{\mathcal{F}}(X_1,\ldots,X_n)}{n}\right],$$

where the proof of the second bound can be found in [36]. It follows that

$$\phi_n(\delta) \leq K\left[\sqrt{2\Big(\inf_{f\in\mathcal{F}} Pf + \delta\Big)\frac{\mathbb{E}\log \Delta^{\mathcal{F}}(X_1,\ldots,X_n)}{n}} + \frac{\mathbb{E}\log \Delta^{\mathcal{F}}(X_1,\ldots,X_n)}{n}\right],$$

which implies, by using the $\sharp$-transform, that with some constant $K$

$$\bar{\delta}_n(t) \leq K\left[\sqrt{\inf_{f\in\mathcal{F}} Pf \frac{\mathbb{E}\log \Delta^{\mathcal{F}}(X_1,\ldots,X_n) + t}{n}} + \frac{\mathbb{E}\log \Delta^{\mathcal{F}}(X_1,\ldots,X_n) + t}{n}\right].$$

We now define

$$\hat{\delta}_n(t) := \hat{K}\left[\sqrt{\inf_{f\in\mathcal{F}} P_n f \frac{\log \Delta^{\mathcal{F}}(X_1,\ldots,X_n) + t}{n}} + \frac{\log \Delta^{\mathcal{F}}(X_1,\ldots,X_n) + t}{n}\right]$$

and

$$\tilde{\delta}_n(t) := \tilde{K}\left[\sqrt{\inf_{f\in\mathcal{F}} Pf \frac{\mathbb{E}\log \Delta^{\mathcal{F}}(X_1,\ldots,X_n) + t}{n}} + \frac{\mathbb{E}\log \Delta^{\mathcal{F}}(X_1,\ldots,X_n) + t}{n}\right].$$

We use the following deviation inequality for shattering numbers due to Boucheron, Lugosi and Massart [12]: with probability at least $1 - e^{-t}$

$$\log \Delta^{\mathcal{F}}(X_1,\ldots,X_n) \leq 2\mathbb{E}\log \Delta^{\mathcal{F}}(X_1,\ldots,X_n) + 2t$$

and

$$\mathbb{E}\log \Delta^{\mathcal{F}}(X_1,\ldots,X_n) \leq 2\log \Delta^{\mathcal{F}}(X_1,\ldots,X_n) + 2t.$$

Using this device together with Lemma 3, it is easy to see that with probability at least $1 - \log_q \frac{q^3 n}{t} e^{-t}$ we have $\bar{\delta}_n(t) \leq \hat{\delta}_n(t) \leq \tilde{\delta}_n(t)$. For instance, to prove the first of the two inequalities, note that, by the above deviation inequality for shattering numbers, on an event of probability at least $1-e^{-t}$ we can replace in the bound on $\bar{\delta}_n(t)$ $\mathbb{E}\log\Delta^{\mathcal{F}}(X_1,\ldots,X_n)$ by $\log\Delta^{\mathcal{F}}(X_1,\ldots,X_n)$. On the other hand, the first bound of Lemma 3 implies that with probability at least $1-\log_q \frac{q^3}{\bar{\delta}_n(t)}e^{-t}$ we have (using $2ab \leq a^2+b^2$)

$$\inf_{\mathcal{F}} Pf \leq \inf_{\mathcal{F}} P_n f + 2\bar{\delta}_n(t) + 2\sqrt{\frac{t}{n}\inf_{\mathcal{F}} Pf/2} + \frac{t}{3n}$$

$$\leq \inf_{\mathcal{F}} P_n f + 2\bar{\delta}_n(t) + \inf_{\mathcal{F}} Pf/2 + \frac{2t}{n},$$

which implies $\inf_{\mathcal{F}} Pf \leq 2\inf_{\mathcal{F}} P_n f + 4\bar{\delta}_n(t) + 4t/n$. Plugging this into the bound on $\bar{\delta}_n(t)$ and replacing $\mathbb{E}\log\Delta^{\mathcal{F}}(X_1,\ldots,X_n)$ by $\log\Delta^{\mathcal{F}}(X_1,\ldots,X_n)$, we easily get (with some constant $K$)

$$\bar{\delta}_n(t) \leq K\left[\sqrt{\inf_{f\in\mathcal{F}} P_n f \frac{\log\Delta^{\mathcal{F}}(X_1,\ldots,X_n)+t}{n}} + \frac{\log\Delta^{\mathcal{F}}(X_1,\ldots,X_n)+t}{n}\right]$$

$$+ 2\sqrt{\frac{\bar{\delta}_n(t)}{2}\frac{K^2\log\Delta^{\mathcal{F}}(X_1,\ldots,X_n)+t}{2n}},$$

which, again using $2ab \leq a^2+b^2$, leads to the following bound (with some $\hat{K}$):

$$\bar{\delta}_n(t) \leq \hat{K}\left[\sqrt{\inf_{f\in\mathcal{F}} P_n f \frac{\log\Delta^{\mathcal{F}}(X_1,\ldots,X_n)+t}{n}} + \frac{\log\Delta^{\mathcal{F}}(X_1,\ldots,X_n)+t}{n}\right]$$

$$= \hat{\delta}_n(t),$$

which holds with probability at least $1-\log_q \frac{q^4}{\bar{\delta}_n(t)}e^{-t}$. The second inequality $\hat{\delta}_n(t) \leq \tilde{\delta}_n(t)$ can be proved similarly. For a sequence $\mathcal{F}_k$ of classes of binary functions, this gives condition (5.2) and allows us to use Theorem 5 to complete the proof. $\square$

PROOF OF LEMMA 5. First note that

$$\phi_n(\delta) = \mathbb{E}\sup_{f,g\in\mathcal{F}(\delta)}|(P_n-P)(f-g)| \leq 2\mathbb{E}\sup_{f\in\mathcal{F}(\delta)}|(P_n-P)(f-\bar{f})|.$$

Also, $f \in \mathcal{F}(\delta)$ implies that

$$\rho_P(f,\bar{f}) \leq \rho_P(f,f_*) + \rho_P(\bar{f},f_*) \leq \sqrt{D(Pf-Pf_*)} + \sqrt{D(P\bar{f}-Pf_*)}$$

$$\leq \sqrt{D(Pf-P\bar{f})} + 2\sqrt{D(P\bar{f}-Pf_*)}$$

$$\leq \sqrt{D\delta} + 2\sqrt{D\Delta} \leq \sqrt{2D(\delta+4\Delta)},$$

where $\Delta := P\bar{f} - Pf_* = \inf_{\mathcal{F}} Pf - Pf_*$. It follows that

$$D(\mathcal{F};\delta) \leq 2\sqrt{D}(\sqrt{\delta} + 2\sqrt{\Delta}) \quad \text{and} \quad \phi_n(\delta) \leq 2\theta_n(\sqrt{2D(\delta + 4\Delta)}).$$

As a consequence, recalling the definition of $\bar{U}_n(\delta;t)$, we easily get with some constant $C > 0$ for all $\varepsilon \in (0,1]$

$$\bar{U}_n(\delta;t) \leq C\theta_n(\sqrt{2D(\delta + 4\Delta)}) + C\sqrt{\frac{D\delta t}{n}} + C\left(\varepsilon\Delta + \frac{Dt}{n\varepsilon}\right)$$
$$=: \psi_1(\delta) + \psi_2(\delta) + \psi_3(\delta),$$

where we used the inequality $2\sqrt{D\Delta \frac{t}{n}} \leq \varepsilon\Delta + \frac{Dt}{n\varepsilon}$ to bound the term $D(\mathcal{F};\delta)\sqrt{\frac{t}{n}}$ involved in $\bar{U}_n(\delta;t)$. Since

$$\bar{\delta}_n(\mathcal{F};t) := \bar{U}_{n,t}^{\sharp,q}\left(\frac{1}{2q^3}\right) \leq \bar{U}_{n,t}^{\sharp}\left(\frac{1}{2q^3}\right),$$

it is enough now to bound the $\sharp$-transform of $\psi_1, \psi_2, \psi_3$ separately and to use property 4 of Section 2.3. Let $u := \frac{1}{6q^3}$. Then, by properties 3, 7, 8 of Section 2.3

$$\psi_1^{\sharp}(u) \leq \frac{1}{2D}\theta_n^{\sharp}\left(\frac{\varepsilon u}{4CD}\right) + 4\varepsilon\Delta.$$

Also, (see property 6 with $\alpha = 1/2$ and property 3) $\psi_2^{\sharp}(u) \leq C^2 Dt/(nu^2)$ and (property 5)

$$\psi_3(u) \leq \frac{C}{u}\left(\varepsilon\Delta + \frac{Dt}{n\varepsilon}\right).$$

As a result, property 4 now yields

$$\bar{\delta}_n(\mathcal{F};t) \leq \frac{1}{2D}\theta_n^{\sharp}\left(\frac{\varepsilon u}{4CD}\right) + \left(4 + \frac{C}{u}\right)\varepsilon\Delta + \left(\frac{C}{u} + \frac{C^2}{u^2}\right)\frac{Dt}{n\varepsilon},$$

which after proper rescaling of $\varepsilon$ and adjusting the constants gives the bound of the lemma. □

PROOF OF THEOREM 11. It is a straightforward consequence of Theorem 6, Remarks 2 and 4 after this theorem and Lemma 5. Note that one should choose $\varphi_k(u) = u^2/D_k$, which implies that $\varphi^*(v) = D_k v^2/4$. The rest is an easy computation. □

PROOF OF LEMMA 6. First of all, note that by Lipschitz condition (7.3) $\forall g_1, g_2 \in \mathcal{G}$

$$P|(\ell \bullet g_1) - (\ell \bullet g_2)|^2 \leq L^2\|g_1 - g_2\|_{L_2(\Pi)}^2.$$

Next, by (7.5), we have for $g \in \mathcal{G}$, $x \in S$, $y \in T$

$$\frac{\ell(y, g(x)) + \ell(y, \bar{g}(x))}{2} \geq \ell\left(y; \frac{g(x) + \bar{g}(x)}{2}\right) + \psi(|g(x) - \bar{g}(x)|^r).$$

Integrating this inequality and observing that $\frac{g+\bar{g}}{2} \in \mathcal{G}$ and hence $P(\ell \bullet (\frac{g+\bar{g}}{2})) \geq P(\ell \bullet \bar{g})$ yields

$$\frac{P(\ell \bullet g) + P(\ell \bullet \bar{g})}{2} \geq P(\ell \bullet \bar{g}) + \Pi \psi(|g - \bar{g}|^r),$$

or

$$P(\ell \bullet g) - P(\ell \bullet \bar{g}) \geq 2\Pi \psi(|g - \bar{g}|^r).$$

Now we can use Jensen's inequality, the monotonicity of $\psi$, and the fact that $|g - \bar{g}| \leq M$ to get

$$\mathcal{E}_P(\mathcal{F}; \ell \bullet g) = P(\ell \bullet g) - P(\ell \bullet \bar{g}) \geq 2\psi(\Pi|g - \bar{g}|^r) \geq 2\psi(M^{r-2}\|g - \bar{g}\|_{L_2(\Pi)}^2),$$

which implies

$$\mathcal{F}(\delta) = \{(\ell \bullet g) : g \in \mathcal{G}, \mathcal{E}_P(\mathcal{F}; \ell \bullet g) \leq \delta\} \subset \{(\ell \bullet g) : g \in \mathcal{G}_\delta\}$$

where $\mathcal{G}_\delta := \{g \in \mathcal{G} : \|g - \bar{g}\|_{L_2(\Pi)}^2 \leq M^{2-r}\psi^{-1}(\delta/2)\}$. Therefore

$$D_P(\mathcal{F}; \delta) \leq L \sup_{g_1, g_2 \in \mathcal{G}_\delta} \|g_1 - g_2\|_{L_2(\Pi)} \leq 2LM^{1-r/2}\sqrt{\psi^{-1}(\delta/2)}.$$

We will now bound $\phi_n(\delta) = \phi_n(\mathcal{F}; \delta)$ in terms of $\theta_n(\delta) = \theta_n(\mathcal{G}; \bar{g}; \delta)$. By the symmetrization inequality,

$$\phi_n(\delta) = \mathbb{E} \sup_{f_1, f_2 \in \mathcal{F}(\delta)} |(P_n - P)(f_1 - f_2)|$$

$$\leq 2\mathbb{E} \sup_{g_1, g_2 \in \mathcal{G}(\delta)} \left|n^{-1}\sum_{i=1}^n \varepsilon_i(\ell(Y_i; g_1(X_i)) - \ell(Y_i; g_2(X_i)))\right|$$

$$\leq 4\mathbb{E} \sup_{g \in \mathcal{G}(\delta)} \left|n^{-1}\sum_{i=1}^n \varepsilon_i(\ell(Y_i; g(X_i)) - \ell(Y_i; \bar{g}(X_i)))\right|,$$

which by the contraction inequality can be bounded further by

$$16L\mathbb{E} \sup_{g \in \mathcal{G}(\delta)} \left|n^{-1}\sum_{i=1}^n \varepsilon_i(g(X_i) - \bar{g}(X_i))\right|$$

$$\leq 16L\mathbb{E} \sup\left\{\left|n^{-1}\sum_{i=1}^n \varepsilon_i(g(X_i) - \bar{g}(X_i))\right| : g \in \mathcal{G}, \|g - \bar{g}\|_{L_2(\Pi)}^2 \right.$$

$$\left. \leq M^{2-r}\psi^{-1}(\delta/2)\right\}.$$

Using now the desymmetrization inequality yields

$$\phi_n(\delta) \leq 32L\mathbb{E}\sup\{|(\Pi_n - \Pi)(g - \bar{g})| : g \in \mathcal{G}, \|g - \bar{g}\|^2_{L_2(\Pi)} \leq M^{2-r}\psi^{-1}(\delta/2)\}$$
$$+ 8L\sqrt{\frac{M^{2-r}\psi^{-1}(\delta/2)}{n}}.$$

As a result, we can bound (with a proper choice of $C$)

$$\bar{U}_n(\delta; t) \leq \bar{W}_n(\delta; t)$$
$$= C\left[L\theta_n(M^{2-r}\psi^{-1}(\delta/2)) + L\sqrt{\frac{M^{2-r}\psi^{-1}(\delta/2)(t+1)}{n}} + \frac{t}{n}\right],$$

and the first bound follows. The second bound is also immediate because of property 2, Section 2.3. □

PROOF OF THEOREM 12. We will apply the lemma with $r = 2$ and $\psi(u) = \Lambda u$. Suppose that $\theta_n$ is upper bounded by a function $\check{\theta}_n$ of strictly concave type. In this case we have

$$\bar{W}_n(\delta; t) \leq C\left[L\check{\theta}_n(\delta/(2\Lambda)) + L\sqrt{\frac{\delta(t+1)}{2\Lambda n}} + \frac{t}{n}\right].$$

Using the basic properties of the $\sharp$-transform it is easy to deduce that with some constant $C$

$$\bar{\delta}_n^W(\mathcal{G}; t) \leq C\left[2\Lambda\check{\theta}_n^\sharp\left(\frac{\Lambda}{L}\right) + \frac{L^2}{\Lambda}\frac{t+1}{n}\right].$$

Since $\mathcal{G} := M\operatorname{conv}(\mathcal{H})$, where $\mathcal{H}$ is a VC-type class of functions from $S$ into $[-1/2, 1/2]$, condition (2.1) holds for $\mathcal{H}$ with envelope $F \equiv 1$. As in Example 4 of Section 2,

$$\theta_n(\delta) \leq \check{\theta}_n(\delta) := C\left[\frac{M^\rho}{\sqrt{n}}\delta^{(1-\rho)/2} \vee \frac{M^{2\rho/(\rho+1)}}{n^{1/(1+\rho)}}\right]$$

with $\rho := \frac{V}{V+2}$. Such a $\check{\theta}_n$ is of strictly concave type and $\theta_n^\sharp(\varepsilon) \leq C\frac{M^{2\rho/(1+\rho)}}{n^{1/(1+\rho)}} \times \varepsilon^{-2/(1+\rho)}$ for $\varepsilon \leq 1$. Therefore,

$$\bar{\delta}_n^W(\mathcal{G}; t) \leq C\left[\Lambda M^{V/(V+1)}\left(\frac{L}{\Lambda} \vee 1\right)^{(V+2)/(V+1)} n^{-\frac{1}{2}\frac{V+2}{V+1}} + \frac{L^2}{\Lambda}\frac{t+1}{n}\right]$$
$$= \pi_n(M, L, \Lambda; t).$$

Assume now that for all $y$, $\ell(y, \cdot)$ is bounded by 1 on the interval $[-M/2, M/2]$. Applying Theorem 2, we get

$$\mathbb{P}\left\{P(\ell \bullet \hat{g}) \geq \min_{g \in \mathcal{G}} P(\ell \bullet g) + \pi_n(M, L, \Lambda; t)\right\} \leq e^{-t}.$$

To get rid of the assumption that $\ell$ is bounded by 1, note that if $\ell$ is bounded by $D$ on the interval $[-M/2, M/2]$, one can replace $\ell$ by $\ell/D$ and also note that $L, \Lambda$ become then $L/D, \Lambda/D$. Since $\pi_n(M, L/D, \Lambda/D; t) = \pi_n(M, L, \Lambda; t)/D$, the result follows by a simple rescaling. $\square$

**Acknowledgment.** The author is thankful to Pavel Cherepanov for pointing out a mistake in an earlier version of the paper and also to an Associate Editor and referees for a large number of helpful suggestions.

DEPARTMENT OF MATHEMATICS AND STATISTICS
UNIVERSITY OF NEW MEXICO
ALBUQUERQUE, NEW MEXICO 87131
USA
E-MAIL: vlad@math.unm.edu
AND
SCHOOL OF MATHEMATICS
GEORGIA INSTITUTE OF TECHNOLOGY
ATLANTA, GEORGIA 30332
USA
E-MAIL: vlad@math.gatech.edu